\documentstyle[12pt,amssymb,epic,eepic,amscd]{amsart}
\def\draft{n}

\theoremstyle{plain}

\newtheorem{theorem}{Theorem}
\newtheorem{proposition}{Proposition}[section]
\newtheorem{lemma}[proposition]{Lemma}
\newtheorem{corollary}[proposition]{Corollary}
\newtheorem{claim}[proposition]{Claim}

\newtheorem{add}{Addendum}
\theoremstyle{definition}
\newtheorem{definition}[proposition]{Definition}

\newtheorem{question}{Question}

\theoremstyle{remark}

\newtheorem{remark}[proposition]{Remark}

\def\printname#1{
	\if\draft y
		\smash{\makebox[0pt]{\hspace{-0.5in}
			\raisebox{8pt}{\tt\tiny #1}}}
	\fi
}

\newlength{\standardunitlength}
\catcode`\@=11
\long\def\@makecaption#1#2{%
    \vskip 10pt
    
\setbox\@tempboxa\hbox{%\ifvoid\tinybox\else\box\tinybox\fi
      \small\sf{\bfcaptionfont #1. }\ignorespaces #2}%
    \ifdim \wd\@tempboxa >\captionwidth {%
        \rightskip=\@captionmargin\leftskip=\@captionmargin
        \unhbox\@tempboxa\par}%
      \else
        \hbox to\hsize{\hfil\box\@tempboxa\hfil}%
    \fi}
\font\bfcaptionfont=cmssbx10 scaled \magstephalf
\newdimen\@captionmargin\@captionmargin=2\parindent
\newdimen\captionwidth\captionwidth=\hsize
\catcode`\@=12

\def\lbl#1{\label{#1}\printname{#1}}

\def\BZ{\Bbb Z}
\def\BQ{\Bbb Q}

\def\BC{\Bbb C}

\def\HQ{H_{\BQ}}
\def\UQ{U_{\BQ}}
\def\VQ{V_{\BQ}}
\def\A{\cal{A}}
\def\B{\cal{B}}

\def\D{\cal{D}}
\def\K{\cal{K}}
\def\J{\cal{J}}
\def\L{\cal L}
\def\G{\cal{G}}

\def\T{\cal{T}}
\def\M{\cal{M}}
\def\F{\cal{F}}

\def\R{\cal{R}}

\def\aa{\alpha}
\def\bb{\beta}

\def\l{\lambda}

\def\as{algebraically split}
     
\def\ihs{integral homology 3-sphere}          

\def\fti{finite type invariant}

\def\upa{\uparrow}

\def\Fas#1{\cal F^{as}_{#1}\cal M}
\def\FT#1{\cal F^{T}_{#1}\cal M}
\def\FK#1{\cal F^{K}_{#1}\cal M}
\def\FL#1{\cal F^{L}_{#1}\cal M}
\def\FHL#1{\cal F^{H\cal L}_{#1}\cal M}

\def\Gas#1{\cal G^{as}_{#1}\cal M}
\def\GT#1{\cal G^{T}_{#1}\cal M}
\def\GK#1{\cal G^{K}_{#1}\cal M}
\def\GL#1{\cal G^{L}_{#1}\cal M}

\def\Tg{\cal T_{g,1}}

\def\Ng{\cal N_{g,1}}

\def\GKg{\Gamma^K_g}

\def\w{\wedge}
\def\L3{\Lambda^3}
\def\lpm{L^{\pm}}
\def\lmp{L^{\mp}}

\def\tbyt#1#2#3#4{ \bigl(\begin{smallmatrix}
                   #1 & #2 \\
                   #3 & #4 
                   \end{smallmatrix} \bigr)}

\def\we{\wedge}
\def\Heg{Heegaard}

\def\PJf{\Phi^J_f}
\def\PTf{\Phi^T_f}
\def\PKf{\Phi^K_f}
\def\PLf{\Phi^L_f}

\def\Ab{\cal A_b(\phi)}
\def\Abt{\tilde{\cal A}_b(\phi)}
\def\At{\tilde{\cal A}(\phi)}

\def\f{\phi}
\def\r{\rho}
\def\a{\alpha}
\def\pa{\partial}
\def\g{\phi^{\rho}}

\def\e{\epsilon}
\def\Ga{\Gamma}

\def\sub{\subseteq}

\def\Lg{{\cal L^L_g}}
\def\Tgg{{\cal T_g}}
\def\Kgg{{\cal K_g}}
\def\Lgg{{\cal L^L_g}}

\def\LL{{\cal L}} 
\def\S{\Sigma}
\def\R{{\Bbb R}}

\def\Z{{\Bbb Z}}
\def\Q{{\Bbb Q}}

\def\sgn{\operatorname{sgn}}

\def\Yw{\cal Y_w}
\def\Yb{\cal Y_b}

\def\Af#1{\cal G_{#1}\cal A(\phi)}
\def\Ac#1{\cal G_{#1}\cal A^{conn}(\phi)}
\def\iso{\cong}
\def\gg{\gamma}
\def\I{\cal I}
\def\hfl{\Phi^L}
\def\fl{\phi^L}
\def\hfk{\phi^K}

\def\fk{\phi^K}
\def\ft{\phi^T}
\def\tg{\T_g}
\def\kg{\K_g}

\begin{document}

%%%%%%%%%%%%%%%%%%%%%%\inplude{page1}

\title[Finite type invariants  and the structure of the Torelli
group I]{Finite type 
3-manifold invariants and the structure of the Torelli group I}

\author{Stavros Garoufalidis}
\address{Department of Mathematics \\
         Harvard University \\
         Cambridge, MA 02138, U.S.A. }
\email{stavros@math.harvard.edu}
\thanks{The  authors were partially supported by NSF grant 
       DMS-95-05105  and DMS-93-03489 respectively. 
This and related preprints can also be obtained 
by
       accessing the WEB in the address 
 {\tt 
http:\linebreak[0]//\linebreak[0]www.\linebreak[0]math.\linebreak
[
0]brown.\linebreak[0]edu/\linebreak[0]$\sim$stavrosg/}}

\author{Jerome Levine}
\address{Department of Mathematics\\
        Brandeis University\\
        Waltham, MA 02254-9110, U.S.A. }
\email{levine@max.math.brandeis.edu}

\date{
This edition: September 10, 1996 \hspace{0.5cm} First edition: September 10, 
1996 
\\
Fax numbers: (617) 495 5132   \hspace{0.5cm} (617) 736 3085 
\\
  Email: {\tt stavros@math.harvard.edu} 
 \hspace{0.5cm} {\tt  levine@max.math.brandeis.edu}
}
\maketitle

\begin{abstract}
Using the recently developed theory of \fti s of \ihs s we study
the structure of the Torelli group of a closed 
surface. Explicitly, we construct (a) natural cocycles of the 
Torelli group (with coefficients in a space of trivalent graphs) 
and
cohomology classes of the abelianized Torelli group; (b)  group 
homomorphisms that detect (rationally) the nontriviality
of the lower central series of the Torelli group.
Our results are motivated by  the appearance of trivalent graphs
in topology and in representation theory and the dual role played 
by the 
Casson invariant in the theory of \fti s of \ihs s and in 
Morita's study 
\cite{Mo1}, \cite{Mo2} of the structure of the Torelli group
Our results generalize those of S. Morita \cite{Mo1}, \cite{Mo2} 
and 
complement the
recent calculation, due to R. Hain \cite{Hain}, of  the $I$-adic 
completion
of the rational group ring of the Torelli group. We also give 
analogous 
results for two other subgroups of the mapping class group.
\end{abstract}

\tableofcontents

%%%%%%%%%%%%%%%%%%%%%%%\inplude{intro}

\newpage
\section{Introduction}

\subsection{Background}
\lbl{sub.history}

The notion of \fti s for oriented \ihs s 
was introduced not long ago by T.Ohtsuki 
\cite{Oh}.  More recently T.T.Q. Le, J. Murakami and T.Ohtsuki
\cite{LMO} used the 
Kontsevich integral to give a complete classification of 
these invariants  in terms of a certain space of trivalent 
graphs.
	
In another very recent paper  \cite{GL3} the present authors gave 
several different formulations of the notion of \fti s. 
In particular we showed that one could use the 
lower central series of the Torelli group (or certain other 
subgroups of the mapping class group), in conjunction with 
Heegaard decompositions, to define them. 

It is the purpose of the present paper to exploit this 
connection, using the classification theorem of \cite{LMO}, 
to investigate the structure of the Torelli  group.
Explicitly, we:
\begin{itemize}
\item
Construct canonical  cocycles of the 
Torelli group (with coefficients in a space of trivalent graphs), 
and
cohomology classes in the abelianized Torelli group.
\item
Show, by very explicit and geometric construction, that 
the (rational) lower central series quotients of the Torelli 
group (and certain other subgroups of the mapping class 
group) map {\em onto} a  space of trivalent graphs.
\end{itemize}

In a recent paper \cite{Hain} R.Hain has given a presentation of
(the Lie algebra associated to) 
the lower central series of the Torelli group using mixed 
Hodge structures. We do not yet understand the relationship 
between our results and his but it certainly will be 
interesting to compare them.

Finally we point out that the relation between trivalent 
graphs (in the theory of \fti s) and the Torelli group has 
been foreshadowed by the work of S. Morita in:
\begin{itemize}
\item
The appearance of  trivalent graphs in invariant theory applied
to the Torelli group,  see \cite{Mo5}, \cite{KM}.
\item
The study of the   {\em Casson 
invariant} in terms of the Torelli group and other subgroups of
the mapping class group, see  \cite{Mo1}, \cite{Mo2}.
\end{itemize}

\subsection{Trivalent graphs in topology and in representation 
theory}
\lbl{sub.trigraphs}

We begin by recalling the appearance of trivalent graphs in 
topology (in the theory of \fti s of \ihs s) and in 
representation
theory (related to invariant tensors of the abelianization of the 
Torelli
group). 
{\em Finite type invariants} of \ihs s were introduced  by T. 
Ohtsuki 
\cite{Oh}, in terms of
a decreasing filtration $\Fas \ast$ on the vector space $\M$ 
(over 
$\BQ$)
 of isomorphism classes of oriented, connected
 \ihs s. A linear map $v:\M \to \BQ$
is called a {\em type} $m$ {\em invariant} of \ihs s if 
$v(\Fas 
{m+1})=0$.
The associated graded quotients  $\Gas \ast$ of the
filtration $\Fas \ast$ has recently been related to trivalent 
graphs in the following way. A graph is called {\em admissible}
if it is  trivalent and  vertex oriented. Let $\A(\phi)$ denote 
the 
$\BQ$ vector space on the set of admissible graphs,
modulo  two relations: the $AS$ and the $ IHX$; see figure \ref{ASIHX}
and \cite{GO1},
\cite{LMO}.
$\A(\phi)$ has a natural grading; the degree of a trivalent 
graph 
is
half the number of its vertices, thus a degree $m$ trivalent 
graph 
has
$3m$ edges and $2m$ vertices.    

\begin{figure}[htpb]
$$ \printname{ASIHX}
	\setlength{\unitlength}{0.03\standardunitlength}
	\begin{array}{c}  \hspace{-1.7mm}
        	\raisebox{-8pt}{\begingroup\makeatletter\ifx\SetFigFont\undefined
% extract first six characters in \fmtname
\def\x#1#2#3#4#5#6#7\relax{\def\x{#1#2#3#4#5#6}}%
\expandafter\x\fmtname xxxxxx\relax \def\y{splain}%
\ifx\x\y   % LaTeX or SliTeX?
\gdef\SetFigFont#1#2#3{%
  \ifnum #1<17\tiny\else \ifnum #1<20\small\else
  \ifnum #1<24\normalsize\else \ifnum #1<29\large\else
  \ifnum #1<34\Large\else \ifnum #1<41\LARGE\else
     \huge\fi\fi\fi\fi\fi\fi
  \csname #3\endcsname}%
\else
\gdef\SetFigFont#1#2#3{\begingroup
  \count@#1\relax \ifnum 25<\count@\count@25\fi
  \def\x{\endgroup\@setsize\SetFigFont{#2pt}}%
  \expandafter\x
    \csname \romannumeral\the\count@ pt\expandafter\endcsname
    \csname @\romannumeral\the\count@ pt\endcsname
  \csname #3\endcsname}%
\fi
\fi\endgroup
\begin{picture}(9396,939)(0,-10)
\thicklines
\path(12,912)(912,912)
\path(462,912)(462,12)
\path(12,12)(912,12)
\path(1812,912)(1812,12)
\path(2712,912)(2712,12)
\path(1812,462)(2712,462)
\path(3612,912)(4512,12)
\path(3612,12)(3987,387)
\path(4137,537)(4512,912)
\path(3837,237)(4287,237)
\path(6462,462)(6462,12)
\path(7812,912)(8262,462)(8712,912)
\path(8262,462)(8262,12)
\path(6012,912)	(6057.918,898.593)
	(6097.434,886.777)
	(6159.896,867.041)
	(6204.660,851.035)
	(6237.000,837.000)

\path(6237,837)	(6289.803,808.922)
	(6354.169,770.355)
	(6416.200,727.611)
	(6462.000,687.000)

\path(6462,687)	(6513.090,620.977)
	(6532.883,579.245)
	(6537.000,537.000)

\path(6537,537)	(6507.300,490.586)
	(6462.000,462.000)

\path(6462,462)	(6424.500,456.195)
	(6387.000,462.000)

\path(6387,462)	(6341.700,490.586)
	(6312.000,537.000)

\path(6312,537)	(6323.396,600.690)
	(6347.793,639.260)
	(6387.000,687.000)

\path(6462,762)	(6517.919,791.245)
	(6559.215,812.321)
	(6612.000,837.000)

\path(6612,837)	(6660.610,856.239)
	(6722.711,878.955)
	(6785.707,899.443)
	(6837.000,912.000)

\path(6837,912)	(6863.186,913.590)
	(6912.000,912.000)

\put(1137,387){\makebox(0,0)[lb]{$=$}}
\put(3162,387){\makebox(0,0)[lb]{$-$}}
\put(7212,387){\makebox(0,0)[lb]{$+$}}
\put(9012,387){\makebox(0,0)[lb]{$=0$}}
\end{picture} }
        	\hspace{-1.9mm}
	\end{array}
 $$
\caption{The $IHX$ and the $AS$ relations on $\cal A(\phi)$.}\lbl{ASIHX}
\end{figure}

One can define a map (for details see section \ref{sub.fti}):
\begin{equation}
\lbl{eq.GA}
\G_{m}  \A(\phi) \to \Gas {3m}
\end{equation}
which was shown in \cite{GO1} to be well defined and onto.
According to the {\em fundamental theorem} of \fti s of \ihs s 
\cite{LMO}, \cite{L} the map \eqref{eq.GA} is one-to-one and
thus a vector space isomorphism.
We wish to think of the above isomorphism as a  relation 
between
 \fti s of \ihs s and {\em trivalent graphs} (decorated by a 
choice
of a vertex orientation, and considered modulo the $AS$ and 
$IHX$
relations).

As it turns out, one can reformulate \cite{GL3} the notion of 
\fti s of
\ihs s in such a way that  makes explicit the dependence  of the 
values
of  \fti s on 
 manifolds obtained by cutting, twisting and gluing of higher 
genus
surfaces. This reformulation, given in \cite{GL3} in terms of six 
filtrations on $\M$, will be used crucially on the present paper.
Three  filtrations on $\M$ were defined in \cite{GL3}
using surgery on special classes of links, and three more 
filtrations by 
using cutting, twisting and gluing along embedded surfaces. In 
the present 
paper we will mainly 
concentrate on one of them, denoted by $\FT \ast$ in 
\cite{GL3} (note however that the proofs of our results will 
require
the use of a few other filtrations from \cite{GL3}).
Following the notation as in \cite{GL3}, we briefly recall  
the definition of $\FT \ast$ .
Let $M$ be an \ihs\ and $f :\Sigma \hookrightarrow M$ an 
embedded,
oriented, connected, closed  genus $g$ surface in $M$. Such a 
surface will be
called {\em admissible} in $M$. Note that an admissible 
surface 
has 
no boundary. Since $M$ is  an \ihs\ it follows that  an 
admissible
surface is {\em separating}, i.e., $M - f(\S)$ is the union 
of two
connected components 
$M_+^o$ and 
$M_-^o$, where the positive 
normal vector to $f(\S)$ points into $M_+^o$. Let $M_\e$  
(for  $\e =\pm$) denote their closures; they are compact 
$3$-manifolds 
with boundary $f(\S)$. There is a natural decomposition 
$H_1 (f(\S) )=L_+ \oplus L_-$, where 
$L_{\e}=\text{Ker}\{ (i_{\e})_{\ast}:H_1 (f(\S ))\to H_1 
(M_{\e})\}$ and $i_{\e}:f(\S )\to M_{\e}$ is the inclusion. 
We refer to $(L_+ ,L_-)$ as the {\em Lagrangian pair} of the 
symplectic space
$H= H_1 (f(\S ))$  associated to 
the admissible surface $f :\Sigma \hookrightarrow M$. 
If $h\in\Ga (f(\S ))$ (the mapping class group of $f(\S)$, 
i.e., 
the
group of isotopy classes of orientation preserving 
diffeomorphisms 
of
$f(\S)$) let $M_h$ denote $M_+ \amalg M_-$ with the 
identifications: 
$i_+ 
(x)\leftrightarrow i_- h(x)$ for every $x\in f(\S)$. The 
notation 
$M_h$
does not explicitly indicate the dependence of $M_h$ on the 
admissible surface
which we keep fixed; we hope that this will not confuse the 
reader.
Note that if $h \in \T (f(\S ))$, the {\em Torelli group} of 
$f(\Sigma )$
(i.e., all elements of the mapping class group that act 
trivially 
on  the 
homology of the surface) and $M$ is
an \ihs\ then the resulting manifold $M_h$ will also 
be an \ihs\ . For a closed surface $\S$ of genus $g$, 
let $\Tgg =\T (\S )$. The assignment $h\to M_{fhf^{-1}}$ 
defines a 
map 
\begin{equation}
\lbl{eq.fun}
\Phi^T_f: \BQ\Tgg\to\M
\end{equation}
 where $\BQ\Tgg$ is  the rational
group ring of $\Tgg$. Let $I \Tgg$ denote  the augmentation 
ideal 
(i.e., the two sided ideal of $\BQ \Tgg$
generated by elements of the form $1-f$, for $f \in \Tgg$). 
We 
{\em define}
 $\FT {m}$ to be the union in $\M$ of the image $\Phi^T_f((I 
\Tgg)^m)$
for all admissible surfaces $f:\S \hookrightarrow M$ in all 
\ihs s 
$M$. Alternatively we can choose a single {\em Heegaard} 
embedding (i.e. $M_+$ and $M_-$ are handlebodies) into $S^3$ 
for 
each genus $g$ and let $\FT {m}$ be the union in $\M$ of the 
images 
$\Phi^T_f((I \Tgg)^m)$ for these embeddings. It is shown in 
\cite{GL3} 
that this gives the same filtration, 
%that $\FT \ast$ is a 
%decreasing 
%filtration on $\M$,
and that $\Fas {3m} = \FT {2m}= \FT {2m-1}$.

Let us now recall one more ingredient, the 
Johnson homomorphism \cite{Jo1}, related to the {\em 
abelianization}
of the Torelli group.   For a more 
detailed description, as well as a summary of properties of 
the 
Johnson 
homomorphism,
see section \ref{sub.johnson}. If $\Sigma$ is a closed
surface, D. Johnson defined a homomorphism
 $ \tau : \Tgg\to U=\L3 H/H$ where $H=H_1(\Sigma, \BZ)$.
There are several versions of Johnson's homomorphism, 
depending
on the surface being closed, or punctured, or with boundary 
components.
The following are three
 important properties of Johnson's homomorphism (and its 
various versions):
\begin{itemize}
\item
It  coincides (modulo torsion) with the abelianization
of the Torelli group.
\item
It is equivariant with respect to the action of the mapping 
class 
group
of the surface.
\item
It is stable with respect to an inclusion of a surface with 
one
boundary component into another.
\end{itemize} 
These properties  have been used extensively by R. Hain 
\cite{Hain} and
S. Morita \cite{Mo1}, \cite{Mo2} to study questions relating 
to the lower central series of the Torelli group and
the mapping class group. From its very definition, the 
image of $\tau$ 
is a quotient $U$ of the third exterior power of $H$.
Furthermore, it turns out that the invariant space of 
$\otimes^{2m} U$
under the symplectic or the general linear group can be 
described 
in
terms of suitably ``decorated'' {\em trivalent graphs}, 
modulo an 
$AS$
relation. For a precise statement, see section \ref{sub.gsp} 
and
especially definition \ref{def.cl}.

It is a natural question to ask whether the above two 
appearances 
of trivalent graphs in the theory of \fti s of \ihs s and 
in the abelianization of the  Torelli group
are related to each other. For a positive
 answer (in terms of a stably onto map
$\Psi_{\lpm,m}$) see theorem \ref{thm.formula}. 

\subsection{The role of the Casson invariant}
\lbl{sub.cassonrole}

As was stated above, a main motivation for the present work was 
the
role of the Casson invariant in the theory of \fti s and in the 
work
of S. Morita \cite{Mo2}.
Explicitly, the Casson invariant $\l$ \cite{AM} has the following 
properties:
\begin{itemize}
\item
$\l$ is a type $3$ invariant of \ihs s, \cite{Oh}. 
\item
Given an admissible genus $g$ surface $f:\S \hookrightarrow M$, 
Morita \cite{Mo2}
used the Casson invariant to construct a map:
\begin{equation}
2 \delta_f: U \otimes U \to \BQ
\end{equation}
where the notation is as in \cite{Mo2}.
\item
Furthermore, given an admissible surface Morita \cite{Mo1}
used the Casson invariant to construct a group homomorphism:
\begin{equation}
\K_g \to \BQ
\end{equation}
where $\K_g$ is the kernel of the Johnson homomorphism $\Tgg \to 
U$.
\end{itemize}

It is a natural question to ask whether one can use \fti s 
of \ihs s to
generalize the two maps constructed above. For a positive answer, 
see
theorems \ref{thm.0}, \ref{thm.00}, \ref{thm.homo}  and 
especially
corollary \ref{cor.type} and theorem \ref{thm.mor}.

\subsection{Statement of the results}
\lbl{sub.results}

In this section we state the main results of the paper. 
For a group $G$, and a positive
integer $n$, let us  inductively define
the {\em lower central series} subgroups of $G$ by 
$G_{n+1}=[G,G_n]$, 
with 
$G_1=G$. Let us also define $G(n)$ to consist of all elements 
of $G$ for which a nonzero power lies in $G_n$. We call $G(n)$
the $n^{th}$ term in the {\em rational lower central series}
of $G$.   
With the notation as in section \ref{sub.history}, we have 
the following:

\begin{theorem}
\lbl{thm.0} 
Let $f :\Sigma \hookrightarrow M$ be an admissible surface.
For every non-negative integer $m$, there is a map 
\begin{equation}
\lbl{eq.map}
C_{f,m}:
\otimes^{2m}U  \to \G_m\A(\phi)
\end{equation}
with the following properties:
\begin{itemize}
\item
The map $C_{f,m}$ is multilinear and
$Sp(H)$ equivariant, i.e.,  satisfies the following
property (for $\alpha_i \in U, h \in \Ga_g$):
\begin{equation}
C_{f,m}(h_\ast \alpha_1 ,  \ldots h_\ast \alpha_{2m} )= 
C_{h f,m}(\alpha_1, \ldots \alpha_{2m})
\end{equation}
\item
$C_{f,m}$ is a $2m$ cocycle of the abelian group $U$ with 
coefficients in the trivial $U$ module 
$\G_m\A(\phi)$. In particular, it represents a cohomology
class $[C_{f,m}] \in H^{2m}(U, \G_m\A(\phi))$
\item
The pullback of the cocycle $C_{f,m}$ to $\Tgg$ and, in fact, 
to 
$\Tgg/\Tgg (3)$,
under the projection maps $\Tgg \to \Tgg/ \Tgg(3)  \to U$,
is a coboundary.
\end{itemize}
\end{theorem}

\begin{add}
\lbl{cor.cor}
With the above notation, given  an admissible surface 
$f: \S \hookrightarrow M$,
 the following diagram commutes:
$$
\begin{CD}
 \Tgg^{2m} @>>> \FT {2m} \\
 @V\tau^{\otimes 2m}VV           @VVV     \\
\otimes^{2m} U @>C_{f,m}>> \G_{m} \A(\phi)
\end{CD}
$$
where the right 
vertical map is the composition of maps (defined in section 
\ref{sub.history}):
$$ \FT {2m} = \Fas {3m} \to \Gas {3m} \simeq \G_m\A(\phi)$$
and the top horizontal map is the map: $(h_1, \ldots, h_{2m}) 
\to 
\Phi^T_f((1-h_1) \ldots (1-h_{2m}))$.
\end{add}

\begin{corollary}
\lbl{cor.type}
If $v$ is a type $3m$ invariant of \ihs s, then 
its associated weight system is an element  $W_v \in  
\G_m\A^\ast(\phi)$, 
\cite{GO1}. Given an admissible surface $f:\S \hookrightarrow 
M$,
we thus get a $2m$ cocycle $ W_v \circ C_{f,m}$ of $U$ with
rational coefficients.
\end{corollary}

\begin{theorem}
\lbl{thm.00}
Let $f :\Sigma \hookrightarrow M$ be an admissible {\em 
Heegaard} 
surface.
Then, for every non-negative integer $m$,
\begin{itemize}
\item
The cocycle $C_{f,m}$ depends only on the associated 
Lagrangian
pair $(L^+, L^-)$ of the admissible surface, and  will thus 
be 
denoted by $C_{\lpm,m}$.
\item Using the natural onto maps $U \to \L3 
(H/L^{\mp})\simeq\L3 L^{\pm}$, 
 $C_{\lpm,m}: \otimes^{2m} U \to \G_m\A(\phi)$ factors though 
a 
$GL(L^+_{\BQ})$-invariant map:
\begin{equation}
\lbl{eq.fct}
\L3 L^+_{\BQ} \otimes ( \otimes^{2m-2} \UQ) \otimes \L3 L^-
_{\BQ}
\to \G_m\A(\phi)
\end{equation}
\item
If we change the orientation of the \ihs\ $M$, this results 
in a 
permutation
of the Lagrangian pair $L^\pm$ and the associated cocycle 
satisfies:
\begin{equation}
C_{\lpm,m}(g_1, g_2, \ldots, g_{2m})= 
(-1)^m C_{L^{\mp},m}(g_{2m}, g_{2m-1}, \ldots, 
g_1)
\end{equation} 
Passing to cohomology classes though, we have
\begin{equation}
[C_{\lpm,m}]=[C_{L^{\mp},m}]
\end{equation}
\end{itemize}
\end{theorem}

\begin{add}
\lbl{cor.corr}
In the case of an admissible Heegaard surface of genus $g$, 
the map $C_{\lpm,m}$
is stable  with respect to an inclusion of one (punctured)  
admissible Heegaard surface
into another. Furthermore, the map $C_{\lpm,m}$ is
stably   (i.e., for $g >> m$)  onto.
In particular, the cocycles $C_{\lpm,m}$ of $U$ are stably 
nontrivial.
\end{add}
See  section~\ref{sub.prf} for a more precise assertion of 
Addendum 
\ref{cor.corr}.

Before we state the next theorem we need to define 
$\G_m\A^{rp,nl,cl}$, a vector space over $\BQ$ on the
set  of ``decorated''  trivalent graphs with $3m$ edges, 
modulo a 
colored antisymmetry relation (see figure \ref{AScl}). These 
 decorations 
involve a choice of {\em ordering for the vertices}, as well
as a choice of {\em vertex orientation} and a choice of $2$-
{\em 
coloring}
for the edges. For a precise definition, as well as 
motivation 
coming from
representation theory of classical Lie groups, see definition 
\ref{def.cl}
and section \ref{sub.gsp}.

\begin{theorem}
\lbl{thm.formula}
Let $f :\Sigma \hookrightarrow M$ be an admissible {\em 
Heegaard} 
surface.
Then, for every non-negative integer $m$ there is an onto 
map:
$\G_m\A^{rp,nl,cl} \to \otimes^{2m} U$, which combined 
with the map $C_{\lpm,m}:\otimes^{2m} U\to \G_m\A(\phi)$, 
induce a (stably 
onto) 
map:
\begin{equation}
\lbl{eq.stonto}
\Psi_{\lpm,m}: \G_m\A^{rp,nl,cl} \to \G_m\A(\phi)
\end{equation}
\end{theorem}

The above map compares trivalent graphs related to the 
abelianization
of the Torelli group (on the left) to trivalent graphs 
related to
\fti s of \ihs s (on the right), and fulfills one of the 
goals of 
the 
present paper. 

In the case of $m=1$, we have an explicit description of the
cocycle $C_{\lpm,1}$ of theorem \ref{thm.00} and corollary 
\ref{cor.type}
and of the map $\Psi_{\lpm,1}$ of theorem \ref{thm.formula}.
We need to recall first that the Casson
invariant $\l$ \cite{AM} is a type $3$ invariant, \cite{Oh}.
 Let $W_\l \in 
\G_1\A(\phi)^\ast$ denote its associated manifold weight 
system as in \cite{GO1}.
Let $\Theta_w \in \G_1\A(\phi)$ denote the trivalent graph  
$\Theta$
with a fixed choice of vertex orientation. 
Let $f:\S\hookrightarrow M$ be an admissible Heegaard
surface, and let 
$C_\Theta: \otimes^6 H \to \BQ$ be given 
%%%%\L3 L^+_{\BQ} \otimes \L3 L^-_{\BQ} 
by:
\begin{equation}
\lbl{eq.ctheta}
C_\Theta(a_1 \otimes a_2 \otimes a_3, b_1 \otimes b_2 \otimes 
b_3) = 
\omega(a_1, b_1)
\omega(a_2, b_2) \omega(a_3, b_3)
\end{equation}
 for $a_i, b_i  \in H$, where $\omega$ is the intersection
pairing on $H$.
%%%L^+$, $b_i \in L^-$.
Recall  the onto maps $\UQ \to \L3 L^{\pm}_{\BQ}$ from theorem 
\ref{thm.00},
their tensor product  $\otimes^2 \UQ \to \L3 L^+_{\BQ} \otimes 
\L3 L^-_{\BQ}$ and the natural inclusion 
 $\L3 L^+_{\BQ} \otimes \L3 L^-_{\BQ} \hookrightarrow \otimes^2 
\L3 \HQ
\hookrightarrow \otimes^6 \HQ$.
Let us denote by $C^U_{\Theta}$ the pullback of $C_{\Theta}$ to
$\otimes^2 \UQ$ under the composition of the above maps. 
Then, we have the following theorem:

\begin{theorem}
\lbl{thm.mor}
Given an admissible Heegaard surface  $f: \S \hookrightarrow M$, 
\begin{itemize}
\item
The map 
$(W_\l) \circ C_{\lpm, 1}: \otimes^2 \UQ \to \BQ $ is given by:
\begin{equation}
(W_\l) \circ C_{\lpm, 1} =
2 C^U_{\Theta} 
%\sum_{i < j < k} \alpha^1_{ijk} \alpha^2_{ijk} 
\end{equation}
\item
the map $C_{\lpm,1}: \otimes^2 U 
\to \G_1\A(\phi)$ is given as follows. For $\alpha_1, \alpha_2 
\in U$ we have:
\begin{equation}
\lbl{eq.mor}
C_{\lpm,1}(\alpha_1, \alpha_2)= - C^U_{\Theta} (\alpha_1, 
\alpha_2) \cdot 
%(\sum_{i < j < k} \alpha^1_{ijk} \alpha^2_{ijk})
 \Theta_w
\end{equation}
\item
The vector space $\G_1\A^{rp,nl,cl}$ is four dimensional, 
with a basis
given in the southeast part of  figure \ref{Basis1}. 
The map \eqref{eq.stonto} of theorem \ref{thm.formula} is 
given as
follows: 
%%$$\eepic{Psi5}{0.03}$$
\item The cocycle $C_{\lpm,1}$ defines a nonzero cohomology 
class $[C_{\lpm,m}] \in H^2 (U;\Q )$, if $\dim H\geq 6$.  Moreover 
$[C_{\lpm,m}]$ 
depends on the 
Lagrangian pair $\lpm$ in the following very strong sense. If 
$K^{\pm}$ is another Lagrangian pair then $[C_{\lpm,m}] 
=[C_{K^{\pm},m}]$
if and only if one of the following holds:
\begin{description}
\item[$\bullet$] $\dim H<6$
\item[$\bullet$] $L^+ =K^+ ,L^- =K^-$, or
\item[$\bullet$] $L^+ =K^- ,L^- =K^+$.
\end{description}
\end{itemize} 
\end{theorem}

\begin{remark}
\lbl{rem.morr}
In coordinates, the map $C^U_{\Theta}$ is given as follows.
Let $\{ x_i \}_{i=1}^{g}$ (respectively, $\{ y_i 
\}_{i=1}^{g}$)
 be basis for $L^+$ (respectively, $L^-$) such that 
$\omega(x_i, y_j)= \delta_{i,j}$. Using the natural projection 
$\L3 \HQ 
\to \UQ$,  consider $\alpha_1, \alpha_2 \in \UQ$ and let
$\bar{\alpha}_1, \bar{\alpha}_2 \in \L3 \HQ$ be their lifts  
written as:
\begin{align*}
\bar{\alpha}_1 &= \sum_{i < j < k} \alpha^1_{ijk} x_i \we x_j \we 
x_k\ +\text{  other terms}\\
  \bar{\alpha}_2 &= \sum_{i < j < k} \alpha^2_{ijk} y_i \we y_j 
\we 
y_k\ +\text{  other terms}
\end{align*}
Then, we have:
\begin{equation}
C^U_{\Theta}(\bar{\alpha}_1, \bar{\alpha}_2)=
\sum_{i < j < k} \alpha^1_{ijk} \alpha^2_{ijk}
\end{equation}
\end{remark}

\begin{remark}
\lbl{rem.morita}
Note that the above map $(W_\l) \circ C_{\lpm, 1}$ coincides
with the map $2 \delta_f$ of \cite[definition 4.1, theorem 
4.3]{Mo2},
and that the first part of the above theorem  was originally 
proven
by Morita \cite[theorem 4.3]{Mo2}.
Morita's result was a starting point for the results of the
present paper. It is interesting to note that the factor of 
$2$
in $2 \delta_f $ in the above mentioned 
paper of Morita was derived from a representation theory 
calculation
(counting irreducible components of $Sp(H)$ representations),
whereas in our context it comes from the identity $\Theta_w= 
2 
Y_w$. 
\end{remark}

The maps $C_{f,m}$ assemble well to define a map:
$C_f :  T_{ev}(U) \to \A(\phi)$, where $T_{ev}(U)=
\oplus_{m=0}^{\infty} (\otimes^{2m} U)$. Recall that
$T_{ev}(U)$ and $\A(\phi)$ are graded Hopf algebras,
where the comultiplication in $T_{ev}(U)$ is given by declaring
$U$ to be the set of primitive elements. 
$C_{f}$ respects the multiplication in the
following sense. For $i=1,2$ let $f_i: \S_{g_i} \to M$  be 
two
 admissible genus $g_i$ 
surfaces  disjointly embedded in an \ihs\ $M$. Without 
loss of generality, let us assume
that $f_i$ are inclusions. Assume  that there is an embedded 
$2$-sphere $S \hookrightarrow M$ with the following 
properties:
\begin{itemize}
\item
The intersections $S \cap \S_{g_i} = D_i$ for $i=1,2$ are
disjoint discs. Let $\S_i= \S_{g_i} - \text{Int} D_i$.
\item
Recall that $S$ separates $M-S$ in two components. Assume 
that
$\S_i$ for $i=1,2$ lie in different components of $M - S$.
\end{itemize}
Then we can form the composite admissible surface
$f_1 \cup f_2:\S= \S_1 \cup_{ S -(D_1 \cup D_2)} \S_2 
\to M$. Considering the homology of $\S_{g_1}, \S_{g_2}$ and $\S$
we get  natural onto maps 
$  U_{g_1+g_2} \to U_{g_i}$ for $i=1,2$ which in turn induce onto 
maps:
$ T_{ev}( U_{g_1+g_2}) \to T_{ev}(U_{g_i})$.
Let $C_{f_i}^i$ denote the pullbacks of the maps $C_{f_i}$ to 
$ T_{ev}( U_{g_1+g_2})$ for $i=1,2$.

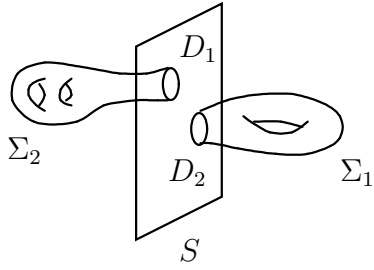
\begin{figure}[htpb]
$$ \printname{gluesurface}
	\setlength{\unitlength}{0.03\standardunitlength}
	\begin{array}{c}  \hspace{-1.7mm}
        	\raisebox{-8pt}{\begingroup\makeatletter\ifx\SetFigFont\undefined
% extract first six characters in \fmtname
\def\x#1#2#3#4#5#6#7\relax{\def\x{#1#2#3#4#5#6}}%
\expandafter\x\fmtname xxxxxx\relax \def\y{splain}%
\ifx\x\y   % LaTeX or SliTeX?
\gdef\SetFigFont#1#2#3{%
  \ifnum #1<17\tiny\else \ifnum #1<20\small\else
  \ifnum #1<24\normalsize\else \ifnum #1<29\large\else
  \ifnum #1<34\Large\else \ifnum #1<41\LARGE\else
     \huge\fi\fi\fi\fi\fi\fi
  \csname #3\endcsname}%
\else
\gdef\SetFigFont#1#2#3{\begingroup
  \count@#1\relax \ifnum 25<\count@\count@25\fi
  \def\x{\endgroup\@setsize\SetFigFont{#2pt}}%
  \expandafter\x
    \csname \romannumeral\the\count@ pt\expandafter\endcsname
    \csname @\romannumeral\the\count@ pt\endcsname
  \csname #3\endcsname}%
\fi
\fi\endgroup
\begin{picture}(5542,3666)(0,-10)
\thicklines
\put(2287,2514){\ellipse{224}{450}}
\put(2682,1888){\ellipse{224}{450}}
\path(3000,3639)(1800,3039)(1800,339)
	(3000,939)(3000,3639)
\path(2700,2139)	(2774.937,2171.229)
	(2839.660,2198.727)
	(2895.266,2221.935)
	(2942.854,2241.293)
	(2983.523,2257.238)
	(3018.371,2270.212)
	(3075.000,2289.000)

\path(3075,2289)	(3111.733,2299.067)
	(3156.845,2310.077)
	(3207.594,2321.461)
	(3261.240,2332.646)
	(3315.044,2343.064)
	(3366.265,2352.142)
	(3412.164,2359.311)
	(3450.000,2364.000)

\path(3450,2364)	(3509.485,2369.284)
	(3581.893,2374.451)
	(3621.617,2376.773)
	(3662.980,2378.804)
	(3705.451,2380.458)
	(3748.500,2381.648)
	(3791.596,2382.285)
	(3834.207,2382.285)
	(3875.805,2381.559)
	(3915.857,2380.020)
	(3989.203,2374.157)
	(4050.000,2364.000)

\path(4050,2364)	(4116.759,2347.632)
	(4154.802,2337.590)
	(4195.178,2326.228)
	(4237.306,2313.481)
	(4280.605,2299.287)
	(4324.492,2283.584)
	(4368.386,2266.309)
	(4411.706,2247.399)
	(4453.869,2226.792)
	(4494.295,2204.425)
	(4532.401,2180.235)
	(4599.329,2126.138)
	(4650.000,2064.000)

\path(4650,2064)	(4677.499,1994.878)
	(4684.373,1955.174)
	(4686.664,1914.000)
	(4684.371,1872.826)
	(4677.496,1833.122)
	(4650.000,1764.000)

\path(4650,1764)	(4611.997,1717.397)
	(4561.801,1676.824)
	(4502.902,1641.906)
	(4438.789,1612.267)
	(4372.952,1587.533)
	(4308.882,1567.328)
	(4250.068,1551.275)
	(4200.000,1539.000)

\path(4200,1539)	(4154.387,1531.502)
	(4099.343,1527.396)
	(4038.065,1526.055)
	(3973.747,1526.858)
	(3909.587,1529.177)
	(3848.778,1532.391)
	(3794.517,1535.873)
	(3750.000,1539.000)

\path(3750,1539)	(3697.521,1543.860)
	(3633.582,1551.484)
	(3562.034,1561.134)
	(3524.612,1566.489)
	(3486.731,1572.075)
	(3448.875,1577.799)
	(3411.524,1583.569)
	(3340.265,1594.879)
	(3276.806,1605.268)
	(3225.000,1614.000)

\path(3225,1614)	(3187.807,1621.092)
	(3142.480,1630.707)
	(3091.703,1641.908)
	(3038.160,1653.761)
	(2984.534,1665.332)
	(2933.510,1675.685)
	(2887.770,1683.886)
	(2850.000,1689.000)

\path(2850,1689)	(2797.838,1690.920)
	(2756.386,1690.440)
	(2700.000,1689.000)

\path(3300,2064)	(3371.164,2012.999)
	(3433.191,1969.853)
	(3487.180,1933.902)
	(3534.229,1904.486)
	(3575.437,1880.948)
	(3611.902,1862.626)
	(3675.000,1839.000)

\path(3675,1839)	(3741.942,1827.556)
	(3782.365,1824.517)
	(3824.790,1823.362)
	(3867.241,1824.162)
	(3907.743,1826.988)
	(3975.000,1839.000)

\path(3975,1839)	(4013.876,1854.931)
	(4060.819,1882.898)
	(4121.102,1926.415)
	(4157.895,1955.105)
	(4200.000,1989.000)

\path(3450,1989)	(3489.180,1990.401)
	(3525.620,1991.615)
	(3590.830,1993.483)
	(3646.728,1994.604)
	(3694.414,1994.977)
	(3734.985,1994.604)
	(3769.541,1993.483)
	(3825.000,1989.000)

\path(3825,1989)	(3870.266,1981.496)
	(3930.667,1967.666)
	(3968.741,1957.830)
	(4013.235,1945.754)
	(4065.028,1931.217)
	(4125.000,1914.000)

\path(2250,2739)	(2190.028,2721.783)
	(2138.235,2707.246)
	(2093.741,2695.170)
	(2055.667,2685.334)
	(1995.266,2671.504)
	(1950.000,2664.000)

\path(1950,2664)	(1912.624,2660.657)
	(1867.061,2658.517)
	(1816.048,2657.457)
	(1762.320,2657.351)
	(1708.615,2658.077)
	(1657.669,2659.510)
	(1612.218,2661.526)
	(1575.000,2664.000)

\path(1575,2664)	(1529.901,2669.088)
	(1475.100,2677.254)
	(1413.854,2687.572)
	(1349.423,2699.115)
	(1285.063,2710.956)
	(1224.034,2722.168)
	(1169.593,2731.825)
	(1125.000,2739.000)

\path(1125,2739)	(1066.095,2748.707)
	(994.398,2761.769)
	(955.032,2769.019)
	(914.004,2776.461)
	(871.827,2783.879)
	(829.012,2791.057)
	(786.072,2797.781)
	(743.519,2803.835)
	(701.864,2809.002)
	(661.621,2813.067)
	(587.416,2817.031)
	(525.000,2814.000)

\path(525,2814)	(483.798,2807.670)
	(435.647,2798.072)
	(383.212,2785.066)
	(329.156,2768.512)
	(276.143,2748.270)
	(226.837,2724.197)
	(183.901,2696.154)
	(150.000,2664.000)

\path(150,2664)	(124.788,2627.449)
	(103.363,2583.989)
	(86.092,2535.601)
	(73.346,2484.262)
	(65.495,2431.954)
	(62.907,2380.655)
	(65.952,2332.344)
	(75.000,2289.000)

\path(75,2289)	(93.894,2243.759)
	(122.203,2198.108)
	(157.876,2153.358)
	(198.866,2110.819)
	(243.123,2071.803)
	(288.597,2037.620)
	(333.239,2009.582)
	(375.000,1989.000)

\path(375,1989)	(437.825,1969.595)
	(511.669,1957.297)
	(551.545,1953.570)
	(592.761,1951.327)
	(634.848,1950.469)
	(677.332,1950.900)
	(719.744,1952.522)
	(761.612,1955.236)
	(802.464,1958.946)
	(841.830,1963.553)
	(914.216,1975.068)
	(975.000,1989.000)

\path(975,1989)	(1017.096,2005.500)
	(1063.942,2031.311)
	(1113.845,2063.440)
	(1165.114,2098.890)
	(1216.056,2134.666)
	(1264.979,2167.774)
	(1310.191,2195.217)
	(1350.000,2214.000)

\path(1350,2214)	(1394.240,2227.615)
	(1448.515,2240.295)
	(1509.522,2251.901)
	(1573.961,2262.296)
	(1638.530,2271.342)
	(1699.927,2278.900)
	(1754.851,2284.832)
	(1800.000,2289.000)

\path(1800,2289)	(1866.235,2291.396)
	(1907.080,2291.496)
	(1950.124,2291.130)
	(1993.152,2290.498)
	(2033.951,2289.799)
	(2100.000,2289.000)

\path(2100,2289)	(2152.035,2289.000)
	(2193.518,2289.000)
	(2250.000,2289.000)

\path(525,2589)	(471.128,2575.048)
	(425.931,2560.604)
	(358.924,2528.483)
	(318.705,2489.120)
	(300.000,2439.000)

\path(300,2439)	(301.329,2364.827)
	(335.756,2294.790)
	(408.555,2221.858)
	(460.992,2182.109)
	(525.000,2139.000)

\path(450,2514)	(485.449,2463.329)
	(508.594,2423.783)
	(525.000,2364.000)

\path(525,2364)	(508.594,2304.217)
	(485.449,2264.671)
	(450.000,2214.000)

\path(900,2589)	(838.338,2540.399)
	(795.128,2501.543)
	(750.000,2439.000)

\path(750,2439)	(735.030,2364.818)
	(737.513,2323.867)
	(750.000,2289.000)

\path(750,2289)	(798.188,2247.724)
	(840.633,2230.797)
	(900.000,2214.000)

\path(825,2514)	(864.207,2466.260)
	(888.604,2427.690)
	(900.000,2364.000)

\path(900,2364)	(880.421,2328.064)
	(825.000,2289.000)

\put(0,1464){\makebox(0,0)[lb]{$\Sigma_2$}}
\put(2400,39){\makebox(0,0)[lb]{$S$}}
\put(4650,1164){\makebox(0,0)[lb]{$\Sigma_1$}}
\put(2400,2889){\makebox(0,0)[lb]{$D_1$}}
\put(2250,1164){\makebox(0,0)[lb]{$D_2$}}
\end{picture} }
        	\hspace{-1.9mm}
	\end{array}
 $$
\caption{Gluing two admissible surfaces  to form a third one. 
Note 
that only
part of the surface $S$ is drawn in the 
figure.}\lbl{gluesurface}
\end{figure}

\begin{proposition}
\lbl{cor.mult}
With the above assumptions  we have the following:
\begin{equation}
C_{f_1}^1 \cdot C_{f_2}^2 = C_{f_1 \cup f_2} 
\end{equation}
\end{proposition}

\begin{remark}
\lbl{rem.operad}
The above proposition makes necessary the existence 
of an operadic formalism of the above cocycles. Such a 
formalism,
which may make more transparent the relation with the ideas 
from
$2D$ gravity \cite{Ko1}, \cite{Ko2}, 
%%\cite{Pe}, 
\cite{Wi1},  \cite{Wi2}
 will be the subject of a future study.
\end{remark}

Before we state the next theorem, we need some notation: for a 
group $G$
and a positive integer $n$ let us {\em denote} by $\G_n G$ the 
(abelian)
quotient $G(n)/G(n+1)$. Let $\A^{conn}(\phi)$ denote the subspace 
of
$\A(\phi)$ consisting of linear combinations of {\em connected} 
admissible 
graphs.
We also define two binary operations:
$[^0x,y]=x\otimes y$ and $[^1 x,y]=-y \otimes x$. Then, we have 
the following 
theorem:

\begin{theorem}
\lbl{thm.homo}
Given an admissible surface $f: \S \hookrightarrow M$, and a 
nonnegative
integer $m$, there is a  linear map:
\begin{equation}
\lbl{eq.ghom}
D_{f,m}: \G_{2m} \Tgg \otimes \BQ \to \G_m\A^{conn}(\phi)
\end{equation}
with the following properties:
\begin{itemize}
\item
$D_{f,m}$ is determined by the cocycle $C_{f,m}$ as follows:
\begin{equation}
\lbl{eq.cd}
D_{f,m}([x_1, [ x_2, \ldots, [x_{2m-1}, x_{2m}]] )=
- \sum_a
C_{f,m}(
[^{a(1)} y_1, [^{a(2)} y_2, \ldots, [^{a(2m-1)} y_{2m-1},
y_{2m}]]]) 
\end{equation}
where $x_i \in \Tgg$, $y_i=\tau(x_i) \in U$,
the summation runs over all functions
$a : \{1,2, \ldots, 2m-1 \} \to \{0,1 \}$, and we set
$[^0x,y]=x\otimes y$ and $[^1 x,y]=-y \otimes x$.
\end{itemize}
Assume in addition that $f$ is an admissible Heegaard surface. 
Then,
\begin{itemize}
\item
$D_{f,m}$ depends only on the associated Lagrangian pair $(L^+, 
L^-)$,
and will be denoted by $D_{\lpm,m}$. 
\item
$D_{\lpm,m}$  satisfies the following
symmetry property:
\begin{equation}
\lbl{eq.symm}
D_{\lpm,m}(a^{-1})=(-1)^m D_{\lmp,m}(a)
\end{equation}
\item
Assume that $f$ is the standard \Heg\ splitting of $S^3$. If $g 
\geq 5m+1$, 
then $D_{\lpm,m}$ is onto.
\item
For an arbitrary admissible \Heg\ surface $f$, let 
$D_m: (\G_{2m}\Tgg \otimes \BQ)^{Sp_g} \to \G_m\A^{conn}(\phi)$
denote the restriction of $D_{\lpm,m}$ on the symplectic 
invariant
part of its domain.
For $m=1$, the composition of  $D_1$ with the weight system
$W_\l : \G_1\A^{conn}(\phi) \to \BQ$ 
coincides with the restriction of $ -\frac{1}{24} d_1: \Gamma_g 
\to \BQ$
of \cite[section 5]{Mo1} to $(\G_2 \Tgg \otimes \BQ )^{Sp(\HQ)}$.
\end{itemize}
\end{theorem} 

\begin{remark}
\lbl{rem.Dexplicit}
The proof of theorem \ref{thm.homo} (and theorem \ref{thm.L} 
below) exhibits  an explicit construction of enough {\em stably } 
non-trivial elements of the lower central series quotients 
$\G_{2m}\Tgg$ of the Torelli group when $g\geq 5m+1$ to prove 
that the map $D_{\lpm,m}$ is onto for a standard \Heg\ splitting
of genus at least $5m+1$. This construction may prove useful in 
further
study of the Torelli group.
\end{remark}

\subsection{Plan of the proof}
\lbl{sub.plan}

In section \ref{sec.pre} we review the definition and a 
few 
essential
properties of the Johnson homomorphism and discuss invariant 
theory 
for the symplectic and general linear group. 
In section \ref{sec.proofs} we prove our main results.
In section \ref{sec.KL} we discuss analogous constructions 
for some other subgroups of the mapping class group.
In section \ref{sec.dis} we discuss related results by Hain 
\cite{Hain}
and Morita \cite{Mo5}. Finally in section \ref{sec.epil} 
we formulate an important question which will be studied in a 
subsequent 
publication.

\subsection{Acknowledgment}
The final part of the paper was written during the authors visit 
at
Waseda University in July 1996; we wish to thank the organizers 
and 
especially S. Suzuki for inviting us.
In addition, we wish to thank D. Vogan and R. Hain for 
encouraging 
conversations. 
We especially wish to thank S. Morita for enlightening and 
clarifying
conversations during the conference in Waseda University.
Finally, we wish to thank 
the {\tt Internet} for providing useful communication
for the two authors.

\section{Preliminaries}
\lbl{sec.pre}

\subsection{Generalities in group theory}
\lbl{sub.gen}

In this section we review some general facts about group 
cohomology
of discrete groups.
Let $G$ be a discrete group and $\Q G$ the 
rational group ring of $G$. Let $IG$ (or simply $I$, in case 
we 
fix
 the group $G$) denote the augmentation ideal of 
$\Q G$ and $I^n$ the $n$-th power of $I$. We first recall 
the definition of the chain complex defined by the bar 
construction. 
Since we will only be dealing with coefficients with trivial 
$G$-action we can define, 
%We will define a 
%sequence 
%of {\em cochains} $\f_n\in C^n (G;I^n )$ where the 
%coefficients 
%$I^n$ are 
%considered as trivial $G$-modules.
 for a {\em trivial} $G$ module $M$, 
$C_n(G,M)=\text{Hom}(C_n(G),M)$, where  (using the bar 
construction)
$C_n (G)$ is the free $\Z$-module generated by $n$-tuples 
$[g_1 |\ldots |g_n ]$, where $g_i \in G$, and the boundary 
operator 
is defined by the formula:
\begin{equation*}
\partial [g_1 |\ldots |g_n ]=[g_2 |\ldots |g_n ]+\sum 
_{i=1}^{n-
1}(-1)^i [g_1 |\ldots |g_i g_{i+1}|\ldots |g_n ] +(-1)^n [g_1 
|\ldots 
|g_{n-1}]
\end{equation*}
As usual, we let $Z^n(G,M)$ (respectively, $B^n(G,M)$) denote 
the 
$n$-{\em cocycles} (respectively, $n$-{\em coboundaries}). 
For a 
cocyle 
$c_n \in Z^n(G,M)$,  let $[c_n] \in H^n(G,M)$ denote the 
associated
cohomology class.

Before proceeding to the main results of this section, we 
will 
prove a lemma which will be needed below in section 
\ref{sub.00}.
Define an involution $\gamma$ of $C_{\ast}(G)$ by the
formula 
\begin{equation*}
\gamma [g_1 |\ldots |g_n ]=(-1)^{\binom n2 }[g_n^{-1}|\ldots 
|g_1^{-1}]
\end{equation*}
We leave it to the reader to check that this is a 
chain map.
We have the following lemma:

\begin{lemma}
\lbl{lemma.involution}
If $G$ is a free abelian group, then $\gamma_{\ast}:H^n 
(G,\BZ)\to H^n 
(G,\BZ)$
is multiplication by $(-1)^n$.
\end{lemma}

\begin{pf}
Notice that for $n=1$ 
this is 
clearly
true. Since $H^1 (G)$ generates $H^{\ast} (G,\BZ)$, as an 
algebra,
 we will be done if we
prove that $\gamma^{\ast}$ preserves cup-products. Recall 
(see 
\cite{Mac})
that the formula
for cup-product, in the context of the bar construction, is:
\begin{equation*}
(\xi\cup\eta )[g_1 |\ldots |g_n |h_1 |\ldots |h_m ]=\xi [g_1 | 
\ldots 
|g_n ]\cdot\eta [h_1 |\ldots |h_m ]
\end{equation*}
where $\xi\in C^n (G ,\BZ), \eta\in C^m (G ,\BZ)$. Now we compute
\begin{align*}
\gamma^{\sharp}(\xi\cup\eta )[g_1 |\ldots |g_m |h_1 |\ldots 
|h_n 
]&
=(-1)^{\binom {m+n}2}(\xi\cup\eta )[h_n^{-1}|\ldots 
|h_1^{-1}|g_m^{-1}|\ldots |g_1^{-1}]\\ 
&=(-1)^{\binom{m+n}2}\xi[h_n^{-1}|\ldots 
|h_1^{-1}]\cdot\eta [g_m^{-1}|\ldots |g_1^{-1}]\\
&=(-1)^{\binom{m+n}2 +\binom m2 +\binom 
n2}\gamma^{\sharp}(\xi )
[h_1 \ldots |h_n ]\cdot\gamma^{\sharp}(\eta )[g_1 |\ldots 
|g_m ]\\
&=(-1)^{mn}(\gamma^{\sharp}(\eta )\cup\gamma^{\sharp}(\xi ))
[g_1 |\ldots |g_m |h_1 |\ldots |h_n ]
\end{align*}
So we see that, for any cohomology classes 
$\alpha\in H^n (G ,\BZ),\beta\in H^m (G ,\BZ)$, 
$\gamma^{\ast}(\alpha\cup\beta 
)=
(-1)^{mn}\gamma^{\ast }(\beta )\cup\gamma^{\ast }(\alpha )$. 
But 
now 
we 
just invoke the commutativity of cup-product on the 
cohomology 
level.
\end{pf}

Turning to our main results, we will now define  for every 
nonnegative integer $n$
a {\em cochain} $\f_n\in C^n (G;I^n )$, where $I^n$ is given 
the 
trivial $G$-module structure, as follows:
\begin{equation*}
\f_n [g_1 |\ldots |g_n ]=(1-g_1)\ldots (1-g_n)
\end{equation*}
   
Let $i_n :I^{n+1}\to I^n$ be the inclusion and $(i_n 
)_{\sharp}$ 
be the 
corresponding coefficient homomorphism of cochains.

\begin{lemma}
\lbl{lemma.fund}
\begin{equation*}
\delta (\f_n )=\begin{cases}
0& \text{$n$ even}\\
(i_n )_{\sharp}(\f_{n+1})& \text{$n$ odd}
\end{cases}
\end{equation*}
\end{lemma}

\begin{pf} We have the following formula for the coboundary:
\begin{equation}
\lbl{eq.del}
\begin{split}
\delta\f_n ([g_1 |\ldots |g_{n+1} ])& =\f_n ([g_2 |\ldots 
|g_{n+1}])+\sum _{i=1}^n (-1)^i \f_n ([g_1 |\ldots |g_i 
g_{i+1}|\ldots 
|g_{n+1}])\\
&\quad +(-1)^{n+1}\f_n ([g_1 |\ldots |g_n ])\\
& =(1-g_2)\ldots (1-g_{n+1})+\sum_{i=1}^n(-1)^i (1-g_1)\ldots 
(1-
g_i 
g_{i+1})\ldots (1-g_{n+1})\\
&\quad +(-1)^{n+1}(1-g_1)\ldots (1-g_n)
\end{split}
\end{equation}
Now making the substitution $1-g_i g_{i+1}=-(1-g_i)(1-
g_{i+1})+
(1-g_i)+(1-g_{i+1})$, the summation term in the above 
equation 
becomes
\begin{align*}
&\sum_{i=1}^n (-1)^i (1-g_1)\ldots (1-g_i g_{i+1})\ldots 
(1-g_{n+1}) \\
=& \sum_{i=1}^n (-1)^i (1-g_1)\ldots (-(1-g_i)(1-g_{i+1})+(1-
g_i)
+(1-g_{i+1}))\ldots (1-g_{n+1})\\
=&\sum_{i=1}^n (-1)^{i+1} (1-g_1)\ldots (1-g_i)(1-
g_{i+1})\ldots 
(1-g_{n+1})\\& +\sum_{i=1}^n (-1)^i (1-g_1)\ldots 
\widehat{(1-g_i 
)}\ldots (1-g_{n+1})+\sum_{i=1}^n (-1)^i (1-g_1)\ldots 
\widehat{(1-g_{i+1})}\ldots (1-g_{n+1})\\
=&(\sum_{i=1}^n (-1)^{i+1} )(1-g_1)\ldots (1-g_{n+1})-(1-
g_2)\ldots 
(1-g_{n+1})+(-1)^n (1-g_1)\ldots (1-g_n)
\end{align*}
Inserting this into equation \eqref{eq.del} we obtain
\begin{align*}
\delta\f_n ([g_1 |\ldots |g_{n+1} ])&=(\sum_{i=1}^n (-
1)^{i+1} 
)(1-
g_1)
\ldots (1-g_{n+1})\\
&=(\sum_{i=1}^n (-1)^{i+1} )\f_{n+1}([g_1 |\ldots |g_{n+1}])
\end{align*}
and the result follows.
\end{pf}

\begin{corollary}
\lbl{cor.trivial}
For  every even non-negative integer
$n$, there is  a well-defined cohomology class $[\f_n]$
in $H^n (G; I^n )$. 
\end{corollary}

\begin{pf}
Immediate by lemma \ref{lemma.fund} above.
\end{pf}

We now point out another useful fact. Recall from section 
\ref{sub.results}
 that for any positive integer 
$q$,  $G_q$ is defined inductively by $G_{q+1}=[G,G_q]$ with 
the
understanding that $G_1=G$. Recall also that $G(q)$ is the 
(normal) subgroup of $G$ that contains all elements of $G$ 
for 
which a 
nontrivial power
 belongs to $G_q$. It is easy to see that $\{G(n) \}_{n \geq 
1}$ is a 
decreasing sequence of normal subgroups of $G$ with the 
property:
$[G(n),G(m)] \subseteq G(n+m)$. 

\begin{lemma}
\lbl{lemma.comm}
If $g_i\in G(q_i)$ then $\f_n ([g_1 |\ldots |g_n ])\subseteq 
I^{q_1 
+\ldots +q_n}$.
\end{lemma}

\begin{pf}
With the notation $[g,h]=gh g^{-1} h^{-1}$, the following  
formula 
\begin{equation}
\lbl{eq.commute}
1-[g,h]=(-(1-g)(1-h)+(1-h)(1-g))g^{-1}h^{-1} \in I^2 
\end{equation}
shows that if $g \in G_q$, then $1-g \in I^q$.
The following formula
$$ 1-g^m =\sum_{i=1}^m(-1)^i \binom mi (1-g)^i$$
and the above shows by induction on $q$ that if $g^m \in G_q$ for 
some 
nonnegative integer $m$,  then 
\begin{equation}
\lbl{eq.power}
m(1-g) \equiv 1-g^m \bmod I^{q+1}
\end{equation}
thus deducing that $ m(1- g) \in I^q$, and since we are using 
rational 
coefficients,  this proves the lemma.
\end{pf}

\begin{corollary}
\lbl{cor.tt}
Given non-negative integers $n,q$,   $\f_n$ 
induces cochains 
\begin{equation}
\f_{n,q}\in C^n (G/G(q) ;I^n /I^{n+q-1})
\end{equation}
Furthermore,  for even $n, \f_{n,q}$ is a cocycle and, for 
odd 
$n$,
 $ \f_{n,2}$ is a cocycle. Moreover, $\f_{n,2}$ is 
multilinear.
\end{corollary}

Now suppose that $\B$ is a  vector space (over $\BQ$),
 carrying a decreasing filtration $\F_\ast \B$  
and $\r :\Q G\to \B$ is a linear 
map preserving the filtration, i.e. $\r (I^n )\subseteq \F_n 
\B$. 
Suppose 
also that the filtration of $\B$ is {\em $p$-step}, for some 
positive 
integer $p$, i.e. $\F_i \B =\F_{i+1} \B$ unless $p$ divides 
$i$. 
Now 
$\f_{n,q}$ induces, via $\r$ a cochain 
$\g_{n,q}\in C^n (G/G(q) ;\F_n \B 
/\F_{n+q-1} \B)$. Let $\G_n\A =\F_{pn} \B/\F_{pn+1} \B$.

We will consider the cochains $\g_{pn-q,q+2}\in C^{pn-
q}(G/G(q+2);\G_n\A )$, for $0\leq q<p$, (since $\F_{pn-
q}\B /\F_{pn+1}\B=\G_n\A$). These cochains 
are cocycles if either $pn-q$ is odd and $q=0$, or
$pn-q$ is even and $0\leq q<p$.

\begin{proposition}
\lbl{prop.coh}
If $pn-q$ is even and $0\leq q<p-1$, then
 $[s^{\sharp}\g_{pn-q,q+2}]= 0\in H^{pn-q}(G/G(q+3);\G_n\A 
)$, 
where $s:G/G(q+3)\to G/G(q+2)$ is the obvious projection.
\end{proposition} 
\begin{pf}
Consider the cochain $\g_{pn-q-1,q+3}\in C^{pn-q-
1}(G/G(q+3);\G_n\A 
)$.
 It follows from lemma~\ref{lemma.fund} that, when $pn-q$ is 
even
\begin{equation}
\lbl{eq.coc}
\delta\g_{pn-q-1,q+3}=s^{\sharp}\g_{pn-q,q+2}
\end{equation}
\end{pf}

\begin{corollary}
\lbl{cor.step}
With the above notation, if $pn$ is even, we have that
$[s^{\sharp}\g_{pn,2}]= 0\in H^{pn}(G/G(3);\G_n\A )$
\end{corollary}

We can also identify a family of {\em secondary} 
cohomology 
classes, although we will not, at this time, explore the 
application of these to our considerations.
 Define, for $pn-q$ odd,  
$\mu_{q,n}=j^{\sharp}\g_{pn-q,q+2}\in 
C^{pn-q}(G(q+1)/G(q+2);\G_n\A)$,
 where $j:G(q+1)/G(q+2)\to G/G(q+2)$ is the obvious 
inclusion.
Since $s\circ j$ is trivial, it follows from 
equation~\eqref{eq.coc} that
$\mu_{q,n}$ is a cocycle.

Clearly $[\mu_{q,n}]\in H^{pn-q}(G(q+1)/G(q+2);\G_n\A )$ and 
$[\g_{pn-q+1,q+1}] \in H^{pn-q+1}(G/G(q+1);\G_n\A )$, for $q>0$
 are related by transgression in
the fibration $G(q+1)/G(q+2)\to G/G(q+2)\to G/G(q+1)$, 
but neither one
is determined by the other. The extra information is encoded 
in the particular
cocycle representatives. If $q=0$, then this gives nothing 
new since $[\mu_{0,n}]=[\g_{pn,2}]$.

We end this section with a lemma that will be used in the proof 
of
theorem \ref{thm.homo}.
Recall first the map $G \to \BQ G$ given by $ g \to 1-g$.
According to lemma \ref{lemma.comm} for every positive integer
$n$, we get an  induced map $G(n) \to I^n$, and thus a linear 
map:
\begin{equation}
\lbl{eq.igr}
\G_n G \otimes \BQ \to I^{n}/I^{n+1}
\end{equation}
Note that addition in $\G_n G \otimes \BQ$ is given by group
multiplication in $G$.

These maps can be assembled together in the following way.
Recall first that $\G G \otimes \BQ$ can be given the structure 
of a 
graded Lie 
algebra (over $\BQ$).  Let $U(\G G \otimes \BQ)$ denote the 
universal 
enveloping 
algebra. Note that $U(\G G \otimes \BQ)$ is a Hopf algebra.
Note also that $\BQ G $ is a filtered algebra with respect to 
powers of the 
augmentation ideal. Let 
 $\G \BQ G$ be the associated graded algebra, i.e. 
$\G_n \BQ G = I^n/I^{n+1}$.
Note that $\BQ G$ is a Hopf algebra with comultiplication defined 
by
$\Delta(g)=g \otimes g$ for $g \in G$. Then the maps of equation 
\eqref{eq.igr} induce a map:
\begin{equation}
\lbl{eq.igra}
U(\G  G \otimes \BQ) \to \G \BQ G
\end{equation}
This map was shown by Jennings (see Quillen \cite{Qu}) to be a Hopf algebra
isomorphism. In particular, the primitive elements of $\G \BQ G$
are isomorphic to the Lie algebra $\G G \otimes \BQ $.

We end the section with  the following lemma:

\begin{lemma}
\lbl{lemma.ci}
For $x_i \in G$ we have the following identity in the graded
quotient $I^n/I^{n+1}$:
\begin{equation}
1-[x_1, \ldots, [x_{n-1},x_n]]=(-1)^{n-1}
\sum_{a}\upa^{a(1)}(z_1, \ldots, \upa^{a(n-1)}(z_{n-1}, z_n))
\end{equation}
where $z_i = 1-x_i \in I$ and the summation is over all functions
$a : \{1,2, \ldots, n-1 \} \to \{0,1 \}$, and we set
$\upa^0(a,b)=ab$ and $\upa^1(a,b)=-ba$.
Furthermore, the map \eqref{eq.igr} is a linear map.
\end{lemma}

\begin{pf}
Using the identity \eqref{eq.commute}
the first part follows by induction on $n$. Indeed, 
\eqref{eq.commute} implies that
\begin{eqnarray*}
1-[x_1, x_2] & \equiv & -(1-x_1)(1-x_2) + (1-x_2)(1-x_1) \bmod 
I^3 \\
             & = & -(\upa^0(1-x_1, 1-x_2) + \upa^1(1-x_1, 1-x_2)) 
\bmod I^3
\end{eqnarray*}
which concludes the proof of the first part for $n=1$. The 
induction step
follows the same way using identity \eqref{eq.commute}.

The second part follows immediately using the following identity:
\begin{equation*}
1-ab=(1-a)b + (1-b)
\end{equation*}
\end{pf}

\subsection{A review of finite type invariants of integral 
homology 3-spheres}
\lbl{sub.fti}

In this section we review some essential properties of \fti s of 
\ihs s
that will be used in the present paper.

We begin by recalling the definition of the map \eqref{eq.GA}
from  \cite{Oh}, \cite{GO1}:  
for an admissible (i.e., trivalent, vertex oriented) graph $G$ 
with $3m$ edges and $2m$ vertices,  let $L_w(G)$ denote the 
(linear combination of $2^{2m}$)
algebraically split links in $S^3$ with framing $f=+1$ on
each component obtained by choosing some of the vertices of 
$G$, 
replacing each non-chosen vertex of  
$G$ by a Borromean
ring, each chosen vertex by a trivial $3$-component link and 
each 
edge of $G$ by a band as in  figure 
\ref{graphlink}. The coefficient of that term is $(-1)^k$, 
where 
$k$ 
is the number of chosen vertices.
Due to the fact that $G$ is an abstract (non-embedded in 
$S^3$ 
graph),
the links whose sum with signs is $L_w(G)$ are  not well 
defined 
(modulo isotopy). Nevertheless, (with the notation of 
\cite{GL1},  \cite{GO1}) one can associate a well defined element 
$[S^3, L_w(G), f] \in \Gas {3m}$ in the associated  graded 
space.  
This map was shown in   \cite{Oh} (see also \cite{GL1})
to be onto. Furthermore, in \cite{GO1} it was shown that
it actually descends to a map $\G_{m}  \A(\phi) \to \Gas {3m}$ 
which, therefore, is also onto. This defines the map 
\eqref{eq.GA}.
According to the {\em fundamental theorem}
 of \fti s of \ihs s 
\cite{LMO}, \cite{L}, the map \eqref{eq.GA} is one-to-one,  
and therefore a vector space isomorphism. 
The isomorphisms of equation \eqref{eq.GA}  can
be assembled together for various $m$. Indeed, $\A(\phi)$ is 
equipped with
a multiplication (induced by the disjoint union of graphs) 
and a 
comultiplication (induced by all ways of splitting a graph 
into its 
connected components), compatible with the grading, thus 
giving 
$\A(\phi)$
the structure of a commutative cocommutative Hopf algebra.
Let $\hat{\A}(\phi)$ denote the completion.
Furthermore, $\M$ is equipped with a multiplication 
(induced by connected sums of \ihs s) and a comultiplication
defined by $\Delta(M)=M \otimes M$ for an \ihs\ $M$, thus
giving $\Fas {\star}$ and $\Gas {\star}$ the structure of
a commutative cocommutative Hopf algebra. 
The above mentioned results of
\cite{LMO}, \cite{L} additionally imply that the maps 
\eqref{eq.GA} combine 
to give an isomorphism of Hopf algebras  $\hat{\A}(\phi) \to \Gas 
\ast$.

\begin{figure}[htpb]
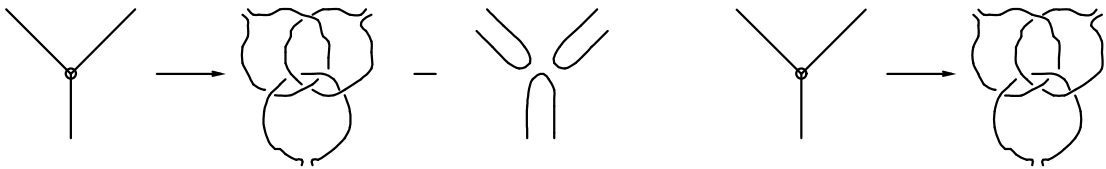

$$ \printname{graphlink2}
	\setlength{\unitlength}{0.03\standardunitlength}
	\begin{array}{c}  \hspace{-1.7mm}
        	\raisebox{-8pt}{\input draws/graphlink2.tex }
        	\hspace{-1.9mm}
	\end{array}
 $$
\caption{Two maps from admissible graphs to (linear combinations 
of) 
algebraically
split links in $S^3$. The map on the left is denoted by $G \to 
L_w(G)$ and the one on the right is denoted by $G \to 
L_b(G)$.}\lbl{graphlink}
\end{figure}

In the rest of this section we recall several facts about
the combinatorics of \fti s that will be used 
{\em exclusively} in the proof of theorem \ref{thm.L}. The reader 
may choose to
postpone them until needed.

%\begin{figure}[htpb]
%$$ \eepic{ASIHXvar}{0.03} $$
%\caption{The $AS$ and $IHX$ relations for vertex oriented graphs
%with univalent and trivalent vertices.}\lbl{ASIHXvar}
%\end{figure}

We begin with an equivalent description of the Hopf algebra 
$\A(\phi)$
taken from \cite{GO1}. It turns out, \cite{GO1} that there is a 
vector space
(over $\BQ$) $\A_b(\phi)$ on the set of vertex oriented graphs 
with univalent
and trivalent vertices only, modulo an appropriate set of 
relations, described
in detail in \cite{GO1}, together with a {\em deframing } map: 
\begin{equation}
\lbl{eq.F}
F: \A(\phi) \to \Ab
\end{equation}
defined as follows: for a vertex oriented trivalent graph $\Ga$  
$$F(\Ga)= \sum_{s: v(\Ga) \to \{0,1 \} } (-1)^{|s^{-1}(1)|} 
\Ga_s$$ where 
$\Ga_s$ is
obtained by splitting $\Ga$ along every vertex $v$ such that 
$s(v)=1$.  
We will not use explicitly the set of relations in $\A_b(\phi)$; 
note however
 \cite{GO1} that  $F$ is a vector space 
isomorphism, thus giving $\Ab$ the structure of a graded Hopf 
algebra.
 For a trivalent vertex oriented graph $\Ga$, let $\Ga_w$ 
(respectively,
$\Ga_b$) denote the associated element  in  $\A(\phi)$ 
(respectively,
$\Ab$). The subscripts $w,b$ denote {\em white} and {\em black} 
vertices 
respectively, the terminology is taken from \cite{GO1}. 
Given a trivalent vertex oriented graph $\Ga$, the associated 
elements
under the maps  $\G\A(\phi) \to \G\M$, $\G\Ab \to \G\M$ 
are shown in the left and right hand of figure 
\ref{graphlink} and are denoted by $\Ga \to [S^3, L_w(\Ga), +1]$ 
and
 $\Ga \to [S^3, L_b(\Ga), +1]$ respectively.
Putting together the isomorphism of equation \eqref{eq.GA} with 
the
isomorphism \eqref{eq.F} we get an isomorphism:

\begin{equation}
\lbl{eq.FF}
F ' : \G_n \Ab \simeq \Gas {3n} 
\end{equation}

$\A(\phi)$ (respectively, $\Ab$) has a naturally defined ideal:
$\Yw= \Theta_w \cdot \A(\phi)$ (respectively, $\Yb= Y_b \cdot 
\Ab$)
where $\Theta_w$ (respectively, $Y_b$) are the obvious generators
of the degree $1$ parts $\G_1\A(\phi)$ (respectively, $\G_1\Ab$).
Note that the ideals $\Yw, \Yb$ correspond under the isomorphism
$F$ of equation \eqref{eq.F}, since 
$F(\Theta_w)
=2 Y_b$ (see \cite{GO1}). Let $\At=\A(\phi)/\Yw  , \Abt=\Ab/\Yb $
denote the  quotient spaces.

We now have the following very useful lemma:

\begin{lemma}
\lbl{lemma.wb}
Let $\Ga$ denote a {\em connected} vertex oriented trivalent 
graph of 
degree $n \neq 1$,
and let $a \in \G_n\A^{conn}(\phi)$ be such that:
\begin{equation}
\Ga_b \bmod \Yb = F(a \bmod \Yw) \in \G_n\Abt 
\end{equation}
Then, we have that $\Ga_w= a$.
\end{lemma}

\begin{pf}
Recall first from \cite{GO1} that 
$F(\Ga_w) \equiv \Ga_b \bmod \Yb$. In fact more is true: namely
$F(\Ga_w)= \Ga_b + k(\Ga) (\coprod_n Y_b)$ for some integer 
$k(\Ga)$. 
This follows
by definition of the map $F$ of equation \eqref{eq.F} and the 
relations in $\G_n\A_b(\phi)$.
Therefore, we have that $F(\Ga_w) \equiv F(a) \bmod \Yb$, thus
$\Ga_w \equiv a \bmod \Yw$, thus $\Ga -a \equiv 0 \in \G_n \At$. 
Recall that $\A(\phi)$ is graded by $\G$, and thus filtered,
where $\F_n \A(\phi)= \oplus_{k \geq n} \G_k\A(\phi)$.
Using the fact that
$\A^{conn}(\phi) = \A(\phi)/(\F_1\A(\phi) \cdot \F_1\A(\phi))$,
we can see that, for any
$k \neq 1$, there is 
 an {\em onto} map $\G_k\At \to \G_k \A^{conn}(\phi)$.
Moreover, the composite map $ \G_k\A^{conn}(\phi) \to \G_k\At  
\to \G_k \A^{conn}(\phi)$ is the identity map on 
$\G_k\A^{conn}(\phi)$.
This finishes the proof of the lemma.
\end{pf}  

\begin{remark}
\lbl{rem.wb}
For $n=1$ the above lemma is obviously not true, since we can
take  $\Ga= \Theta$ and $a=0$.
\end{remark}

We now recall a few essential facts from \cite{GL3} relating the 
various filtrations on $\M$. These, too, will only be used in the 
proof
of theorem \ref{thm.L}.
 Following the notation of \cite{GL3}, 
consider   $f:\S_g\hookrightarrow M$  an admissible genus $g$ 
surface.
We need to recall from \cite[section 1.3]{GL3} an important 
subgroup $\Lgg$ of
the mapping class group. We call a Lagrangian $L \subseteq H$ 
$f$-{\em compatible} if $L=(L\cap L^+ )+(L\cap L^- )$, where
$(L^+, L^-)$ is the associated Lagrangian pair of the admissible 
surface $f$. For example, $L^+$ and $L^-$ themselves are
$f$-compatible.  For any 
$f$-compatible Lagrangian $L$, let $\Lgg$ denote the subgroup of 
the mapping
class group generated by {\em Dehn twists} along simple closed 
curves that
homologically represent elements of $L$.
%%An example of an $f$-compatible Lagrangian is $L^+$ (and also 
%%$L^-$);
%%we denote by $\cal L_g^{L^+}$ the corresponding subgroup of the 
%%mapping class group.  
Let 
 $\J$ denote any of the subgroups $\Tgg, \Kgg, \Lgg$ of the 
mapping
class group. 
 Consider the maps $\PJf: \BQ \J \to \M$ defined 
the same way as the map $\PTf$ of equation \eqref{eq.fun}.
Let $\G\PJf: \G\BQ\J \to \G^J\M$ denote the associated graded 
maps. Recall from
section \ref{sub.gen} that $\G\BQ\J$ is a coalgebra, and so is 
$\G^J\M$
(with the comultiplication of $\G^J\M$ induced by the one on 
$\M$).
Then, with the above conventions we have the following lemma:

\begin{lemma}
\lbl{lemma.coalg}
For $\J$ as above, the maps $\G\PJf : \G\BQ\J \to \G^J\M$ are 
maps
of coalgebras.
\end{lemma}

\begin{pf}
Recalling that  the coproduct
on $\BQ\J$  is defined by $\Delta(g)=g \otimes g$,
and the coproduct on $\M$ defined by $\Delta(M)=M \otimes M$, it 
follows that
the map $\PJf : \BQ \J \to \M$ preserves the coalgebra structure. 
Therefore,
the associated graded map preserves the coalgebra structure as 
well.
\end{pf}

\begin{corollary}
\lbl{cor.prim}
For $\J$ as above, we get an induced map:
\begin{equation}
\lbl{eq.prim}
\phi^J_f : \G \J \otimes \BQ \to \G\A^{conn}(\phi)
\end{equation}
\end{corollary}

\begin{pf}
We need to recall from \cite{GL3} the following filtrations on 
$\M$: for
a nonnegative integer $m$, 
$\FT m$ (respectively, $ \FK m, \FL m$) are defined to be the 
span
of the images (over all admissible surfaces $f$) of $\PTf (I 
\Tgg)^m$
(respectively, $\PKf (I \Kgg)^m$, $\PLf (I \cal L_g^{L^+})^m$).
Let $\FHL m$ denote the span of the images, for Heegaard surfaces 
$f$ and $f$-compatible Lagrangians $L$ of $\PLf ((I \Lgg)^m)$.  
In \cite{GL3}
we showed that $\FHL m = \FL m$, and from now on we will identify
these two filtrations.
The  filtrations considered above
 can be compared to the $\F^{as}$ filtration on $\M$
as follows \cite[corollary 1.20]{GL3}:

\begin{equation*}
\FK {m}  \subseteq \FT {2m} = \FT {2m-1} = \FL {3m} = \FL {3m-2}= 
\Fas {3m} 
%\FT {2m} & \hookrightarrow & \Fas {3m} \\
% \FL {3m} & \hookrightarrow & \Fas {3m} 
\end{equation*}
inducing associated graded maps

\begin{equation*}
\GK {m}  \to \GT {2m} = \GT {2m-1} = \GL {3m} = \GL {3m-2}= \Gas 
{3m}
\end{equation*}

%\begin{eqnarray*}
%\GK {m} & \to & \Gas {3m} \\
%\GT {2m} & \to & \Gas {3m} \\
%\GL {3m} & \to & \Gas {3m} 
%\end{eqnarray*}

Using  the isomorphism of Hopf algebras $\G_{\star}\A(\phi) 
\simeq
\Gas {\star} $
the above graded maps and lemma \ref{lemma.coalg} show that there
are coalgebra maps $\G\PJf : \G \BQ\J \to \G\A(\phi)$, which 
induce 
maps $\phi^J_f$ on the primitive elements.
Recall finally from section \ref{sub.gen} that the subspace of 
primitive elements of 
$\G\BQ\J$ is the Lie algebra $\G \J \otimes \BQ$,
and that the subspace of primitive elements of $\G\A(\phi)$ is 
$\G\A^{conn}(\phi)$.
\end{pf}

\begin{remark}
\lbl{rem.LKT}
$\G^J\M$ is  a Hopf algebra and the map $\G^J\M\to\G\A(\phi)$
 is a Hopf algebra isomorphism for $\J=\Tgg$ and $\Lgg$ and a 
Hopf algebra epimorphism for $\J =\Kgg$, see \cite{GL2}.  
Furthermore, $\G \BQ \J$ is also
a Hopf algebra, see section \ref{sub.gen}. The map $\G\PJf$ 
however is {\em  not} an algebra map.
\end{remark}

We close the section with the following lemma which will be 
used in the proof
of theorem \ref{thm.homo}. 
\begin{lemma}   
\lbl{lemma.incl}
Given an admissible surface and compatible
Lagrangian $L$,  we have 
inclusions $\Lgg \supseteq \Kgg \subseteq \Tgg$.
Then:
\begin{equation}
\lbl{eq.restr}
\PLf |_{\BQ \Kgg}=\PTf |_{\BQ \Kgg}
\end{equation}
\end{lemma}

\begin{pf}
This is a straightforward consequence of the  definitions of the 
maps.
\end{pf}

\subsection{Johnson's homomorphism and representation theory}
\lbl{sub.johnson}

In  this section 
 we review well known properties of the Johnson homomorphism, 
\cite{Jo1}, and some essential facts about representation 
theory
of the symplectic group.

We begin with the following: 
\begin{remark}
\lbl{rem.bou}
Even though in the present paper we are
interested mainly in closed surfaces embedded in closed 
3-manifolds, we will,
for a variety of reasons, also consider surfaces {\em with 
boundary}.
These reasons include $(a)$ historical traditions \cite{Jo1}, 
\cite{Mo1},
\cite{Mo2}, $(b)$ technical reasons (the fact that the 
fundamental
group of a closed surface is not free, whereas that of one 
with 
boundary
is. Also, one can glue surfaces along boundary to increase 
the 
genus
and discuss stability problems, whereas there is no canonical 
way 
of 
increasing
the genus of closed surfaces), and $(c)$ modern 
interpretation in 
terms
of open string field theory \cite{Ko1}, \cite{Ko2}, 
\cite{Wi1}, 
\cite{Wi2}. 
For all of the above reasons, we usually first decorate 
surfaces 
by 
boundary
components
(or punctures), and only afterwards do we discuss closed 
surfaces. 
At any rate, 
the reader should keep in mind that there are
exact sequences that relate invariants of decorated surfaces 
to 
invariants
of closed surfaces.
\end{remark}

Let $\S_g$ denote a closed, oriented surface of genus $g$, 
and let 
$D 
\subseteq \S_g$ be a fixed embedded disk. Let $\S_{g,1}$ 
denote 
the 
associated surface $\S - \text{Int} D$ with one boundary 
component. 
Let
$\Ga_g$ (resp. $\Ga_{g,1}$) denote the mapping class group, 
i.e., 
the
 group of isotopy classes of 
orientation preserving surface diffeomorphisms (resp. that 
are 
identity
on the boundary). Let  $\T_g$ (resp. $\T_{g,1}$)  (the Torelli 
group) denote the 
subgroup
of $\Ga_g$ (resp. $\Ga_{g,1}$) of elements that act trivially on 
the homology 
of 
the 
surface. Let $H= H_1(\S_g, \BZ)$, and 
$\omega$ 
be the intersection form. 
Note that the inclusion $\S_{g,1} \hookrightarrow \S_g$ 
induces a
canonical isomorphism $H_1(\S_{g,1}, \BZ) \simeq H$, and in 
this 
section 
we will identify $H_1(\S_{g,1}, \BZ)$ with $H$.
The groups $\Ga_g, \Ga_{g,1}, \T_g, \T_{g,1}$ are related in 
the 
following (exact) commutative diagram \cite{Jo1}: 
$$
\begin{CD}
%%%@. @. 1 @. 1 \\
%%%@. @. @VVV @VVV \\
1 @>>> \pi_1(T_\S) @>>> \Tg @>>> \Tgg @>>> 1 \\
@.         @|         @VVV      @VVV      \\
1 @>>> \pi_1(T_\S) @>>> \Ga_{g,1} @>>> \Ga_g @>>> 1 \\
@.         @.           @VVV      @VVV      \\
@. @.                   Sp(H) @>>> Sp(H) 
\end{CD}
$$
where $T_\S$ denotes the unit tangent  bundle of the surface 
$\S_g$.
Note that the two rightmost vertical sequences are also short 
exact.

With the above notation, we can recall a few essential facts 
from 
representation theory. For proofs, we refer the reader to 
\cite{FH}.
Recall that for an abelian group $A$, we let $A_{\BQ}$ denote 
the 
rational 
vector space $A \otimes_{\BZ} \BQ$.
All the linear maps to be described in this section will be 
$Sp(\HQ)$ equivariant.
Recall first that the intersection form $ \omega \in 
\Lambda^2 
H^\ast$
induces an isomorphism $H \simeq H^\ast$.
 Let $H \to \L3 H$ denote the map defined by $x \to x \we 
\omega $
(where we think of $\omega \in \Lambda^2 H$, via the 
isomorphism
$H \simeq H^\ast$). In terms of a symplectic basis $\{x_i, 
y_i \}$
of $H$, the above map is given by $x \to \sum_i x \we x_i \we 
y_i$.
Let $U= \L3 H/H$ denote the quotient.
Note that our notation differs from \cite{Mo5} (Morita 
denotes 
$\frac{1}{2}
\L3 H/H$ by $U$). Since both Morita  and we are
only dealing with rational results, this is not really a 
problem.

We can think of $\L3 \HQ$ as a  quotient module of $\otimes^3 
\HQ$ 
(in the natural way), or
as a submodule of $\otimes^3 \HQ$ as follows:
we let  $\L3 \HQ \hookrightarrow \otimes^3 \HQ$ be the 
(one-to-one) map:
$ x_1 \we x_2 \we x_3 \to \sum_{s \in \text{Sym}_3} \sgn(s)
x_{s(1)} \we x_{s(2)} \we x_{s(3)}$, where $\text{Sym}_3$ 
denotes 
the
symmetric group and $sgn$ denotes the sign homomorphism. Note
that we do not divide out the above map by $1/6$. As a 
result, the
composite map 
$ \L3 \HQ \to \otimes^3 \HQ \to \L3 \HQ$ is multiplication by 
$6$.
Let $\otimes^3 H \to H$ denote the map $x_1 \otimes x_2 
\otimes 
x_3
\to \omega(x_1, x_2) x_3$, and
let  $ \L3 H \to H$ denote the composite map $\L3 H 
\hookrightarrow
\otimes^3 H \to H$. More explicitly, the above composite map
is the following: 
\begin{equation}
x_1 \we x_2 \we x_3 \to 2 (\omega(x_1, x_2) x_3 -
\omega(x_1, x_3) x_2 + \omega(x_2, x_3) x_1) 
\end{equation}
It is easy to see that there is a rational isomorphism
$\UQ \simeq \text{Ker}( \L3 \HQ \to \HQ)$.
%%Furthermore, given an admissible pair $(L^+, L^-)$ of 
%%Lagrangian 
%%subspaces of $H$, there are  onto maps $\UQ \to \L^3 
%%L^{\pm}_{\BQ}$
%%defined as follows:

In his pioneering work \cite{Jo1}, \cite{Jo2}
 D. Johnson described a homomorphism
$ \tau: \Tg \to \L3 H$. 
We briefly summarize its  properties here:
\begin{itemize}
\item $\tau$ is onto.
\item
$ \tau$ is equivariant with respect to the conjugation 
action of the mapping class 
group $\Ga_{g,1}$ on $\Tg$ and the natural action of the
symplectic group $Sp(H)$ on $\L3 H$.
\item
$\tau$ coincides, modulo $2$ torsion, with the abelianization 
of 
the
Torelli group. In fact, $[\Tg, \Tg] \subset \text{Ker}( \Tg 
\to
\L3 H)$ is a normal subgroup with quotient a $2$-group.
%%In particular, we have that $\Kgg=\Tgg(2)$ and
%%that $\tau: \Tgg/\Tgg(2) \simeq U$.
\item
$\tau $ factors through a map  $\tau: \Tgg \to U$ (denoted by 
the 
same 
name) making the following diagram
commute:
$$
\begin{CD}
1 @>>> \pi_1(T_\S) @>>> \Tg @>>> \Tgg @>>> 1 \\
@. @VVV    @V \tau VV    @V \tau VV \\
0 @>>> H @>>> \L3 H @>>> U @>>> 1
\end{CD}
$$
Furthermore, $\Tgg \to U$ is equivariant, onto, and coincides 
(modulo
$2$ torsion) with the abelianization of $\Tgg$. From this and the
preceding property  we have that $\Kgg=\Tgg(2)$ and
that $\tau: \Tgg/\Tgg(2) \simeq U$.
\item
$\tau$ is stable with respect to an inclusion $\S_{g,1} 
\hookrightarrow
\S_{g+h,1}$.
\end{itemize}

\subsection{Representations of the symplectic and the general 
linear
group}
\lbl{sub.gsp}

In this section we review a few essential facts about 
representations of
the symplectic and the general linear group.
The main result is proposition \ref{prop.gl} which will be 
used in 
the proof
of theorem \ref{thm.formula}.

Let $(H,\omega)$ denote a symplectic space, and assume a 
given 
splitting
$H=L^+ \oplus L^-$ into two Lagrangian subspaces. An example  
is
given by the Lagrangian pair of an admissible surface in an 
\ihs\ 
.
Let us denote $L^+$ by $V$. Then, the symplectic form induces
isomorphisms $L^- \simeq V^\ast$ and $H \simeq V \oplus 
V^\ast$.
In this section we will identify $L^+, L^-$ and $ H$ with $V, 
V^\ast$
and $V \oplus V^\ast$, respectively.

 Consider the subgroup of the symplectic group $Sp(H)$ formed 
by 
all
matrices of the form $ \tbyt {A} {0} {0} {(A^t)^{-1}} $, 
where 
$A \in GL(V)$ and $A^t$ stands for the transpose of a matrix 
$A$.
This subgroup of $Sp(H)$ is obviously isomorphic to the group 
$GL(V)$. Note that the action of $GL(V)$ on $H$ preserves the
decomposition $H= V \oplus V^\ast$, and, as a subgroup of 
$Sp(H)$, 
extends to an action on
$\L3 H$ and $U$.

We will be mainly concerned with describing a generating set 
for
the vector space of invariants $(\otimes^{2m} \UQ 
)^{GL(\VQ)}$.
First, however, we need to recall several ideas about 
irreducible 
representations of $Sp(\HQ)$ and $GL(\VQ)$. For more details 
see 
\cite{FH}
and \cite{KK}.

We begin by recalling some results about the invariant theory
of the symplectic group as formulated by Morita \cite{Mo5}.
There is a one-to-one correspondence 
between
irreducible $Sp(H_{\BC})$ representations and {\em dominant 
integral weights}
(with respect to a standard choice of a Weyl chamber). Let us 
denote by
$V(\l)$ the rational representation  associated to weight 
$\l$.
Dominant integral weights are parametrized as follows:
if $dim(\HQ)=2n$ and $\{\e_i \}_{i=1}^n$ is the set of 
dominant 
weights,
then every dominant integral weight can be written uniquely 
in the
form $\sum_{i} f_i \e_i$ for integers $f_i$ 
such that $f_1 \geq f_2 \ldots \geq f_n \geq 0$.
It is often customary to denote the representation $V(\l)$ by 
a 
{\em Young diagram} of $n$ rows with $f_i$ boxes on each row.
Due to  typographical limitations though, we will not denote 
them 
by Young diagrams.
In this language we have the following identifications (as 
$Sp(\HQ)$ 
representations):
\begin{eqnarray*}
\HQ & = & V(\e_1) \\
\L3 \HQ & = & V(\e_3) + V(\e_1) \\
\UQ & = & V(\e_3) 
\end{eqnarray*}

Before we state the next result, we need a few definitions:
a degree $m$ {\em linear chord diagram}  is an involution on 
the 
set
$\{1,\ldots, 2m\}$ without fixed points \cite{B-N}. Let 
$\G_m\D^l$ 
denote 
the vector
space over $\BQ$ on the set of  linear chord diagrams with 
$m$ 
chords. 
It is easy to see
that $\G_m\D^l$ is a vector space of dimension
 $(2m-1)!!=1.3.5. \ldots (2m-1)$.

The symplectic form gives a $Sp(\HQ)$ invariant map
 $\omega: \HQ \otimes \HQ \to \BQ$, thus given a degree $m$ 
linear 
chord 
diagram,
we get an induced $Sp(\HQ)$ invariant map $\otimes^{2m} \HQ 
\to 
\BQ$,
and dually (using the isomorphism $H \simeq H^\ast$ induced 
by  
the
symplectic form) a $Sp(\HQ)$ invariant element in 
$\otimes^{2m} 
\HQ$. The definition is clear from the example shown in 
figure \ref{tricd}. 
Thus we have a map:
\begin{equation}
\lbl{eq.gld}
\G_m\D^l \to (\otimes^{2m} \HQ)^{Sp(\HQ)}
\end{equation}
According to the {\em first fundamental theorem} of
representation theory (see \cite[p.167]{W}), the above map 
\eqref{eq.gld} is onto, and
according to the {\em second fundamental theorem} of 
representation
theory \cite[p.168]{W}), provided $dim(\HQ)=2n \geq 2m$, 
\eqref{eq.gld} is one-to-one
and therefore a vector space isomorphism.

\begin{figure}[htpb]
$$ \printname{tricd}
	\setlength{\unitlength}{0.03\standardunitlength}
	\begin{array}{c}  \hspace{-1.7mm}
        	\raisebox{-8pt}{\begingroup\makeatletter\ifx\SetFigFont\undefined
% extract first six characters in \fmtname
\def\x#1#2#3#4#5#6#7\relax{\def\x{#1#2#3#4#5#6}}%
\expandafter\x\fmtname xxxxxx\relax \def\y{splain}%
\ifx\x\y   % LaTeX or SliTeX?
\gdef\SetFigFont#1#2#3{%
  \ifnum #1<17\tiny\else \ifnum #1<20\small\else
  \ifnum #1<24\normalsize\else \ifnum #1<29\large\else
  \ifnum #1<34\Large\else \ifnum #1<41\LARGE\else
     \huge\fi\fi\fi\fi\fi\fi
  \csname #3\endcsname}%
\else
\gdef\SetFigFont#1#2#3{\begingroup
  \count@#1\relax \ifnum 25<\count@\count@25\fi
  \def\x{\endgroup\@setsize\SetFigFont{#2pt}}%
  \expandafter\x
    \csname \romannumeral\the\count@ pt\expandafter\endcsname
    \csname @\romannumeral\the\count@ pt\endcsname
  \csname #3\endcsname}%
\fi
\fi\endgroup
\begin{picture}(9368,2488)(0,-10)
\thicklines
\put(1800.000,1239.000){\arc{1800.000}{3.1416}{6.2832}}
\put(2700.000,1239.000){\arc{600.000}{3.1416}{6.2832}}
\put(2700.000,1239.000){\arc{1800.000}{3.1416}{6.2832}}
\put(6449,1700){\ellipse{1530}{1530}}
\path(600,1239)(3900,1239)
\path(5670,1689)(7170,1689)
\path(5850,1539)	(5920.898,1575.097)
	(5967.188,1609.785)
	(5992.383,1646.581)
	(6000.000,1689.000)

\path(6000,1689)	(5992.383,1731.419)
	(5967.188,1768.215)
	(5920.898,1802.903)
	(5850.000,1839.000)

\path(5971.741,1817.109)(5850.000,1839.000)(5947.117,1762.395)
\path(7050,1539)	(6979.102,1575.097)
	(6932.812,1609.785)
	(6907.617,1646.581)
	(6900.000,1689.000)

\path(6900,1689)	(6907.617,1731.419)
	(6932.812,1768.215)
	(6979.102,1802.903)
	(7050.000,1839.000)

\path(6952.883,1762.395)(7050.000,1839.000)(6928.259,1817.109)
\put(5250,1689){\makebox(0,0)[lb]{$1$}}
\put(7425,1614){\makebox(0,0)[lb]{$2$}}
\put(0,39){\makebox(0,0)[lb]{$v_1 \otimes v_2 \otimes v_3 \otimes v_4 \otimes 
v_5 \otimes v_6 \to \omega(v_1, v_4) \omega(v_2, v_6) \omega(v_3, v_5)$}}
\end{picture} }
        	\hspace{-1.9mm}
	\end{array}
 $$
\caption{A  degree $3$ linear chord 
diagram and its associated trivalent graph. The trivalent 
graph 
$G$
comes equipped with an ordering $od_{V(G)}$ of its vertices,
as well as with a vertex orientation $or_{V(G)}$, indicated 
by a 
choice
of cyclic order of the edges around each vertex. The linear 
chord 
diagram corresponds (under the map \eqref{eq.gld}) to the 
contraction shown $\otimes^6 H \to \BQ$ shown  in  the 
figure.}\lbl{tricd}
\end{figure}

Our next goal is to describe the invariant spaces 
$(\otimes^{2m} \L3 \HQ)^{Sp(\HQ)}$.
In order to do so, we need one  more definition.
Given a degree $3m$ linear chord diagram, its associated 
trivalent
graph is defined  as follows: the set of vertices is the 
quotient
set $\{1,2, \ldots, 6m \}/\sim$\ , modulo the relation $3j-2 
\sim 3j-1 \sim 3j$
(for $1 \leq j \leq m$), and the set of edges is given by the 
quotient map
of the chord diagram. There is an orientation at every vertex,
induced by the ordering $3j-2 < 3j-1 < 3j$ of the edges around 
it.
  For an example, see figure \ref{tricd}. 
The trivalent graphs constructed above have extra data: they 
come
equipped with an ordering $od_{V(G)}$ of the vertices. Let 
$\G_m\A^{rp}$
denote the vector space over $\BQ$ on the set of isomorphism 
classes
of tuples $(G, od_{V(G)}, or_{V(G)} )$ (where $G$ is a 
trivalent 
graph
with $3m$ edges and $2m$ vertices, and $od_{V(G)}$ is an 
ordering 
of 
its
vertices, and $or_{V(G)}$ is  a vertex orientation of $G$, 
i.e., a
choice of cyclic order for the $3$ edges 
that emanate through each vertex of $G$) {\em divided out} by 
the 
usual
{\em antisymmetry} relation (denoted by $AS$ in figure 
\ref{ASIHX}).
The above discussion defines a map:
\begin{equation}
\lbl{eq.cdgr}
\G_{3m}\D^l \to \G_m\A^{rp}
\end{equation}
Due to the projection  %$\L3 H \to \otimes^3 H$ 
$\otimes^3 H \to \L3 H$ of section \ref{sub.johnson},
and the choice of cyclic order for the above mentioned 
trivalent graphs,
the map of equation \eqref{eq.gld} factors through the map of
equation \eqref{eq.cdgr}, thus inducing a map:
\begin{equation}
\lbl{eq.gal}
\G_m\A^{rp} \to (\otimes^{2m} \L3 \HQ)^{Sp(\HQ)}
\end{equation}
Furthermore, there are
natural inclusion maps $\G_m\A^{rp} \to \G_{3m} \D^l$ and \newline
$ (\otimes^{2m} \L3 \HQ)^{Sp(\HQ)} \to \G_m\A(\phi)$, such that
the composite maps $\G_m\A^{rp} \to \G_{3m} \D^l \to \G_m\A^{rp}$
and $ (\otimes^{2m} \L3 \HQ)^{Sp(\HQ)} \to (\otimes^{6m} \HQ)^{Sp(\HQ)}
\to (\otimes^{2m} \L3 \HQ)^{Sp(\HQ)} $ are multiplication by a nonzero scalar. 
Moreover, we  have two commutative diagrams:
$$
\begin{CD}
 \G_{3m} \D^l@>>>  (\otimes^{6m} \HQ)^{Sp(\HQ)} \\
@VVV            @VVV \\
\G_m\A^{rp} @>>> (\otimes^{2m} \L3 \HQ)^{Sp(\HQ)}
\end{CD}
\text{ and }
\begin{CD}
\G_m\A^{rp} @>>> (\otimes^{2m} \L3 \HQ)^{Sp(\HQ)} \\
@VVV            @VVV \\
 \G_{3m} \D^l@>>> (\otimes^{6m} \HQ)^{Sp(\HQ)} 
\end{CD}
$$
where the vertical maps on the left diagram are onto, and on the right
diagram are one-to-one.  
Using the fact that  the map of equation \eqref{eq.gld} is 
an isomorphism (provided that $n \geq m$) and the above commutative
diagrams  
it is easy to see that the map \eqref{eq.gal} is a vector
space isomorphism, provided that $n \geq 3m$. 

Finally, we describe a generating set for the invariant 
vector
space  $(\otimes^{2m} \UQ)^{Sp(\HQ)}$: let $\G_m\A^{rp,nl}$ 
denote 
the 
quotient space $\G_m\A^{rp}$ divided out by the subspace of 
tuples 
$(G, od_{V(G)}, or_{V(G)} )$, where $G$ contains a loop. For 
an example of a loop, see figure \ref{loop}.   

\begin{figure}[htpb]
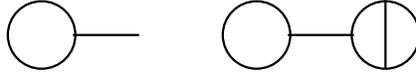

$$ \printname{loop}
	\setlength{\unitlength}{0.03\standardunitlength}
	\begin{array}{c}  \hspace{-1.7mm}
        	\raisebox{-8pt}{\input draws/loop.tex }
        	\hspace{-1.9mm}
	\end{array}
 $$
\caption{An example of a loop on the left and of a trivalent
trivalent graph containing a loop on the right.}\lbl{loop}
\end{figure}

Due to the projection $\L3 \HQ \to \UQ$
%$\L3 \HQ$ via the identification 
%$\UQ \simeq \text{Ker}(\L3 \HQ \to \HQ)$, 
it is easy to see that the map 
of equation \eqref{eq.gal} induces a  map:
\begin{equation}
\lbl{eq.gau}
\G_m\A^{rp,nl} \to (\otimes^{2m} \UQ)^{Sp(\HQ)}
\end{equation}
Using the  fact of an inclusion $\UQ \simeq \text{Ker}(\L3 \HQ \to \HQ)
\subseteq \L3 \HQ$ and the same reasoning as that  of $\L3 \HQ$ above,
implies that equation \eqref{eq.gau} is an isomorphism
provided  $n \geq 3m$. The above 
isomorphism has already been discussed in previous work of Morita 
\cite[page 11]{Mo5} which has been a source of inspiration 
for us.
The reason that we recall it here in detail is to show the 
similarities and differences between the invariant theory of 
$Sp(\HQ)$ and the invariant theory of the general linear group 
$GL(\VQ)
\hookrightarrow Sp(\HQ)$,
to which we now turn.

We begin by recalling  a few standard facts about representations 
of 
$GL(\VQ)$.
In this case too, it turns out that there is a one-to-one
correspondence between irreducible $GL(V_{\BC})$ 
representations 
and 
{\em
dominant integral weights} (with respect to a standard choice 
of a
Weyl chamber). Let us denote by $V(\l)$ the representation 
(over 
$\BQ$)
 associated
to weight $\l$.  Dominant integral weights are parametrized 
as
follows: if $dim(\VQ)=n$ and $\{\e_i \}_{i=1}^n$ is the set 
of 
dominant
weights, then every dominant integral weight can be written 
uniquely
in the form $\sum_{i} f_i \e_i$ for integers $f_i$ such that 
$f_1 
\geq
f_2 \ldots \geq f_n $.  (In case $f_n \geq 0$, such 
representations 
are
called {\em polynomial}, and a usual graphical way of 
representing
them is by a Young diagram of $n$ rows with $f_i$ boxes on 
each 
row.
However, if $V(\l)$ is not polynomial, there is {\em no} 
graphical 
way
to represent it. Since we will be dealing with non-polynomial
representations of $GL(\VQ)$, we will not use the Young 
diagram 
method.)
In this language we have the following lemma:

\begin{lemma}
\lbl{lemma.rep}
As representations of $GL(\VQ)$, we have the following 
decomposition:
\begin{eqnarray*}
\HQ & = & V(\e_1) + V(-\e_n) \\
\L3 \HQ & = &  V(\e_1 + \e_2 + \e_3 ) + V(\e_1)+
 V(\e_1 + \e_2 -\e_n) \\
      & + & V(-\e_n) + V(\e_1 -\e_{n-1} -\e_n)
+ V(-\e_{n-2} - \e_{n-1} -\e_n)\\
\UQ & = & V(\e_1 + \e_2 + \e_3 ) + V(\e_1 + \e_2 
-\e_n)+V(\e_1 
-\e_{n-1} -\e_n)
+ V(-\e_{n-2} - \e_{n-1} -\e_n)
\end{eqnarray*}
\end{lemma}
 
\begin{pf}
The first part is obvious. The second one follows using the 
identity
\begin{equation*}
\L3 (V + V^\ast)= \L3 V + \Lambda^2 V \otimes V^\ast + V 
\otimes
\Lambda^2 V^\ast + \L3 V^\ast
\end{equation*}
together with the facts that
\begin{eqnarray*}
\L3 V & = & V(\e_1 + \e_2 + \e_3)   \\
\Lambda^2 V \otimes V^\ast & = &
V(\e_1)+ V(\e_1 + \e_2 -\e_n)
\end{eqnarray*}
 The third part follows from first and the second.
We thank D. Vogan for a crash course on representation theory 
of 
$GL(\VQ)$.
\end{pf}

Before we state the main proposition of this section, we need 
a 
few
more definitions:
\begin{definition}
\lbl{def.cl}
A $2$-{\em coloring} of a linear chord diagram (respectively, 
of a  graph)  is an orientation
for each of the edges of the chord diagram (respectively, 
graph). For examples, see figure \ref{2Color}.
Let $\G_m\D^{l,cl}$ denote the vector
space over $\BQ$ on the set of $2$-colored  linear chord 
diagrams with $m$ chords. It is easy to see
that $\G_m\D^{l,cl}$ is a vector space of
dimension  $(2m-1)!! \ldots 2^{m}=1.3.5. \ldots (2m-1) 2^{m}$.
Let $\G_m\A^{rp,cl}$  denote the
 vector space over $\BQ$ on the set of isomorphism classes
of tuples $(G, od_{V(G)}, or_{V(G)}, cl_{E(G)} )$ 
(where $G$ is a trivalent graph
with $3m$ edges and $2m$ vertices, and $od_{V(G)}$ is an 
ordering of its
vertices,  $or_{V(G)}$ is a vertex orientation of $G$, i.e., 
a
choice of cyclic order for the $3$ edges 
that emanate through each vertex of $G$ and $cl_{E(G)}$ is a 
$2$-
coloring of
the edges $E(G)$ of $G$) divided out by the {\em colored
antisymmetry} relation (denoted by $AS^{cl}$) of figure 
\ref{AScl}.
Let $\G_m\A^{rp,nl,cl}$ denote the quotient space of 
$\G_m\A^{rp,cl}$ divided out by all graphs which contain a 
loop.
As is maybe apparent from the notation, the superscripts on 
$\A$ 
are explained as follows:  $rp$ stands for representation 
theory,
$nl$  stands for no loops and $cl$ stands for $2$-coloring 
(of the 
edges).
\end{definition}

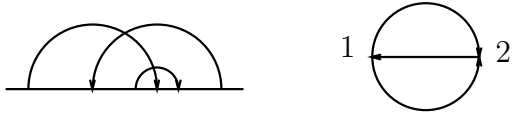
\begin{figure}[htpb]
$$ \printname{2Color}
	\setlength{\unitlength}{0.03\standardunitlength}
	\begin{array}{c}  \hspace{-1.7mm}
        	\raisebox{-8pt}{\begingroup\makeatletter\ifx\SetFigFont\undefined
% extract first six characters in \fmtname
\def\x#1#2#3#4#5#6#7\relax{\def\x{#1#2#3#4#5#6}}%
\expandafter\x\fmtname xxxxxx\relax \def\y{splain}%
\ifx\x\y   % LaTeX or SliTeX?
\gdef\SetFigFont#1#2#3{%
  \ifnum #1<17\tiny\else \ifnum #1<20\small\else
  \ifnum #1<24\normalsize\else \ifnum #1<29\large\else
  \ifnum #1<34\Large\else \ifnum #1<41\LARGE\else
     \huge\fi\fi\fi\fi\fi\fi
  \csname #3\endcsname}%
\else
\gdef\SetFigFont#1#2#3{\begingroup
  \count@#1\relax \ifnum 25<\count@\count@25\fi
  \def\x{\endgroup\@setsize\SetFigFont{#2pt}}%
  \expandafter\x
    \csname \romannumeral\the\count@ pt\expandafter\endcsname
    \csname @\romannumeral\the\count@ pt\endcsname
  \csname #3\endcsname}%
\fi
\fi\endgroup
\begin{picture}(7113,1533)(0,-10)
\thicklines
\path(2082.000,429.000)(2112.000,309.000)(2142.000,429.000)
\put(1212.000,309.000){\arc{1800.000}{3.1416}{6.2832}}
\path(2382.000,429.000)(2412.000,309.000)(2442.000,429.000)
\put(2112.000,309.000){\arc{600.000}{3.1416}{6.2832}}
\path(1182.000,429.000)(1212.000,309.000)(1242.000,429.000)
\put(2112.000,309.000){\arc{1800.000}{3.1416}{6.2832}}
\path(6582.000,879.000)(6612.000,759.000)(6642.000,879.000)
\put(5862.000,759.000){\arc{1500.000}{3.1416}{6.2832}}
\path(6642.000,639.000)(6612.000,759.000)(6582.000,639.000)
\put(5862.000,759.000){\arc{1500.000}{6.2832}{9.4248}}
\path(12,309)(3312,309)
\path(5232.000,789.000)(5112.000,759.000)(5232.000,729.000)
\path(5112,759)(6612,759)
\put(4662,759){\makebox(0,0)[lb]{$1$}}
\put(6837,684){\makebox(0,0)[lb]{$2$}}
\end{picture} }
        	\hspace{-1.9mm}
	\end{array}
 $$
\caption{An example of a $2$-coloring of a linear chord 
diagram and of its associated trivalent graph.}\lbl{2Color}
\end{figure}

\begin{figure}[htpb]
$$ \printname{AScl}
	\setlength{\unitlength}{0.03\standardunitlength}
	\begin{array}{c}  \hspace{-1.7mm}
        	\raisebox{-8pt}{\begingroup\makeatletter\ifx\SetFigFont\undefined
% extract first six characters in \fmtname
\def\x#1#2#3#4#5#6#7\relax{\def\x{#1#2#3#4#5#6}}%
\expandafter\x\fmtname xxxxxx\relax \def\y{splain}%
\ifx\x\y   % LaTeX or SliTeX?
\gdef\SetFigFont#1#2#3{%
  \ifnum #1<17\tiny\else \ifnum #1<20\small\else
  \ifnum #1<24\normalsize\else \ifnum #1<29\large\else
  \ifnum #1<34\Large\else \ifnum #1<41\LARGE\else
     \huge\fi\fi\fi\fi\fi\fi
  \csname #3\endcsname}%
\else
\gdef\SetFigFont#1#2#3{\begingroup
  \count@#1\relax \ifnum 25<\count@\count@25\fi
  \def\x{\endgroup\@setsize\SetFigFont{#2pt}}%
  \expandafter\x
    \csname \romannumeral\the\count@ pt\expandafter\endcsname
    \csname @\romannumeral\the\count@ pt\endcsname
  \csname #3\endcsname}%
\fi
\fi\endgroup
\begin{picture}(8346,939)(0,-10)
\thicklines
\path(238.154,486.231)(312.000,387.000)(293.538,509.308)
\put(462.000,449.500){\arc{325.000}{2.7468}{6.6780}}
\path(2430.462,584.308)(2412.000,462.000)(2485.846,561.231)
\put(2262.000,524.500){\arc{325.000}{2.7468}{6.6780}}
\path(4992.000,882.000)(5112.000,912.000)(4992.000,942.000)
\put(5112.000,612.000){\arc{600.000}{1.5708}{4.7124}}
\put(5112.000,612.000){\arc{600.000}{4.7124}{7.8540}}
\put(7212.000,612.000){\arc{600.000}{1.5708}{4.7124}}
\path(7332.000,942.000)(7212.000,912.000)(7332.000,882.000)
\put(7212.000,612.000){\arc{600.000}{4.7124}{7.8540}}
\path(12,912)(462,462)
\path(912,912)(462,462)
\path(1812,912)(2262,462)
\path(2712,912)(2262,462)
\path(2262,462)(2262,12)
\path(462,462)(462,12)
\path(5112,312)(5112,12)
\path(7212,312)(7212,12)
\put(1212,387){\makebox(0,0)[lb]{$+$}}
\put(3012,387){\makebox(0,0)[lb]{$=0$}}
\put(6012,387){\makebox(0,0)[lb]{$+$}}
\put(7962,462){\makebox(0,0)[lb]{$=0$}}
\end{picture} }
        	\hspace{-1.9mm}
	\end{array}
 $$
\caption{The colored antisymmetry relation on $2$-colored, 
vertex
oriented  trivalent graphs. All edges are $2$-colored by the 
choice 
of
an arrow. The left hand side corresponds to $8$ identities 
(for 
all 
possible $2$-colorings of the edges of the $Y$ graph). 
Similarly, 
the right
hand side corresponds to $2$ identities.}\lbl{AScl}
\end{figure}
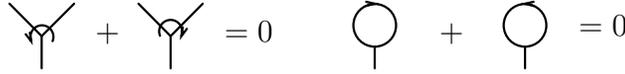

With the notation as in the above definition, we are ready to 
state
the main result of this section:

\begin{proposition}
\lbl{prop.gl}
For every non-negative integer $m$, there are onto maps:
\begin{eqnarray*}
\G_m\D^{l,cl} & \to & (\otimes^{2m}\HQ)^{GL(\VQ)} \\ 
\G_m\A^{rp,cl} & \to & (\otimes^{2m}\L3 \HQ)^{GL(\VQ)} \\
\G_m\A^{rp,nl, cl} & \to & (\otimes^{2m}\UQ)^{GL(\VQ)}
\end{eqnarray*}
which are isomorphisms provided that $n \geq m$, $ n \geq 3m$ 
and
$n \geq 3m$ respectively.
\end{proposition}

\begin{pf}
Recall first that $H= V \oplus V^\ast$, thus we get two 
$GL(\VQ)$
 invariant maps:
$\HQ \otimes \HQ \to \BQ$, obtained as follows:
\begin{eqnarray*}
\HQ \otimes \HQ \to & \VQ \otimes \VQ^\ast & \to \BQ \\
\HQ \otimes \HQ \to & \VQ^\ast \otimes \VQ & \to \BQ
\end{eqnarray*}
by projecting $\HQ \to \VQ$ or $\HQ \to \VQ^\ast$ 
appropriately. 
A choice between each of the above mentioned invariant maps 
will 
be 
associated with a $2$-coloring.
Therefore, given a degree $m$ $2$-colored linear chord 
diagram, we 
get
an invariant map $\otimes^{2m} \HQ \to \BQ$, and dually 
(using the
isomorphism $H \simeq H^\ast$) a $GL(\VQ)$ invariant element
in $\otimes^{2m}\HQ$. This defines a map:
$\G_m\D^{l,cl} \to (\otimes^{2m}\HQ)^{GL(\VQ)}$, which 
according 
to the
{\em first fundamental theorem} of invariant theory is onto, 
and according
to the {\em second fundamental theorem} (provided that $n \geq m$)
is one-to-one and thus a vector
space isomorphism. This shows the first part of proposition 
\ref{prop.gl}.
Next, recall the map of equation \eqref{eq.cdgr} that assigns 
to each
degree $3m$ linear chord diagram its associated degree $m$ 
trivalent graph.
This map by definition respects the $2$-coloring of linear 
chord diagrams
and trivalent graphs and therefore descends to onto maps:
\begin{equation}
\lbl{eq.cdgrcl}
\G_{3m}\D^{l,cl} \to \G_m\A^{rp,cl} \text{ and }
\G_{3m}\D^{l,cl} \to \G_m\A^{rp,cl,nl}
\end{equation}

Arguing as in equations \eqref{eq.cdgr} and \eqref{eq.gau}, 
due to the projections
 $\otimes^3 \HQ \to \L3 \HQ$ and $\otimes^3 \HQ \to \L3 \HQ$
we  the map $ \G_m\D^{l,cl} \to (\otimes^{2m}\HQ)^{GL(\VQ)}$
 induces quotient  maps:
\begin{equation}
\lbl{eq.hope}
\G_m\A^{rp,cl}  \to  (\otimes^{2m}\L3 \HQ)^{GL(\VQ)}
\text{ and }
\G_m\A^{rp,nl, cl}  \to  (\otimes^{2m}\UQ)^{GL(\VQ)}
\end{equation}

The same reasoning as that of equations \eqref{eq.cdgr} and
\eqref{eq.gau} together with the isomorphism of the first part
of the proposition implies the rest of  proposition \ref{prop.gl}. 
\end{pf}
 
%We claim that the isomorphism of the first part of the 
%proposition, 
%restricted to
%$\otimes^{2m} \L3 \HQ \hookrightarrow \otimes^{6m} \HQ$ 
%(respectively, $ \UQ \hookrightarrow \L3 \HQ $) descends to 
%maps:
%$\G_m\A^{rp,cl} \to (\otimes^{2m}\L3 \HQ)^{GL(\VQ)}$
%(respectively, $\G_m\A^{rp,nl,cl} \to (\otimes^{2m}
%\UQ)^{GL(\VQ)}$), which  are onto.
%Indeed, the first one is well defined and onto by the defining 
%property
%of $\L3 \HQ$. In order to show that the second one is well 
%defined,
%consider  a degree $m$ tuple
%$(G, od_{V(G)}, or_{V(G)}, cl_{E(G)} )$, where $G$ is a graph
%that contains a loop, and its associated invariant element
%$ \otimes^{3m} \HQ \to \BQ$. The restriction of this
%element to $\otimes^{2m} \UQ \to \BQ$ vanishes, due to
%the facts that $\UQ = \text{Ker}(\L3 \HQ \to \HQ)$
%and the fact that $G$ contains a loop, and thus
%contraction along the loop vanishes. This shows that
%the map $\G_m\A^{rp,nl,cl} \to (\otimes^{2m}\UQ)^{GL(\VQ)}$ is 
%well defined.
%The same argument also shows that it is onto.

%Since the map $\G_m\D^{l,cl} \to (\otimes^{2m}\HQ)^{GL(\VQ)}$ 
%is one-to-one provided that $n \geq m$, 
%it follows 
%from their definitions that the maps $\G_m\A^{rp,nl,cl} \to 
%(\otimes^{2m}\L3 
%\HQ)^{GL(\VQ)}$ and $\G_m\A^{rp,nl,cl} \to 
%(\otimes^{2m}\UQ)^{GL(\VQ)}$ are
%also one-to-one provided that $n \geq 3m$, \c\ {\bf It seems 
%again that one-one requires some more argument} thus finishing 
%the
%proof of proposition \ref{prop.gl}.

\begin{corollary}
\lbl{cor.gl}
In particular, for $m=1$, we have the following table of 
dimensions for the various invariant spaces:
$$
\vbox{
\offinterlineskip
\halign{\strut#&&\vrule#&\quad\hfil#\hfil\quad\cr
\noalign{\hrule}
&& $A$  && $dim(A)^{Sp(\HQ)}$ && $dim(A)^{GL(\VQ)}$     &\cr
\noalign{\hrule}
&& $\otimes^6 \HQ$  && 
$5!!=15 $ && $2^3 . 5!!= 120$  &\cr
\noalign{\hrule}
&& $ \otimes^2 \L3 \HQ $ && $2$ && $6$  &\cr
\noalign{\hrule}
&& $ \otimes^2 \UQ $ && $1$ && $4$   &\cr 
\noalign{\hrule}
}}
$$
In terms of the isomorphisms of proposition \ref{prop.gl}, 
the 
graphs
in figure \ref{Basis1} form basis for the invariant spaces 
$(\L3 \HQ)^{Sp(\HQ)},
(\UQ)^{Sp(\HQ)}, (\L3 \HQ)^{GL(\VQ)}$ and $(\UQ)^{GL(\VQ)}$.
\end{corollary}

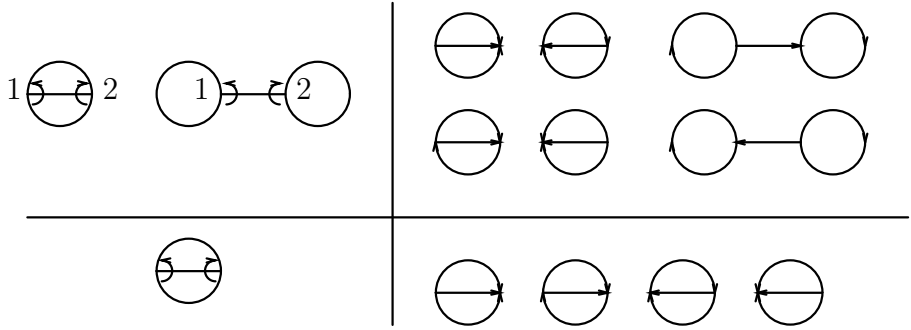
\begin{figure}[htpb]
$$ \printname{Basis1}
	\setlength{\unitlength}{0.03\standardunitlength}
	\begin{array}{c}  \hspace{-1.7mm}
        	\raisebox{-8pt}{\begingroup\makeatletter\ifx\SetFigFont\undefined
% extract first six characters in \fmtname
\def\x#1#2#3#4#5#6#7\relax{\def\x{#1#2#3#4#5#6}}%
\expandafter\x\fmtname xxxxxx\relax \def\y{splain}%
\ifx\x\y   % LaTeX or SliTeX?
\gdef\SetFigFont#1#2#3{%
  \ifnum #1<17\tiny\else \ifnum #1<20\small\else
  \ifnum #1<24\normalsize\else \ifnum #1<29\large\else
  \ifnum #1<34\Large\else \ifnum #1<41\LARGE\else
     \huge\fi\fi\fi\fi\fi\fi
  \csname #3\endcsname}%
\else
\gdef\SetFigFont#1#2#3{\begingroup
  \count@#1\relax \ifnum 25<\count@\count@25\fi
  \def\x{\endgroup\@setsize\SetFigFont{#2pt}}%
  \expandafter\x
    \csname \romannumeral\the\count@ pt\expandafter\endcsname
    \csname @\romannumeral\the\count@ pt\endcsname
  \csname #3\endcsname}%
\fi
\fi\endgroup
\begin{picture}(12612,4539)(0,-10)
\thicklines
\put(9751.000,3911.000){\arc{900.000}{3.1416}{6.2832}}
\put(11551.000,3911.000){\arc{900.000}{6.2832}{9.4248}}
\put(11551.000,2561.000){\arc{900.000}{6.2832}{9.4248}}
\put(9751.000,2561.000){\arc{900.000}{3.1416}{6.2832}}
\put(752.000,3235.000){\arc{900.000}{3.1416}{6.2832}}
\put(752.000,3235.000){\arc{900.000}{6.2832}{9.4248}}
\put(2551.000,761.000){\arc{900.000}{3.1416}{6.2832}}
\put(2551.000,761.000){\arc{900.000}{6.2832}{9.4248}}
\path(6870.000,4032.000)(6900.000,3912.000)(6930.000,4032.000)
\put(6450.000,3912.000){\arc{900.000}{3.1416}{6.2832}}
\path(6930.000,3792.000)(6900.000,3912.000)(6870.000,3792.000)
\put(6450.000,3912.000){\arc{900.000}{6.2832}{9.4248}}
\path(6871.000,2681.000)(6901.000,2561.000)(6931.000,2681.000)
\put(6451.000,2561.000){\arc{900.000}{3.1416}{6.2832}}
\path(8371.000,4031.000)(8401.000,3911.000)(8431.000,4031.000)
\put(7951.000,3911.000){\arc{900.000}{3.1416}{6.2832}}
\path(11971.000,4031.000)(12001.000,3911.000)(12031.000,4031.000)
\put(11551.000,3911.000){\arc{900.000}{3.1416}{6.2832}}
\path(11971.000,2681.000)(12001.000,2561.000)(12031.000,2681.000)
\put(11551.000,2561.000){\arc{900.000}{3.1416}{6.2832}}
\path(6871.000,581.000)(6901.000,461.000)(6931.000,581.000)
\put(6451.000,461.000){\arc{900.000}{3.1416}{6.2832}}
\path(6931.000,341.000)(6901.000,461.000)(6871.000,341.000)
\put(6451.000,461.000){\arc{900.000}{6.2832}{9.4248}}
\path(8371.000,581.000)(8401.000,461.000)(8431.000,581.000)
\put(7951.000,461.000){\arc{900.000}{3.1416}{6.2832}}
\path(9871.000,581.000)(9901.000,461.000)(9931.000,581.000)
\put(9451.000,461.000){\arc{900.000}{3.1416}{6.2832}}
\path(7531.000,3791.000)(7501.000,3911.000)(7471.000,3791.000)
\put(7951.000,3911.000){\arc{900.000}{6.2832}{9.4248}}
\path(6031.000,2441.000)(6001.000,2561.000)(5971.000,2441.000)
\put(6451.000,2561.000){\arc{900.000}{6.2832}{9.4248}}
\path(7531.000,2441.000)(7501.000,2561.000)(7471.000,2441.000)
\put(7951.000,2561.000){\arc{900.000}{6.2832}{9.4248}}
\path(7471.000,2681.000)(7501.000,2561.000)(7531.000,2681.000)
\put(7951.000,2561.000){\arc{900.000}{3.1416}{6.2832}}
\path(7531.000,341.000)(7501.000,461.000)(7471.000,341.000)
\put(7951.000,461.000){\arc{900.000}{6.2832}{9.4248}}
\path(9031.000,341.000)(9001.000,461.000)(8971.000,341.000)
\put(9451.000,461.000){\arc{900.000}{6.2832}{9.4248}}
\path(10471.000,581.000)(10501.000,461.000)(10531.000,581.000)
\put(10951.000,461.000){\arc{900.000}{3.1416}{6.2832}}
\path(10531.000,341.000)(10501.000,461.000)(10471.000,341.000)
\put(10951.000,461.000){\arc{900.000}{6.2832}{9.4248}}
\path(9331.000,3791.000)(9301.000,3911.000)(9271.000,3791.000)
\put(9751.000,3911.000){\arc{900.000}{6.2832}{9.4248}}
\path(9331.000,2441.000)(9301.000,2561.000)(9271.000,2441.000)
\put(9751.000,2561.000){\arc{900.000}{6.2832}{9.4248}}
\path(495.000,3417.000)(375.000,3387.000)(495.000,3357.000)
\put(375.000,3237.000){\arc{300.000}{4.7124}{7.8540}}
\path(1005.000,3357.000)(1125.000,3387.000)(1005.000,3417.000)
\put(1125.000,3237.000){\arc{300.000}{1.5708}{4.7124}}
\path(3195.000,3417.000)(3075.000,3387.000)(3195.000,3357.000)
\put(3075.000,3237.000){\arc{300.000}{4.7124}{7.8540}}
\path(3705.000,3357.000)(3825.000,3387.000)(3705.000,3417.000)
\put(3825.000,3237.000){\arc{300.000}{1.5708}{4.7124}}
\path(2805.000,882.000)(2925.000,912.000)(2805.000,942.000)
\put(2925.000,762.000){\arc{300.000}{1.5708}{4.7124}}
\path(2295.000,942.000)(2175.000,912.000)(2295.000,882.000)
\put(2175.000,762.000){\arc{300.000}{4.7124}{7.8540}}
\put(2551,3236){\ellipse{900}{900}}
\put(4351,3236){\ellipse{900}{900}}
\path(5400,4512)(5400,12)
\path(300,1512)(5400,1512)
\path(3001,3236)(3901,3236)
\path(302,3235)(1202,3235)
\path(2101,761)(3001,761)
\path(5400,1512)(12600,1512)
\path(6000,3912)(6900,3912)
\path(6780.000,3882.000)(6900.000,3912.000)(6780.000,3942.000)
\path(6001,2561)(6901,2561)
\path(6781.000,2531.000)(6901.000,2561.000)(6781.000,2591.000)
\path(6001,461)(6901,461)
\path(6781.000,431.000)(6901.000,461.000)(6781.000,491.000)
\path(7501,461)(8401,461)
\path(8281.000,431.000)(8401.000,461.000)(8281.000,491.000)
\path(7621.000,3941.000)(7501.000,3911.000)(7621.000,3881.000)
\path(7501,3911)(8401,3911)
\path(7621.000,2591.000)(7501.000,2561.000)(7621.000,2531.000)
\path(7501,2561)(8401,2561)
\path(9121.000,491.000)(9001.000,461.000)(9121.000,431.000)
\path(9001,461)(9901,461)
\path(10621.000,491.000)(10501.000,461.000)(10621.000,431.000)
\path(10501,461)(11401,461)
\path(10320.000,2592.000)(10200.000,2562.000)(10320.000,2532.000)
\path(10200,2562)(11100,2562)
\path(10200,3912)(11100,3912)
\path(10980.000,3882.000)(11100.000,3912.000)(10980.000,3942.000)
\put(0,3162){\makebox(0,0)[lb]{$1$}}
\put(1350,3162){\makebox(0,0)[lb]{$2$}}
\put(2625,3162){\makebox(0,0)[lb]{$1$}}
\put(4050,3162){\makebox(0,0)[lb]{$2$}}
\end{picture} }
        	\hspace{-1.9mm}
	\end{array}
 $$
\caption{Basis for the  invariant spaces corresponding to the 
last 
two rows
and columns of the table of corollary \ref{cor.gl}. Note that 
all
graphs are vertex ordered,  vertex oriented and $2$-colored. 
For 
simplicity, 
the ordering
and the orientation of the vertices is indicated only in the 
northwest part
of the figure.  }\lbl{Basis1}
\end{figure}

\begin{pf}
For $m=1$, there are only two trivalent graphs with $3$ edges 
and 
no 
``decorations''. After including the possible decorations and 
taking 
into
account the colored $AS$ relation of figure \ref{AScl}, we 
arrive 
at
the table of figure \ref{Basis1}.
\end{pf}

Two remarks are in order:
 
\begin{remark}
\lbl{rem.schur}
An alternative way of counting the dimensions of the 
invariant 
spaces
$(\L3 \HQ)^{Sp(\HQ)},
(\UQ)^{Sp(\HQ)}, (\L3 \HQ)^{GL(\VQ)}$ and $(\UQ)^{GL(\VQ)}$
is by decomposing into irreducible representations and using 
Schur's
lemma.
Indeed, since $\L3 \HQ^\ast \simeq \L3 \HQ$, we obtain:
\begin{eqnarray*}
(\otimes^2 \L3 \HQ)^{GL(\VQ)} & = &
( \L3 \HQ^\ast \otimes \L3 \HQ)^{GL(\VQ)} \\
& = & \text{Hom}_{GL(\VQ)}(\L3 \HQ,  \L3 \HQ) 
\end{eqnarray*}
Using lemma \ref{lemma.rep}, we see that
$ \L3 \HQ$ is a sum of six irreducible nonisomorphic 
representations
of $GL(\VQ)$, therefore, it follows by {\em Schur's lemma} 
that
the dimension of  the invariant space $\text{Hom}_{GL(\VQ)}(\L3 
\HQ,  
\L3 
\HQ)$
is six.
We can deduce the rest of the dimensions of the invariant 
spaces
the same way as above. Note however, that this alternative 
way
can only provide us with a dimension count but not with a 
choice 
of
basis for the above mentioned invariant spaces.
\end{remark}

\begin{remark}
On the one hand, weights (or Young diagrams) are a  classical 
and 
convenient way of
parametrizing irreducible representations of classical Lie 
groups. 
On the other hand, trivalent graphs seem to be a very 
convenient 
way
of parametrizing invariant tensors of representations of 
classical
groups.
\end{remark}

\section{Proofs}
\lbl{sec.proofs}

\subsection{Proof of  theorem \ref{thm.0}}
\lbl{sub.prf}

This section is devoted to the proof of theorem \ref{thm.0},
 addenda~\ref{cor.cor}, \ref{cor.corr}, and 
proposition~\ref{cor.mult}. 

\begin{pf}[of theorem \ref{thm.0}]
We fix an admissible surface $f: \Sigma \to M$ and 
consider the (discrete) group  $G=\Tgg$ as  in section  
\ref{sub.gen}. 
 For every nonnegative integer $m$, consider the following 
cocycle
as in corollary \ref{cor.tt}
$$
\f_{2m,2} \in C^{2m} (\Tgg/\Tgg (2);I^{2m} /I^{2m+1})
$$
where $I$ denotes the augmentation ideal of the group ring 
$\BQ \Tgg$.
According to the properties of the Johnson homomorphism 
reviewed in
section \ref{sub.johnson} we have that the abelianization 
(modulo torsion)
of the Torelli group is isomorphic, via the Johnson 
homomorphism, 
with $U$, i.e., that  $ \tau: \Tgg/\Tgg(2)\simeq U$.
Recall from equation \eqref{eq.fun} the linear map
$\Phi^T_f: \BQ \Tgg \to \M $ preserving (by definition) the 
filtration $\{ I^n \}$ of $\BQ \Tgg$ and  
$\FT {\ast}$ of $\M$. In \cite{GL3} it was shown 
that
the filtration $\FT {\ast}$ is $2$-step (following the 
definition 
of
section \ref{sub.gen}), i.e., we have that $\FT {2m-1}=\FT 
{2m}$.
In \cite{GL3} it was additionally  shown that $\FT {2m}=\Fas 
{3m}$,
and that $\Fas {\ast}$ is $3$-step, from which  follows the 
equality
$\GT {2m}=\Gas {3m}$ of the associated graded spaces.
According to the {\em fundamental theorem} of \cite{LMO} and 
\cite{L} it
follows that there is an isomorphism $\Gas {3m} \simeq 
\G_m\A(\phi)$.
Thus we get a linear map $I^{2m} /I^{2m+1} \to \G_m\A(\phi)$, 
and
putting everything together, 
 according to corollary \ref{cor.tt}, we get a $2m$ cocycle  
\begin{equation}
C_{f, m} \in C^{2m}(U, \G_m\A(\phi))
\end{equation} 
Of course, this cocycle depends on the choice of an 
admissible 
surface
$f: \S \hookrightarrow M$, as the notation indicates.
It follows from corollary \ref{cor.step} that $C_{f,m}$ is 
multilinear. In order to show that $C_{f,m}$
is equivariant, consider a diffeomorphism  $h : \S \to \S$.
It follows by definition of the map $\Phi^T_f$ of equation
\eqref{eq.fun} that for $a \in \T(\S)$
we have: $\Phi_{hf}(a)=\Phi_{f}(h^{-1} a h)$. Thus $\Gamma_g$
acts by conjugation on the graded quotients 
$(I \Tgg)^{n}/(I \Tgg)^{n+1}$, and since the action of $\Tgg
\subseteq \Gamma_g$ is trivial, and the quotient $\Gamma_g/\Tgg$
is the symplectic group $Sp(H)$, equivariance follows as stated 
in the
first part of theorem \ref{thm.0}.
Furthermore, it follows by proposition \ref{prop.coh} that
the pullback of the cocycle $C_{f,m}$ to $\Tgg/\Tgg(3)$ (and 
therefore,
to $\Tgg$) represents a trivial cohomology class.
The proof of theorem \ref{thm.0} is complete.
\end{pf}

\begin{pf}[of addendum \ref{cor.cor}]
The commutativity of the diagram follows by definition of the 
map
$C_{f,m}$.
\end{pf}

\begin{pf}[of addendum \ref{cor.corr}]
We explain the statement of this addendum more fully. 
Suppose, 
for each $g$, we choose a Heegaard embedding of a closed 
surface of 
genus $g$ into the $3$-sphere. Fixing a disk, and considering 
surface
diffeomorphisms on the disc complement that pointwise fix a 
neighborhood of the 
boundary,
we obtain a map 
$\Q\Tg\to\M$ by $h\to S^3_{\hat{h}}$, where $\hat{h}$ is the 
obvious extension of 
$h$ to a diffeomorphism of the closed surface. Moreover, with 
respect to an 
inclusion of
such a surface with boundary in another one, we obtain an 
inclusion
$\Tg\sub\T_{g+1,1}$ and we can arrange that the following diagram 
is commutative:

$$
\begin{CD}
\Q\Tg @>>> \M \\
@VVV      @|   \\
\Q\T_{g+1,1} @>>> \M
\end{CD}
$$
This is, in fact, a special case of proposition 
\ref{cor.mult}. We can define $\T =\lim_{g\to\infty}\Tg$ and 
thus combine 
these into a single map

$$\Q\T \to\M$$

Then it is proved in \cite{GL3} that this map sends  
$(I\T)^n$  onto $\FT n$, thus we have epimorphisms 
$(I\T)^n /(I\T)^{n+1}\ \to\GT n$. Thus it follows 
that the stable cocycle $C_{2m}\in C^{2m}(\T / 
\T(2);\G_m\A(\phi))$ 
is onto. From this it is clear, since 
$\G_m\A(\phi)$ is finitely-generated, that $C_{f,m}$ is onto 
for 
large enough $g$.
\end{pf}

\subsection{Proof of  theorem \ref{thm.00}}
\lbl{sub.00}

%%This section is devoted to the proof of theorem \ref{thm.00}.

\begin{pf}%%[of theorem \ref{thm.00}]
For the convenience of the reader, we divide 
the
proof of theorem \ref{thm.00} into several lemmas.
Recall
 first that for a diffeomorphism $\theta$ of a manifold, we 
denote
by $\theta_\ast$ the action of it on the homology of the 
manifold.

\begin{lemma} 
\lbl{lem.hndl}
Let $Q$ be a handlebody and $L=\text{Ker}\{ 
i_{\ast}:H_1 (\partial Q, \BZ)\to H_1 (Q, \BZ)\}$. Suppose 
$\a$ is 
a 
symplectic automorphism of $H_1 (\pa Q, \BZ)$ such that $\a
(L)=L$. Then there exists a diffeomorphism $h$ of $Q$ 
such that $(h|\pa Q)_{\ast}=\a$.
\end{lemma}
\begin{pf}
Let $\a_Q$ be the automorphism of $H_1 (Q,\BZ)\cong H_1 (\pa 
Q,\BZ)/L$ induced by $\a$. Then there exists an orientation
preserving diffeomorphism 
$\bar h$ of $Q$ such that $\bar h_{\ast}=\a_Q$. Now 
consider the symplectic automorphism $\bb =(\bar h |\pa 
Q)\circ\a^{-1}$ of $H_1(\pa Q, \BZ)$. 
If we write $H_1 (Q,\BZ)=L\oplus L'$, where 
$L'$ is a complementary Lagrangian to $L$, then a matrix 
representative of $\bb$ has the form 
$\pmatrix I&C\\0&X 
\endpmatrix$
 Since $\bb$ is symplectic, it follows that $X=I$ and $C$ 
is symmetric. It suffices to see that any such matrix can be 
realized by a diffeomorphism of $Q$. But this is proved in
 \cite{GL3}.
\end{pf}

\begin{lemma}
\lbl{lem.c1}
 Suppose that $f_1,f_2:\S\hookrightarrow M$ are admissible 
Heegaard
surfaces in an \ihs\ $M$ satisfying:
\begin{itemize}
\item $f_1(\S )=f_2(\S )$
%%\item The associated decomposition of $M$ is Heegaard, i.e. 
%%$M_+$ and $M_-$ are handlebodies
\item $(f_{1,\ast})^{-1}(L_{\e})=(f_{2,\ast})^{-1}(L_{\e}) 
\subseteq H_1(\S, \BZ)$
where $L_{\e}$ is the Lagrangian pair in $H_1(f_1(\S), \BZ)$
as in section \ref{sub.history}.
\end{itemize}
Then $C_{m,f_1} =C_{m,f_2}$.
\end{lemma}

\begin{pf} 
Consider the diffeomorphism $g=f_2  f_1^{-1}$ 
of $\f_1(\S)$. Since $g_{\ast}$ preserves $L_{\e}$, for 
$\e =\pm$, we can apply lemma~\ref{lem.hndl} to deduce the 
existence of a diffeomorphism $\hat{h}_{\e}$ of $M_{\e}$ 
which 
induces the same automorphism of $H_1 (f_1(\S), \BZ)$ as $g$. 
In other words we can write $f_2=h_{\e}  f_1  g_{\e}$, 
where $g_{\e}\in\T (\S )$ and $h_{\e}$ is the restriction of 
$\hat{h}_{\e}$
to $\partial M_{\e}$. Recalling the map $\Phi^T_{f_i}$
of equation \eqref{eq.fun}, for $g_i \in \T(\S)$ we have:

\begin{align*}
C_{f_2,k} (g_1 , \ldots ,g_k ) & =\Phi^T_{f_2}((1-g_1 )\ldots 
(1-
g_k ))\\
& =M_{f_2(1-g_1 )\ldots (1-g_k )(f_2)^{-1}}\\
& =M_{h_+f_1 g_+(1-g_1 )\ldots (1-g_k )g_-^{-1}f_1^{-1}h_-^{-
1}}\\
& =M_{f_1g_+(1-g_1 )\ldots (1-g_k )g_-^{-1}f_1^{-1}}\\
&=\Phi^T_{f_1} (g_+ (1-g_1 )\ldots (1-g_k )g_-^{-1})
\end{align*}
But $g(1-g_1 )\ldots (1-g_k )\equiv (1-g_1 )\ldots (1-g_k 
)g\equiv (1-g_1 )\ldots (1-g_k ) \bmod I^{k+1}$ for any 
$g\in\T (\S )$ and so, if $k=2m$, this last term is the 
same as $\Phi^T_{f_1} ((1-g_1 )\ldots (1-g_k ))\in 
\G_{m}\A(\phi)$.
\end{pf}

\begin{lemma}
\lbl{lem.c2}
Suppose that $f_i:\S\hookrightarrow M_i$ for $i=1,2$
are admissible Heegaard surfaces in \ihs s
$M_i$ satisfying:
\begin{itemize}
%%\item The associated decompositions of $M, M'$ are Heegaard 
%%(necessarily of the same genus),
\item $(f_{\ast})^{-1}(L_{1,\e})=(f_{2,\ast})^{-1}(L_{2,\e}) 
\subseteq H_1(\S, \BZ)$, 
where $L_{i,\e}$ for $i=1,2, \e=\pm$ are  the Lagrangian 
pairs 
 associated to $(M_i, f_i)$ respectively.
\end{itemize}
Then $C_{m,f_1} =C_{m,f_2}$.
\end{lemma}

\begin{pf}
We reduce this to the preceding lemma~\ref{lem.c1} by the 
following observation. Let $N$ be an \ihs , 
$f:\S\hookrightarrow N, \S '=f(\S )$. 
Let $h\in\T (\S ')$ and $N'=N_h$. Define 
$f':\S\hookrightarrow N'$ from $f$ by identifying $N_+\sub N$ 
with 
$N_+\sub N_+ \cup_h N_- =N_h$. Then 
$\Phi^T_{f'}(g)=N'_{f'g(f')^{-1}}=N_{hfgf^{-1}}=N_{f(f^{-
1}hf)gf^{-1}}=\Phi^T_f (h'g)$, where $h'=f^{-1}hf\in\T (\S 
)$. Therefore, when $k=2m$: 
\begin{align*}
C_{f',k} (g_1 , \ldots g_k ) &=\Phi^T_{f'}((1-g_1 )\ldots 
(1-g_k ))\\
&=\Phi^T_f (h'(1-g_1 )\ldots (1-g_k ))\\
&\equiv \Phi^T_f (1-g_1 )\ldots (1-g_k ) \bmod \FT {k+1}\\
&=C_{f,k} (g_1 |\ldots |g_k )\in \G_{m}\A(\phi)
\end{align*}
since $h'(1-g_1 )\ldots (1-g_k )\equiv (1-g_1 )\ldots (1-g_k 
)
 \bmod I^{k+1}$.
Thus $C_{f',k} =C_{f,k}$.

Now, since $M_1$ and $M_2$ are endowed with Heegaard 
decompositions of the same genus by $f_1, f_2$, we may 
identify $M_2$ as $(M_1)_h$, for some $h\in\Ga (\S ')$, 
where $\S '=f(\S )$. We may even assume that $h\in\T (\S 
')$ since we can increase the genus of the Heegaard 
decompositions without destroying the hypotheses and every 
\ihs\  has some Heegaard decomposition where the gluing map
is an element of the 
Torelli group. The above observation  and 
lemma~\ref{lem.c1} complete the proof of the present lemma.
\end{pf}

Lemma \ref{lem.c2} implies that the cocycle $C_{f,m}$ in the
case of a Heegaard admissible surface depends only on the 
Lagrangian
pair $(L^+, L^-)$, and will thus be denoted by $C_{\lpm,m}$. In
addition the $Sp(H)$-equivariance of $C_{f,m}$ shows that
$C_{\lpm,m}$ is $GL(L^+_{\BQ})$-invariant. 

We next recall a  useful  lemma due to S. Morita \cite{Mo2}.
Let $\Ng^+$ (respectively, $\Ng^-$)
 be the subgroup of the mapping class group $\Ga_{g,1}$ of
$\Sigma$ that extend to $M^+$ (respectively, $M^-$).
Let $W^+$ (respectively, $W^-$) denote the quotient space:
$ \L3 H/\L3 L^+$ (respectively, $\L3 H/ \L3 L^-$). Note that
we can identify $W^+$ (and similarly, $W^-$) with a subgroup 
 of $\L3 H$ in the following way: Let 
$$\eta_{\pm}:\L3 H\to\L3 (H/L^{\mp})\simeq\L3 L^{\pm} $$
be the projection followed by the natural isomorphism. Then
the projection $\L3 H\to W^{\pm}$ induces an isomorphism  
$W^{\pm}\simeq\text{Ker}\eta_{\pm}$. Alternatively, if we 
choose 
$\{ x_i \}_{i=1}^{g}$ (respectively, $\{ y_i \}_{i=1}^{g}$)
a basis for $L^+$ (respectively, $L^-$) such that 
$\omega(x_i, y_j)=
\delta_{i,j}$, (where $\omega$ is the symplectic form on $H$)
 then $W^+$ is isomorphic to the subgroup of $\L3 H$ 
generated by
all elements of the form: $x_i \w x_j \w y_k, x_i \w y_j \w 
y_k,
y_i \w y_j \w y_k$ for all $1 \leq i < j < k \leq g$.

Recall the Johnson homomorphism $\tau$ of section 
\ref{sub.johnson}.
Then, we have the following lemma:

\begin{lemma}\cite[lemma 4.6]{Mo2}
\lbl{lemma.mor}
With the above notation, identifying $W^{\pm}$ with a 
subgroup of
$\L3 H$, we have the following:
\begin{equation}
\tau(\Tg \cap \Ng^{\pm})= W^{\pm}
\end{equation}
\end{lemma}

We can now show that the cocycle $C_{\lpm,m}$ factors through
a map as in equation \eqref{eq.fct}. 

First notice that, the map $H \to \L3 H$ from section 
\ref{sub.johnson} followed by the projection $\L3 H \to 
\L3 L^{\pm}$ vanishes, and therefore induces   
onto
maps $U= \L3 H/H \to \L3 L^{\pm}$.
%%These maps have the following property: using the 
%%identification
%%$\L3 H \simeq W^{\pm} \oplus \L3 L^{\pm}$, the maps 
%%$\UQ \subseteq \L3 \HQ W^{\pm}_{\BQ} \oplus \L3 
%%L^{\pm}_{\BQ} 
%%\to \L3 L^{\pm}_{\BQ}$ defined above 
%%are given by multiplication by a nonzero rational number
%%on the $\L3 L^{\pm}_{\BQ}$ piece.
   
Let us now consider elements $g_i, h_i \in \Tgg$
(for $i=1,2, \ldots, 2m$)
such that $[\tau(g_1)]=[\tau(h_1)]  \in  \L3 L^+$, and
$[\tau(g_{2m})]=[\tau( h_{2m})] \in \L3 L^-$,
and $g_i=h_i$ for $2 \leq i \leq 2m-1$.
The notation is as follows:
recall that $\tau(g_i), \tau(h_i) \in U$, and temporarily 
denote
both maps $U \to \L3 L^{\pm}$ by $x \to [x]$. 

Now it is clear that $W^{\pm}\cap U=\text{Ker}\{ U\to\L3 
L^{\pm}\}$. 
Since $h_1 g_1^{-1}\in W^+ \cap U$ and $h_{2m}g_{2m}^{-1}\in 
W^- 
\cap U$, we can choose liftings $\tilde{g}_1, \tilde{h}_1, 
\tilde{g}_{2m}, 
\tilde{h}_{2m} 
\in \Tg$ so that  and $\tau(\tilde{h}_{2m} \tilde{g}_{2m}^{-
1}) 
\in
W^-$. Using lemma \ref{lemma.mor} above, there exist $b^+ \in 
\Tg \cap 
\Ng^+$ 
(respectively, $b^-  \in \Tg \cap \Ng^-$) such that
$\tilde{h}_1= b^+ \tilde{g}_1$ (respectively, 
$\tilde{h}_{2m}=  
\tilde{g}_{2m} b^-$).
Since $f:\S \hookrightarrow M$ is a Heegaard embedding, it 
follows
by definition of $\Ng^\pm$ that $M_{b^+ \tilde{h}}=M_{b^+ h}
=M_h = M_{h b^-}= M_{\tilde{h} b^-}$, and 
therefore
that $M_{(1-g_1) \ldots(1- g_{2m})}=M_{(1- b^+ g_1)(1- g_2) 
\ldots
(1 -g_{2m-1})(1- g_{2m} b^-)}$, and therefore that
\begin{equation}
C_{\lpm,m}(g_1, \ldots, g_{2m})=C_{\lpm,m}(h_1, \ldots, 
h_{2m})
\end{equation}
This completes the second part  of theorem 
\ref{thm.00}.

In order to show the third part of theorem \ref{thm.00},
recall first that all \ihs s are oriented. The change of 
orientation 
of
an \ihs\ induces an involution on $\M$, and thus on $\Gas m$ 
for 
every
$m$. Recalling the isomorphism $\Gas m \simeq \G_m\A(\phi)$, 
the
above involution on $\G_m\A(\phi)$ is simply multiplication 
with
$(-1)^m$, \cite[proposition 5.2]{LMO}.
On the other hand, given an admissible Heegaard
 surface $f: \S \hookrightarrow M$
in an \ihs\ $M$, let $\overline{f}: \S \hookrightarrow 
\overline{M}$
denote the same embedding but with different orientation on 
the 
ambient
space $M$.

It is easy to see that the associated change to the set of
Lagrangian pairs is given by $(L^+, L^-) \to  (L^-, 
L^+)$.

Furthermore, note  that for an element $g$ of the Torelli 
group of $\S$, we have the following identity: 
$(\overline{M})_g =\overline{M_{g^{-1}}}$, see also  figure 
\ref{RevOrient}.
Thus, by the above discussion,
we deduce that
 $C_{\overline{f},m}(g_1 ,\ldots ,g_{2m})=(-1)^m 
C_{f,m}(g_{2m}^{-
1},\ldots 
,g_1^{-1})$.
Since passing from $f$ to $\overline{f}$ interchanges the 
Lagrangians, 
and since
$C_{f,m}(g_{2m}^{-1},\ldots ,g_1^{-1})= C_{f,m}(g_{2m},\ldots 
,g_1)$
(due to the fact that $C_{f,m}$ is  multilinear, and the fact 
that
we use multiplicative notation here to denote addition)
this proves the third part on the cocycle level.
The assertion about the cohomology class follows from lemma 
\ref{lemma.involution}, using the fact that 
$[C_{\overline{f},m}]=(-1)^{\binom {2m}2 
+1}\gamma^{\ast}[C_{f,m}]$ and 
$\binom {2m}2 \equiv m \pmod 2$.
This completes the proof  of theorem \ref{thm.00}.
\end{pf}

\begin{figure}[htpb]
$$ \printname{RevOrient}
	\setlength{\unitlength}{0.03\standardunitlength}
	\begin{array}{c}  \hspace{-1.7mm}
        	\raisebox{-8pt}{\begingroup\makeatletter\ifx\SetFigFont\undefined
% extract first six characters in \fmtname
\def\x#1#2#3#4#5#6#7\relax{\def\x{#1#2#3#4#5#6}}%
\expandafter\x\fmtname xxxxxx\relax \def\y{splain}%
\ifx\x\y   % LaTeX or SliTeX?
\gdef\SetFigFont#1#2#3{%
  \ifnum #1<17\tiny\else \ifnum #1<20\small\else
  \ifnum #1<24\normalsize\else \ifnum #1<29\large\else
  \ifnum #1<34\Large\else \ifnum #1<41\LARGE\else
     \huge\fi\fi\fi\fi\fi\fi
  \csname #3\endcsname}%
\else
\gdef\SetFigFont#1#2#3{\begingroup
  \count@#1\relax \ifnum 25<\count@\count@25\fi
  \def\x{\endgroup\@setsize\SetFigFont{#2pt}}%
  \expandafter\x
    \csname \romannumeral\the\count@ pt\expandafter\endcsname
    \csname @\romannumeral\the\count@ pt\endcsname
  \csname #3\endcsname}%
\fi
\fi\endgroup
\begin{picture}(5703,3029)(0,-10)
\thicklines
\put(1125,1507){\ellipse{1200}{3000}}
\put(4125,1507){\ellipse{1200}{3000}}
\path(2025,1507)(3225,1507)
\path(3105.000,1477.000)(3225.000,1507.000)(3105.000,1537.000)
\path(525,1807)(1725,1807)
\path(3525,1207)(4725,1207)
\put(900,757){\makebox(0,0)[lb]{$M^+$}}
\put(900,2257){\makebox(0,0)[lb]{$M^-$}}
\put(3825,607){\makebox(0,0)[lb]{$M^-$}}
\put(3825,2032){\makebox(0,0)[lb]{$M^+$}}
\put(0,1732){\makebox(0,0)[lb]{$g$}}
\put(4950,1132){\makebox(0,0)[lb]{$g^{-1}$}}
\end{picture} }
        	\hspace{-1.9mm}
	\end{array}
 $$
\caption{An orientation reversing of $M$ corresponds to a 
reflection
along the $x$-axis.}\lbl{RevOrient}
\end{figure}
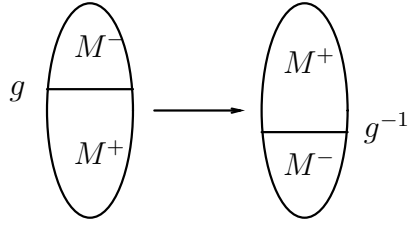

\subsection{Proof of  theorem \ref{thm.formula}}
\lbl{sub.formula}

\begin{pf}%%[of theorem \ref{thm.formula}]

Let $f:\S \hookrightarrow M$ be an admissible Heegaard 
surface, and 
$(L^+,L^-)$  the associated Lagrangian pair of the symplectic 
vector
space $(H,\omega)$ as in section \ref{sub.history}.  
Consider the cocycle $C_{\lpm,m}: \otimes^{2m}U \to 
\G_m\A(\phi)$. Recall from section \ref{sub.gsp} the subgroup of 
the 
symplectic group
$Sp(H)$ isomorphic to GL$(L^+)$. Since this group acts on $H$
preserving the Lagrangian pair $(L^+, L^-)$, it follows from 
theorem 
\ref{thm.00} that $C_{\lpm,m}$ factors through an invariant 
map:
$(\otimes^{2m}U)^{GL(L^+)} \to \G_m\A(\phi)$.
Composing with  the onto map $\G_m\A^{rp,nl,cl} \to 
(\otimes^{2m}U)^{GL(L^+)}$
of proposition \ref{prop.gl},   we get a composite map:
\begin{equation}
\Psi_{\lpm,m}: \G_m\A^{rp,nl,cl} \to \G_m\A(\phi)
\end{equation}
thus finishing the proof of theorem \ref{thm.formula}.
\end{pf}

\subsection{Proof of  theorem~\ref{thm.mor}}
\lbl{sub.mor}

%%This section is devoted to the proof of theorem~\ref{thm.mor}.

\begin{pf}
Let $f:\S \hookrightarrow M$ be an admissible Heegaard 
surface, and 
let
$(L^+,L^-)$ be the associated Lagrangian pair of the 
symplectic 
vector space $(H,\omega)$ as in section \ref{sub.history}.

Let $\l$ denote the Casson invariant, and $W_\l$ its 
associated
manifold weight system. Consider the associated $2$-cocycle
$W_\l \circ C_{\lpm,1}$ of $U$ with
coefficients in $\BQ$ as in corollary \ref{cor.type}. 
According to
theorem \ref{thm.00}, $W_\l \circ C_{\lpm,1}$ factors through 
a
$GL(L^+_{\BQ})$ invariant map:
$\L3 L^+_{\BQ} \otimes \L3 L^-_{\BQ} \to \BQ$. 
According to corollary \ref{cor.gl}, the vector space of such 
invariant maps
is $1$ dimensional, and a nonzero such  map is 
the restriction of the map $C_{\Theta}$ of equation 
\eqref{eq.ctheta}
to $ \L3 L^+_{\BQ} \otimes \L3 L^-_{\BQ} \hookrightarrow 
\otimes^2 \L3 \HQ
\hookrightarrow \otimes^6 \HQ  $, which we will denote by the 
same
name. 
%%Recall the composite map $C^U_{\Theta}$ in the discussion 
before the
%%statement of theorem \ref{thm.mor}. 
%%Using the  $\UQ \hookrightarrow \L3 L_{\BQ} \hookrightarrow
%%\otimes^6 \HQ $ let
%% $C^U_{\Theta}$ denote the restriction of $C_\Theta$ to 
$\otimes^2 \UQ$.
Using the definition of the map $C^U_{\Theta}$ (given before the 
statement
of theorem \ref{thm.mor}) and the above discussion, we deduce 
that 
$W_\l \circ C_{\lpm,1} = c_g C^U_\Theta$ for some rational number 
$c_g$
depending on the genus $g$ of $\S$. According to addendum 
\ref{cor.corr},
 $C_{\lpm,1}$ is stable (with respect to an inclusion of a 
surface in another) and so is $C^U_\Theta$,
therefore $c_g=c$ independent of the genus.

\begin{figure}[htpb]
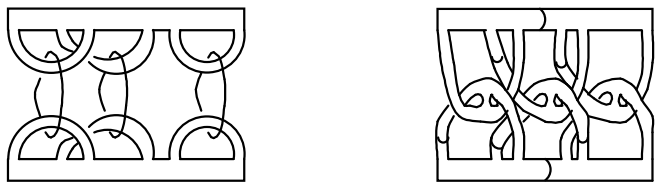

 $$ \printname{2genus16}
	\setlength{\unitlength}{0.03\standardunitlength}
	\begin{array}{c}  \hspace{-1.7mm}
        	\raisebox{-8pt}{\input draws/2genus16.tex }
        	\hspace{-1.9mm}
	\end{array}
 $$
 \caption{On the left, a special case of a $2$-pair blink $L_{bl} 
'$
(bounding the disjoint union of two genus $2$ surfaces)
union a $3$-component \as\ link $L '$. 
 On the right the boundary surface $\S_9$ of a genus
$9$ handlebody. }
 \lbl{2genus16}
\end{figure}

To determine the value of $c$, we need to calculate a particular 
example.
Consider a 2-pair blink $L_{bl} '$ and a 3-component 
algebraically
split link $L '$ in $S^3$ shown in the left part of figure 
\ref{2genus16}.
Part of this figure appeared first in \cite[section 3, figure 
39]{GL3}. 
Note that 
$L'$ is a trivial 3-component link, bounding a disjoint union of
obvious disks.
Choose a unit Seifert-framing for the blink $L_{bl}'$, and for $L 
'$.
 Perform a Dehn twist on each of the three disks
that $L '$ bounds, and let 
$L_{bl}$ denote the image of $L_{bl} '$. After performing the 
twists,
thicken the surface that $L_{bl}$ bounds in order to get two 
disconnected 
genus $3$  solid surfaces, and join them along three tubes to 
form
a genus $9$ surface $\S_9$, see the right part of figure
\ref{2genus16}. We can assume that $L_{bl}$ lies in $\S_9$.
It is easy to see that $\S_9 $ is a genus
$9$ \Heg\ splitting of $S^3$. Let $L^i_{bl}$ (for $i=1,2$)
denote the two pairs of the blink $L_{bl}$; each gives rise  to
an element of the Torelli group of $\S_9$, see \cite{GL3}. Let 
$\alpha^i$  (for $i=1,2$) denote the image in $U$ of each of the 
above mentioned elements of the Torelli group under the Johnson 
homomorphism.
Using the definition of $C^U_{\Theta}$ and the definition of the 
Johnson
homomorphism, it is easy to show that $C^U_{\Theta}(\alpha^1, 
\alpha^2)=-1$.

The first part follows from the following lemma:

\begin{lemma}
\lbl{lemma.Thw}
With the above normalizations, we have the following equalities:
\begin{eqnarray*}
W_{\l} \circ C_{\lpm,1}(\alpha^1, \alpha^2) & = & -2 \\
W_{\l}(\Theta_w) & = & -2
\end{eqnarray*}
\end{lemma}

\begin{pf}
For the first part, note first that for $f: \S_9 \subseteq S^3$
the admissible \Heg\ surface we have: 
\begin{eqnarray*}
\Phi^T_f((1-\aa^1)(1-\aa^2)) & = & S^3 - S^3_{\aa^1}- S^3_{\aa^2} 
+ S^3_{\aa^1
 , \aa^2} \\
& = & S^3 -S^3_{\aa^1 , \aa^2} \\
& = & -[S^3, L_{bl},f] \in \FT 2 = \Fas 3 
\end{eqnarray*}
where the first equality follows from the fact that $S^3_{\aa^1}= 
S^3_{\aa^2}
=S^3$. We also have that:
\begin{equation*}
-[S^3, L_{bl},f] = [S^3, L_{bl} ' \cup L ',f] = \Theta_w =  2 Y_b 
=
2( S^3 - S^3_{\text{Trefoil}, +1})
\end{equation*}
where the first equality follows from the fact that $L '$ is an 
unlink,
and the second follows from the fact that surgery on each of the 
pairs
of $L_{bl} '$ corresponds to the alternating sum of cutting or 
not
 each vertex of the $\Theta$ graph, see also  
\cite[section 3, figure 36]{GL3}. The third equality follows
from the main identities of \cite{Oh}, \cite{GL2}, and the last
by the fact that $S^3_{\text{Trefoil}, +1}$ equals the result of 
Dehn
surgery on a Borromean ring of three components with framing 
$+1$.
Note that Trefoil is a right (or left)handed trefoil in $S^3$ 
depending on the  
Borromean ring; either case does not affect the 
validity of our calculation.

 This and the definition of $W_{\l} \circ C_{\lpm,1}$ imply that
\begin{eqnarray*}
W_{\l} \circ C_{\lpm,1}(\alpha^1, \alpha^2) & = & 
\l(\Phi^T_f((1-\aa^1)(1-\aa^2))) \\
& = & \l(\Theta_w) \\
& = & 2 \l(S^3) - 2\l( S^3_{\text{Trefoil}, +1})
\end{eqnarray*}
Using the normalizations of the Casson invariant it follows  that 
$\l(S^3)=0$ and $\l(S^3_{\text{Trefoil}, +1})=1$
which proves  the lemma.
\end{pf}

In order to show the second part, since $\G_1\A(\phi)$ is $1$ 
dimensional,
the first part implies that: 
\begin{equation}
\lbl{eq.wl}
C_{\lpm,1} = c^{'}_g C^U_\Theta \cdot
\Theta_w
\end{equation}
 for some rational number $c^{'}_g$. The
stability of $C_{\lpm,1}$ implies that $c^{'}_g=c^{'}$ 
independent 
of the
genus. Composing equation \eqref{eq.wl} with $W_\l: 
\G_1\A(\phi)\to 
\BQ$
we obtain that $ W_\l \circ C_{\lpm,1} = c^{'} W_\l(\Theta_w) 
C^U_\Theta$, 
which (due to the first part of theorem \ref{thm.mor}) 
implies that $ 2= c^{'} W_\l(\Theta_w)$.
Using lemma  \ref{lemma.Thw}, the second part of 
theorem \ref{thm.mor} follows.

In order to show the third part of theorem \ref{thm.mor},
recall first from theorem \ref{thm.00} that $C_{\lpm,1} : 
\otimes^2 U \to
\G_1\A(\phi)$ factors through a $GL(L^+)$ invariant map 
$ \L3 L^+_{\BQ} \otimes \L3 L^-_{\BQ} \to \G_1\A(\phi)$.
Thus using the basis of the four dimensional vector
space  $\G_1\A^{rp,nl,cl}(\phi)$ shown in the southeast part of 
figure \ref{Basis1} and the definition of $\Psi_{\lpm,1}$ and 
corollary
\ref{cor.gl}, the second part of theorem \ref{thm.mor} implies 
the
third part.

In order to show the fourth part of theorem \ref{thm.mor},  
namely
that $C_{\lpm,1}$ represents a non-zero cohomology 
class, we interpret it as a cup-product. Consider the 
elements 
$\xi^{\pm}\in H^1 (\L3 H;\L3 L^{\pm})$ defined by the 
homomorphisms $\L3 H\to\L3 (H/L^{\mp})\cong\L3 L^{\pm}$. Then 
the 
intersection pairing on $H$ defines a non-singular pairing on 
$\L3 H$ which induces a non-singular pairing $\L3 L^+ 
\otimes\L3 
L^-\to\Q$. Using this pairing we have a cup-product 
$$ H^1 (\L3 H;\L3 L^+)\otimes H^1 (\L3 H;\L3 L^-)\to H^2 (\L3 
H;\Q )
$$
Now it is straightforward to check, from 
equation~\eqref{eq.mor}, 
that $j^{\ast}[C_{\lpm,1}]=\xi^+ \cup\xi^-$, where $j:\L3 
H\to U$ is 
the projection.

Recall that, for any finitely-generated abelian group $A$, 
there 
is an isomorphism of $\Q$-algebras:
$$ H^{\ast} (A;\Q )\cong \Lambda^{\ast}(A^{\ast})\otimes\Q
$$
where $A^{\ast}$ is the dual of $A$. The product structure on 
$\Lambda^{\ast}(A^{\ast})$ is defined as follows:
$$ \phi^p \cdot\psi^q (v_1\wedge\ldots\wedge v_{p+q})= 
\sum_{\pi}\operatorname{sgn}(\pi )\phi^p (v_{\pi 
1}\wedge\ldots\wedge 
v_{\pi p})\psi^q (v_{\pi (p+1)}\wedge\ldots\wedge v_{\pi 
(p+q)})
$$
where the sum ranges over all {\em shuffle } permutations 
$\pi$.

Thus $\xi^+ \cup\xi^-\in H^2 (\L3 H;\Q )\cong\Lambda^2 (\L3 
H^{\ast})$ is defined by 
$$u\wedge v\mapsto \omega(\xi^+ (u), \xi^- (v))-\omega(\xi^+ 
(v),\xi^- (u))
$$
for $u,v\in\L3 H$. To see this is non-trivial we note, for 
example, 
that, if 
$u\in\L3 L^+$ and $v\in\L3 L^-$,\ then $\xi^+ \cup\xi^- 
(u\wedge v)=\omega(u,v)$.

Now suppose that $K^{\pm}$ is another Lagrangian pair with 
associated classes $\eta^{\pm}\in H^1 (\L3 H;K^{\pm})$, and 
suppose 
that $\xi^+ \cup\xi^- =\eta^+ \cup\eta^-$. We first point out 
that 
$\L3 L^+ +\L3 L^- =\L3 K^+ +\L3 K^- \sub\L3 H$. This follows 
from 
the observation that $\xi^+ \cup\xi^- (u\wedge w)=0$ for all 
$v\in\L3 H$ if and only if $\omega(u, \L3 L^+ +\L3 L^- )=0$.

We can assume that $\dim H\geq 6$. Let $p:K^+ \to L^+$ be the 
restriction of the projection of $H$ onto $L^+$ with kernel 
$L^-$. 
Choose any basis $v_1 ,\ldots ,v_n$ of $K^+$ such that
$p(v_i )=0$ for $i\leq r$, and $p(v_i )$ are linearly 
independent in 
$L^+$ for $i>r$. Write $v_i =v_i^+ +v_i^-$ for $i>r$, where 
$v_i^{\pm}\in L^{\pm}$. We consider several cases.

$\bold{1<r<n}$. Then $v_1 ,v_2 \in L^-$ and $v_1 \wedge v_2 
\wedge 
v_n =v_1 \wedge v_2 \wedge v_n^+ +v_1 \wedge v_2 \wedge v_n^-
$. If 
we now consider the direct sum decomposition:
\begin{equation}
\lbl{eq.dec}
\L3 H=(\L3 L^+ \oplus\L3 L^- )\oplus (L^+ \otimes\Lambda^2 
L^- )\oplus 
(L^- \otimes\Lambda^2 L^+ )
\end{equation}
then we see that the component of $v_1 \wedge v_2 \wedge v_n$ 
in 
$L^+ \otimes\Lambda^2 L^- $ is non-zero. But $v_1 \wedge v_2 
\wedge v_n 
\in\L3 K^+ \sub\L3 K^+ +\L3 K^- =\L3 L^+ +\L3 L^- $, which 
means this 
component should be zero.

$\bold{r=1}$. For any $1<i<j$, we examine $v_1 \wedge v_i 
\wedge 
v_j$. The component in $L^+ \otimes\Lambda^2 L^- $ is $v_i^+ 
\wedge (v_j^- 
\wedge v_1 )+v_j^+ \wedge (v_1 \wedge v_i^- )$. Since $v_i^+ 
,v_j^+$ 
are linearly independent, this means $v_1 \wedge v_i^- =v_1 
\wedge 
v_j^- =0$. Thus $v_i^- =c_i v_1$, for some scalars $c_i$, 
$i>1$. But 
then $v_1 \wedge v_i \wedge v_j =v_1 \wedge v_i^+ \wedge 
v_j^+$ 
lying in $L^- \otimes\Lambda^2 L^+ $ in the decomposition of 
equation~\eqref{eq.dec}.

$\bold{r=0}$. Suppose $K^+ \not= L^+$. Then we can assume 
without loss
of generality that 
$v_1^- \not= 0$. If we consider the element $v_1 \wedge v_i 
\wedge 
v_j$, then the component in $L^+ \otimes\Lambda^2 L^- $ is 
$v_1^+ \wedge 
(v_i^- \wedge v_j^- )+v_i^+ \wedge (v_j^- \wedge v_1^- 
)+v_j^+ 
\wedge (v_1^- \wedge v_i^- )$. Since  $v_1^+ ,v_i^+ ,v_j^+$ 
are 
linearly independent, we conclude that $v_1^- \wedge v_i^- 
=v_1^- 
\wedge v_j^- =0$ and so $v_i^- =c_i v_1^-$, for suitable 
scalars 
$c_i$. Now let us replace each $v_i$ by $v_i -c_i v_1$. In 
other 
words we can assume that each $v_i \in L^+$ for $i>1$. But 
now we 
can see that $v_1 \wedge v_i \wedge v_j$ has component $v_1^- 
\wedge 
v_i \wedge v_j$ in $L^- \otimes\Lambda^2 L^+ $ and so $v_1^-$ 
would have 
to be zero.

The conclusion from these arguments is that either $r=n$, in 
which 
case $K^+ =L^-$, or, from the last case, that $K^+ =L^+$. 
Similarly,
we see that $K^- =L^+$ or $L^-$.
 Finally we need to check that 
$\xi^+ \cup\xi^- =\eta^+ \cup\eta^-$ in case $K^+ =L^-$ and 
$K^- =L^+$. Using the orthogonal direct sum decomposition 
\eqref{eq.dec}, we can write $u=u^+ + u^- + u^0, v=v^+ +v^- 
+v^0$, where $u^+ ,v^+ \in\L3 L^+ , u^- ,v^- \in\L3 L^-$ and 
$u^0 ,v^0 \in (L^+ \otimes\Lambda^2 L^- )\oplus 
(L^- \otimes\Lambda^2 L^+ )$. Then we have 
\begin{align*} 
\xi^+ \cup\xi^- (u\wedge v) &=\omega(u^+ , v^-) -\omega(v^+ , 
u^-) 
\\
\eta^+ \cup\eta^- (u\wedge v)&=\omega(u^- , v^+) -\omega(v^- 
, 
u^+)
\end{align*}
The skew-symmetry of the symplectic pairing implies that 
these are equal. 
\end{pf}

\subsection{Proof of proposition \ref{cor.mult}}
\lbl{sub.cor.mult}

In this section we give the proof of proposition 
\ref{cor.mult}.

\begin{pf}[of proposition \ref{cor.mult}]
Under the assumptions of corollary \ref{cor.mult}, we are 
given an
embedded sphere $S \hookrightarrow M$ in an \ihs\ $M$ 
which 
separates
$M$ into two components. Therefore, we have that $M$ is
a connected sum of two \ihs s $M_1, M_2$, along the 
separating sphere $S$,
i.e., $M=M_1 \sharp M_2$. Furthermore, by assumption,
the admissible surfaces $f_i:\S_{g_i}
\hookrightarrow M$ belong to different components of $M-S$. 
Recall
the composite surface $f_1 \cup f_2: \S_{g_1 + g_2}=
\S_{g_1} \cup_{\partial S-(D_1 \cup D_2)}
\S_{g_2} \hookrightarrow M$. Therefore,
for $h_i \in \T(\S_{g_i,1})$ ($i=1,2$) we get an element $h_1 
\cup 
h_2
\in \T(\S_{g_1 + g_2})$, and an isomorphism:
\begin{equation*}
M_{h_1 \cup h_2} \simeq (M_1)_{h_1} \sharp (M_2)_{h_2} 
\end{equation*}
Since $\M$ (respectively, $\A(\phi)$)
is a commutative algebra with multiplication the operation
of connected sum on \ihs s, (respectively, the disjoint union 
of 
vertex 
oriented trivalent graphs) and since the isomorphism
 $\Gas {3m} \simeq \G_m\A(\phi)$ preserves the algebra 
structures,
the result follows.
\end{pf}

\subsection{Proof of theorem \ref{thm.homo}}
\lbl{sub.homo}

In this section we give a proof of theorem \ref{thm.homo}.
Since the proof  combines several rather different 
techniques, for the convenience of the reader we separate it into 
several
steps.

\begin{pf}[of theorem \ref{thm.homo}]

\begin{itemize}
\item{{\bf Step 1}} 
\hspace{0.5cm}
The definition of the map $D_{f,m}$ of equation \eqref{eq.ghom}.
\end{itemize}

%%For a group $G$ and a positive integer $m$
%% recall first the map $\G_{2m}G   \otimes \BQ \to 
%%(IG)^{2m}/(IG)^{2m+1}$
%%from equation \eqref{eq.igr}. According to lemma
%%\ref{lemma.ci} this is a linear map.
%%
%%Given an admissible genus $g$ surface $\S \hookrightarrow M$,
%%consider  the above map for $G=\Tgg$.
%%Composing with the induced graded map $\PTf:
%%(I\Tgg)^{2m}/(I\Tgg)^{2m+1} \to \GT {2m} $ from equation 
%%\eqref{eq.fun}
%%together with the isomorphism $\GT {2m} \simeq \G_m\A(\phi)$
%%we get an induced linear map: 
%%\begin{equation}
%%D_{f,m}: \G_{2m} \Tgg \otimes \BQ \to \G_m\A(\phi)
%%\end{equation}
%%
%%In order to finish the first step, we need to show that
%%the image of $D_{f,m}$ is contained in the subspace 
%%$\G_m\A^{conn}(\phi)$
%%of $\G_m\A(\phi)$ spanned by linear combinations of connected 
%%graphs.
%%This follows by corollary \ref{cor.prim}.
Just set $D_{f,m}=\phi_f^{\Tgg}$ as defined in  corollary \ref{cor.prim}.

\begin{itemize}
\item{{\bf Step 2}} 
\hspace{0.5cm}
$D_{f,m}$ is determined by $C_{f,m}$.
\end{itemize}

Indeed,  equation \eqref{eq.cd} follows from lemma \ref{lemma.ci} 
and the above definition of  $D_{f,m}$. It remains to show that
$\G_{2m} \Tgg \otimes \BQ$ is spanned by elements of the form
$[x_1, [x_2, \ldots ,[x_{2m-1}, x_{2m}]]]$ for $x_i \in \Tgg$.
This follows from several applications of the Jacobi identity:
$[[a,b],c]=[a,[b,c]] -[b,[a,c]]$ for $a \in \Tgg(n_1), b \in 
\Tgg(n_2),
c \in \Tgg(n_3)$ with $n_1 + n_2 + n_3 = 2m$.

In case $f: \S \hookrightarrow M$ is an admissible Heegaard 
surface,
$C_{f,m}$, and thus 
$D_{f,m}$, depends only on the associated Lagrangian pair 
$(L^+,L^-)$.
In that case we will denote $D_{f,m}$ by $D_{\lpm,m}$.
Assume from now on that $f$ is an admissible \Heg\ surface.

\begin{itemize}
\item{{\bf Step 3}} 
\hspace{0.5cm}
$D_{\lpm,m}$ satisfies the symmetry property of equation 
\eqref{eq.symm}.
\end{itemize}

Indeed, figure \ref{RevOrient} shows that
$D_{f,m}(a^{-1})=D_{\overline{f},m}(a)$ where $\overline{f}$ is
the surface in the {\em orientation reversed} 3-manifold 
$\overline{M}$.
Since the involution of reversing the orientation in the ambient 
3-manifold
is multiplication by $(-1)^m$ on $\G_m\A(\phi)$, the result 
follows.
%%Equation \eqref{eq.symm} and the fact that $D_{\lpm,m}$ is a
%%linear map imply that $D_{\lpm,m}=0$ for $m=\text{ even.}$

\begin{itemize}
\item{{\bf Step 4}} 
\hspace{0.5cm} Assume now that $f$ is the standard genus $g$ 
\Heg\ 
splitting of $S^3$. Then, for $g \geq 5m+1$, $D_{\lpm,m}$ is 
onto.
\end{itemize}

This follows by corollary \ref{cor.L} of theorem \ref{thm.L}
whose proof is given below. The proof of theorem \ref{thm.L} and
corollary \ref{cor.L} given below is long and technical; 
furthermore it is
logically independent from the rest of the proof of theorem  
\ref{thm.homo}.

\begin{itemize}
\item{{\bf Step 5}} 
\hspace{0.5cm}
The case of $m=1$.
\end{itemize}

We now describe explicitly the map $D_{\lpm,1}$. Assume that
we are given an admissible \Heg\ genus $g$ surface $f$. Recall 
first that
$\G_m \Tgg \otimes \BQ$ is a  finite dimensional, stable with 
respect to the genus, representation of $Sp(H)$. It follows by a 
theorem
of Quillen \cite{Qu} (see also \cite{Hain1}) that it is a 
rational
representation of $Sp(\HQ)$.  It is a {\em  
very interesting} question to analyze the structure of the above
representation. Motivated by the above question Morita 
\cite{Mojo}
developed a theory of {\em higher Johnson homomorphisms}, known 
to form
a Lie algebra. The structure of this and related Lie algebras 
have been
analyzed in the pioneering work of Morita \cite{Mo1} and Hain 
\cite{Hain}.
In case of $\G_2 \Tgg \otimes \BQ$ the answer is known and we 
describe it
here.

It turns out that for $ g \geq 6$ we have the following 
decomposition
as representations of $Sp(\HQ)$ \cite{Mo1}, \cite{Hain}:
\begin{equation}
\lbl{eq.J3}
\G_2 \Tgg \otimes \BQ = V(0) + V(2 \e_2)
\end{equation}

Recall from section \ref{sub.johnson} that $\Tgg(2)=\Kgg$, thus
we have that $\G_2 \Tgg= \Kgg/\Tgg(3)$.
Morita \cite[section 5]{Mo1} using his theory of {\em secondary 
characteristic classes} defined a  group homomorphism 
$d_1: \Kgg \to \BQ$ which vanishes on $\Tgg(3)$, thus inducing a 
map
$\G_2\Tgg \otimes \BQ \to \BQ$. Furthermore, Morita \cite[section 
1]{Mo1}
defined a higher version of Johnson's homomorphism
 $ \tau_3: \Kgg \to V(2 \e_2)$, which also vanishes
on $\Tgg(3)$ thus inducing a map $\G_2 \Tgg \otimes \BQ \to V(2 
\e_2)$.
The above maps are $Sp(\HQ)$ invariant, and stable, and realize 
the
decomposition of \eqref{eq.J3} as a $Sp(\HQ)$ module.
Moreover, Morita \cite{Mo1} using a \Heg\ splitting $f$, defined 
a 
map $q_f : V(2 \e_2) \to \BQ$ and showed in his main result 
\cite[theorem 6.1]{Mo1} that:
\begin{equation}
W_{\l} \circ D_{\lpm,1}= - \frac{1}{24} d_1 - q_f
\end{equation}
Note that the change in sign from \cite[theorem 6.1]{Mo1} to the 
above
equation is due to the fact that Morita uses the map $\Tgg \to 
\BQ$
to be $a \to \l(S^3_a ) - \l(S^3)$; however  we use the map
$\Tgg \to \BQ$ to be $ a \to  \l(S^3)- \l(S^3_a ) $.
From this, it follows immediately that 
$W_{\l} \circ D_1= - \frac{1}{24} d_1$, thus
finishing step 5 and the proof of theorem \ref{thm.homo}.
\end{pf}

 The proof of theorem \ref{thm.L} and corollary \ref{cor.L} 
occupies the 
rest of this section. The proof is long and technical, and 
consists of 
combinatorial
as well as geometric topology arguments. We urge the reader to
keep in mind the figures.

Let $f:\S_g \hookrightarrow S^3$ be the standard genus $g$ \Heg\ 
splitting of $S^3$,
which we keep fixed for the rest of this section. We follow the 
notation
and terminology of section \ref{sub.fti}.  For $L$
an $f$-compatible Lagrangian  consider the maps
$\G\PLf : \G\BQ\Lgg \to \G^L\M$ and 
$\phi^L_f : \G\Lgg \otimes \BQ \to \G\A^{conn}(\phi)$,
which, for simplicity, we denote  by $\Phi^L$ and $\phi^L$
respectively.
Then, we have the following theorem:

\begin{theorem}
\lbl{thm.L}
Suppose $g\geq 5n+1$. Then there exists an $f$-compatible 
Lagrangian 
$L\sub H_1 (\S_g )$ so that $\fl (\G_{3n}\Lgg \otimes\Q ) 
=\Ac n$.
\end{theorem}

\begin{pf}

Due to the length of the proof, for the convenience of the reader
we provide the proof in six (or perhaps seven) steps.

\begin{itemize}
\item{{\bf Step 0}} 
\hspace{0.5cm}
A non proof.
\end{itemize}

We first give a ``too good to be true'' proof. $\G\PLf$ is a map 
of coalgebras
and according to the results of \cite{GL3} reviewed in section 
\ref{sub.fti},
$\G\PLf$ is stably onto. If it {\em were} the case that $\G\PLf$ 
was a Hopf
algebra map, the induced map on the primitives would be onto 
by a dimension count using the Poincar\' e-Birkhoff-Witt theorem.
Unfortunately,
$\G\PLf$ does not preserve the product structure,  see
remark \ref{rem.LKT}.

\begin{itemize}
\item{{\bf Step 1}} 
\hspace{0.5cm}
A reduction to ``chord diagrams''.
\end{itemize}

We begin by recalling  the following definition: a degree $n$
chord diagram (on a cirle embedded in the standard way in the plane) 
is a
collection of
$n$  chords with 
$2n$ distinct end points; 
for an example see figure \ref{A}. Note that a chord diagram can
be  thought of as
a connected vertex oriented trivalent graph (with the 
counterclockwise orientation at each vertex), and 
furthermore this way
a degree $n$ chord diagram gives rise to an element of $\Ac n$.
Note that a  degree $n$ chord diagram has $2n$ {\em external} 
edges 
(the ones on the circle) and $n$ {\em internal} ones. 
Degree $n$ chord diagrams  are rather special elements of $\Ac 
n$; however 
we have the following:

\begin{claim}
\lbl{claim.semi}
$\Ac n$ is generated by chord diagrams as above.
\end{claim}

A proof, using the $IHX$ relation,  can be found in \cite[lemma 
3.2]{GL2}.

With the notation of section \ref{sub.fti},
%using  lemma \ref{lemma.wb} and the above claim it follows that
%we can reformulate theorem \ref{thm.L} in the following way:
we make the following:

\begin{claim}
\lbl{claim.eL}
For every  $n > 1$
there is an $f$-compatible Lagrangian $L$ with the following 
property:
For each degree $n$ chord diagram $\Ga$, there is an element
$\xi^{\Ga} \in \Lg(3n)$ such that: 
\begin{equation}
\lbl{eq.LL}
F(\phi^L(\xi^{\Ga})) \equiv \Ga_b \bmod \Yb
\end{equation}
\end{claim}

Using lemma \ref{lemma.wb} the above claim implies that
for each degree $n > 1$ chord diagram $\Ga$, there is
an element $\xi^{\Ga} \in \Lg(3n)$ such that 
$\phi^L(\xi^{\Ga}) = \Ga_w$. This,   together with claim \ref{claim.semi}
implies theorem \ref{thm.L} for $ n > 1$.

Thus, we will prove theorem
\ref{thm.L} by proving the special case of theorem \ref{thm.L} for
$n=1$, and proving  claim \ref{claim.eL} for $n  > 1$. Note that
theorem \ref{thm.L} is obvious for $n=0$.
In the rest of the proof, we will be working in the graded space 
$\G_{3n}\Ab$
which is isomorphic to $\Gas {3n}$.

\begin{figure}[htpb]
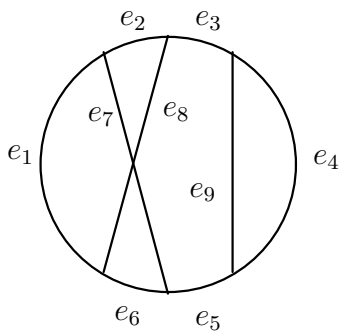

$$ \printname{A}
	\setlength{\unitlength}{0.03\standardunitlength}
	\begin{array}{c}  \hspace{-1.7mm}
        	\raisebox{-8pt}{\input draws/A.tex }
        	\hspace{-1.9mm}
	\end{array}
 $$
\caption{A degree $3$ chord diagram on a circle.}\lbl{A}
\end{figure}
 
\begin{itemize}
\item{{\bf Step 2}} 
\hspace{0.5cm}
An important construction.
\end{itemize}

We begin with a construction which will be an important component 
of, as 
well as a warm-up for, the general construction. Let $\S_2 =\pa 
H$ be 
the surface of genus two in $S^3$, where $H=T_1 \coprod T_2$ is 
the 
boundary 
connected sum of two solid tori in $\R^3$. We distinguish three 
simple 
closed curves on $\S_2$. Let $a$ be a meridian in $\pa T_1$, i.e. 
the boundary of a meridian disk in $T_1$, and $c$ be a longitude 
in 
$\pa T_2$. Finally let $b$ be a band sum of two disjoint 
meridians 
in $\pa T_2$, where the band passes once, longitudinally, around 
$\pa 
T_1$. See figure \ref{B}. The orientations of these curves can be 
chosen 
arbitrarily for the moment. Let $\aa ,\bb .\gg$ be the 
diffeomorphisms 
of $\S_2$ defined by Dehn twists along $a,b,c$, respectively. 
Note 
that $\gg\in (\kg )_2$, since $b$ bounds in $\S_2$, but $\aa$ and 
$\bb$ are not even in the Torelli group. But there is a 
Lagrangian $L\sub H_1 (\S_2 )$, compatible  with the given 
embedding
 of $\S_2$ in $S^3$, so that $\aa ,\bb\in\LL^L_2$. In fact 
we can just take $L$ to be the subgroup generated by the homology 
classes of $a$ and $c$.

\begin{figure}[htpb]
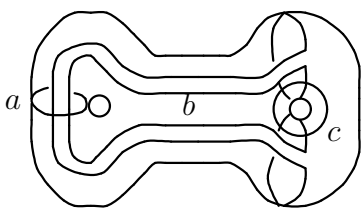

$$ \printname{B}
	\setlength{\unitlength}{0.03\standardunitlength}
	\begin{array}{c}  \hspace{-1.7mm}
        	\raisebox{-8pt}{\input draws/B.tex }
        	\hspace{-1.9mm}
	\end{array}
 $$
\caption{A genus $2$ handlebody $\Sigma_2$ together with $3$ 
curves $a,b,c$
on it.}\lbl{B}
\end{figure}

Consider the element $[[\aa ,\bb ],\gg ]\in (\LL^L_2 )_3$. 
Then, we have the following:

\begin{itemize}
\item{{\bf Step 3}} 
\hspace{0.5cm} $\fl [[\aa ,\bb ],\gg ]$ is a generator of $\Gas 
3\iso\Af 
3=\Ac 1\iso\Q$. Thus theorem \ref{thm.L} holds for $n=1$.
\end{itemize}

\begin{pf}
Consider the element $(1-\aa)(1-\bb)(1-\gg)\in I= 
(I\LL^L_2 )^3$. Then $\hfl ((1-\aa)(1-\bb)(1-\gg))$  represents 
the same
linear combination in $\M$  as $[S^3 ,K]$, where the 3-component 
link 
$K$ is constructed as follows. Take three concentric copies of 
$\S_2$ in $\R^3$ and place $a$ in the outer copy, $b$ in the 
middle copy and $c$ in the inner copy. Then $K$ is given by these 
three disjoint curves in $\R^3$. See figure \ref{C}.
 We refer the reader 
to \cite{GL3} for the 
explanation of this. 

\begin{figure}[htpb]
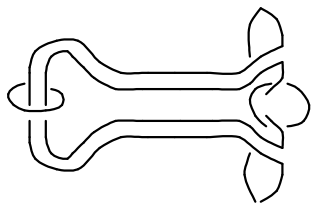

$$ \printname{C}
	\setlength{\unitlength}{0.03\standardunitlength}
	\begin{array}{c}  \hspace{-1.7mm}
        	\raisebox{-8pt}{\input draws/C.tex }
        	\hspace{-1.9mm}
	\end{array}
 $$
\caption{After taking $3$ concentric copies of the surface of 
figure \ref{B},
and placing $a, b, c$ on each copy in that order, we arrive at 
the 
$3$-component Borromean link $K$ shown in the figure above.  
}\lbl{C}
\end{figure}

Note that $K$ is just the Borromean rings and so
 represents a generator of $\Gas 3$. The following claim 
\ref{claim.bor}
completes the proof of this step.

\begin{claim}
\lbl{claim.bor}
We have:
$$ \hfl (1-[[\aa 
,\bb ],\gg ])=\hfl ((1-\aa)(1-\bb)(1-\gg))$$
\end{claim}

 To prove the claim recall the following identity from lemma 
\ref{lemma.ci}:
\begin{eqnarray*}
1-[[\aa ,\bb ],\gg ] & \equiv & (1-\aa)(1-\bb)(1-\gg)-(1-\bb)(1-
\aa)(1-\gg) \\
& & -(1-\gg)(1-\aa 
)(1-\bb)+(1-\gg)(1-\bb)(1-\aa)\mod I^4
\end{eqnarray*}

Now for any permutation $\aa_1 ,\aa_2 ,\aa_3$ of $\aa ,\bb ,\gg$, 
$\hfl (\aa_1 -1)(\aa_2 -1)(\aa_3 -1)=[S^3 ,K']$, where $K'$ is 
the 
3-component link  defined by placing $a,b,c$ in 
concentric copies of $\S_2$, just as above, except that they are 
placed in the permuted order. But it is easy to see that, for any 
{\em 
non-trivial }permutation, the resulting link $K'$ is trivial. See
figure \ref{D}. Thus, 
from the above formula we see that $\hfl (1-[[\aa ,\bb ],\gg ])=
\hfl ((1-\aa)(1-\bb)(1-\gg))$. This concludes the proof of the 
claim and
of step 3.  \end{pf}

\begin{figure}[htpb]
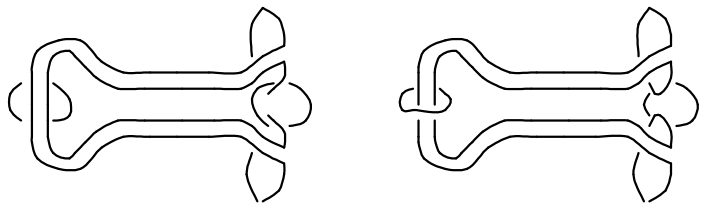

$$ \printname{D}
	\setlength{\unitlength}{0.03\standardunitlength}
	\begin{array}{c}  \hspace{-1.7mm}
        	\raisebox{-8pt}{\input draws/D.tex }
        	\hspace{-1.9mm}
	\end{array}
 $$
\caption{The links shown above associated to any nontrivial 
permutation of 
$a,b,c$ are trivial.}\lbl{D}
\end{figure}

\begin{itemize}
\item{{\bf Step 4}} 
\hspace{0.5cm} 
The definition of the $f$-admissible Lagrangian $L$.
\end{itemize}

Now let $\Ga$ be any chord diagram with $2n$ vertices. We 
associate 
to $\Ga$ a Heegaard surface $\S_{\Ga}\sub\R^3$ as follows. Choose 
an 
embedding of $\Ga\sub\R^3$. Then, at every vertex $v$ of $\Ga$, 
place a copy $\S (v)$ of $\S_2$ and, for every edge $e$ of $\Ga$ 
with vertices $v_e ,v_e '$, take a connected sum of $\S (v_e )$ 
with 
$\S (v_e ')$ using a {\em tube} $T(e)$ running along the edge 
$e$. 
See figure \ref{E}.

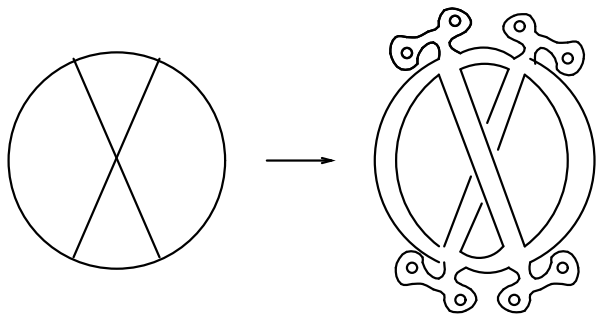
\begin{figure}[htpb]
$$ \printname{E}
	\setlength{\unitlength}{0.03\standardunitlength}
	\begin{array}{c}  \hspace{-1.7mm}
        	\raisebox{-8pt}{\begingroup\makeatletter\ifx\SetFigFont\undefined
% extract first six characters in \fmtname
\def\x#1#2#3#4#5#6#7\relax{\def\x{#1#2#3#4#5#6}}%
\expandafter\x\fmtname xxxxxx\relax \def\y{splain}%
\ifx\x\y   % LaTeX or SliTeX?
\gdef\SetFigFont#1#2#3{%
  \ifnum #1<17\tiny\else \ifnum #1<20\small\else
  \ifnum #1<24\normalsize\else \ifnum #1<29\large\else
  \ifnum #1<34\Large\else \ifnum #1<41\LARGE\else
     \huge\fi\fi\fi\fi\fi\fi
  \csname #3\endcsname}%
\else
\gdef\SetFigFont#1#2#3{\begingroup
  \count@#1\relax \ifnum 25<\count@\count@25\fi
  \def\x{\endgroup\@setsize\SetFigFont{#2pt}}%
  \expandafter\x
    \csname \romannumeral\the\count@ pt\expandafter\endcsname
    \csname @\romannumeral\the\count@ pt\endcsname
  \csname #3\endcsname}%
\fi
\fi\endgroup
\begin{picture}(8207,4279)(0,-10)
\thicklines
\put(5573,3632){\ellipse{150}{150}}
\put(6248,4082){\ellipse{150}{150}}
\path(6023,3557)	(6029.531,3612.126)
	(6031.707,3653.158)
	(6023.000,3707.000)

\path(6023,3707)	(5994.417,3752.296)
	(5948.000,3782.000)

\path(5948,3782)	(5905.671,3778.310)
	(5863.798,3759.226)
	(5798.000,3707.000)

\path(5798,3707)	(5752.618,3635.461)
	(5734.390,3592.870)
	(5723.000,3557.000)

\path(5723,3557)	(5726.885,3520.430)
	(5723.000,3482.000)

\path(5723,3482)	(5679.191,3439.674)
	(5621.919,3410.739)
	(5558.937,3398.684)
	(5498.000,3407.000)

\path(5498,3407)	(5439.228,3444.076)
	(5394.136,3502.194)
	(5363.477,3568.965)
	(5348.000,3632.000)

\path(5348,3632)	(5346.982,3688.832)
	(5358.226,3752.270)
	(5383.107,3811.823)
	(5423.000,3857.000)

\path(5423,3857)	(5458.655,3869.268)
	(5500.318,3867.905)
	(5573.000,3857.000)

\path(5573,3857)	(5647.539,3853.869)
	(5689.433,3853.477)
	(5723.000,3857.000)

\path(5723,3857)	(5778.080,3878.656)
	(5818.802,3900.488)
	(5873.000,3932.000)

\path(6323,3632)	(6252.102,3668.091)
	(6205.812,3702.778)
	(6180.617,3739.575)
	(6173.000,3782.000)

\path(6173,3782)	(6202.918,3829.509)
	(6248.000,3857.000)

\path(6248,3857)	(6280.651,3874.169)
	(6323.600,3891.755)
	(6398.000,3932.000)

\path(6398,3932)	(6450.226,3997.798)
	(6469.310,4039.671)
	(6473.000,4082.000)

\path(6473,4082)	(6452.243,4131.650)
	(6413.592,4174.828)
	(6367.146,4209.091)
	(6323.000,4232.000)

\path(6323,4232)	(6272.541,4246.951)
	(6211.730,4254.455)
	(6150.304,4250.731)
	(6098.000,4232.000)

\path(6098,4232)	(6047.420,4161.594)
	(6032.441,4118.292)
	(6023.000,4082.000)

\path(6023,4082)	(6026.015,4045.276)
	(6023.000,4007.000)

\path(6023,4007)	(5974.820,3965.716)
	(5932.373,3948.791)
	(5873.000,3932.000)

\put(5573,3632){\ellipse{150}{150}}
\put(6248,4082){\ellipse{150}{150}}
\path(6023,3557)	(6029.531,3612.126)
	(6031.707,3653.158)
	(6023.000,3707.000)

\path(6023,3707)	(5994.417,3752.296)
	(5948.000,3782.000)

\path(5948,3782)	(5905.671,3778.310)
	(5863.798,3759.226)
	(5798.000,3707.000)

\path(5798,3707)	(5752.618,3635.461)
	(5734.390,3592.870)
	(5723.000,3557.000)

\path(5723,3557)	(5726.885,3520.430)
	(5723.000,3482.000)

\path(5723,3482)	(5679.191,3439.674)
	(5621.919,3410.739)
	(5558.937,3398.684)
	(5498.000,3407.000)

\path(5498,3407)	(5439.228,3444.076)
	(5394.136,3502.194)
	(5363.477,3568.965)
	(5348.000,3632.000)

\path(5348,3632)	(5346.982,3688.832)
	(5358.226,3752.270)
	(5383.107,3811.823)
	(5423.000,3857.000)

\path(5423,3857)	(5458.655,3869.268)
	(5500.318,3867.905)
	(5573.000,3857.000)

\path(5573,3857)	(5647.539,3853.869)
	(5689.433,3853.477)
	(5723.000,3857.000)

\path(5723,3857)	(5778.080,3878.656)
	(5818.802,3900.488)
	(5873.000,3932.000)

\path(6323,3632)	(6252.102,3668.091)
	(6205.812,3702.778)
	(6180.617,3739.575)
	(6173.000,3782.000)

\path(6173,3782)	(6202.918,3829.509)
	(6248.000,3857.000)

\path(6248,3857)	(6280.651,3874.169)
	(6323.600,3891.755)
	(6398.000,3932.000)

\path(6398,3932)	(6450.226,3997.798)
	(6469.310,4039.671)
	(6473.000,4082.000)

\path(6473,4082)	(6452.243,4131.650)
	(6413.592,4174.828)
	(6367.146,4209.091)
	(6323.000,4232.000)

\path(6323,4232)	(6272.541,4246.951)
	(6211.730,4254.455)
	(6150.304,4250.731)
	(6098.000,4232.000)

\path(6098,4232)	(6047.420,4161.594)
	(6032.441,4118.292)
	(6023.000,4082.000)

\path(6023,4082)	(6026.015,4045.276)
	(6023.000,4007.000)

\path(6023,4007)	(5974.820,3965.716)
	(5932.373,3948.791)
	(5873.000,3932.000)

\put(7823,3557){\ellipse{150}{150}}
\put(7148,4007){\ellipse{150}{150}}
\path(7373,3482)	(7366.469,3537.126)
	(7364.292,3578.158)
	(7373.000,3632.000)

\path(7373,3632)	(7401.586,3677.296)
	(7448.000,3707.000)

\path(7448,3707)	(7490.329,3703.310)
	(7532.202,3684.226)
	(7598.000,3632.000)

\path(7598,3632)	(7643.382,3560.461)
	(7661.610,3517.870)
	(7673.000,3482.000)

\path(7673,3482)	(7669.111,3445.430)
	(7673.000,3407.000)

\path(7673,3407)	(7716.809,3364.674)
	(7774.081,3335.739)
	(7837.063,3323.684)
	(7898.000,3332.000)

\path(7898,3332)	(7956.771,3369.076)
	(8001.860,3427.194)
	(8032.519,3493.965)
	(8048.000,3557.000)

\path(8048,3557)	(8049.018,3613.832)
	(8037.774,3677.270)
	(8012.893,3736.823)
	(7973.000,3782.000)

\path(7973,3782)	(7937.345,3794.268)
	(7895.682,3792.905)
	(7823.000,3782.000)

\path(7823,3782)	(7748.461,3778.869)
	(7706.567,3778.477)
	(7673.000,3782.000)

\path(7673,3782)	(7617.920,3803.656)
	(7577.198,3825.488)
	(7523.000,3857.000)

\path(7073,3557)	(7143.898,3593.091)
	(7190.188,3627.778)
	(7215.383,3664.575)
	(7223.000,3707.000)

\path(7223,3707)	(7193.082,3754.509)
	(7148.000,3782.000)

\path(7148,3782)	(7115.349,3799.169)
	(7072.400,3816.755)
	(6998.000,3857.000)

\path(6998,3857)	(6945.774,3922.798)
	(6926.690,3964.671)
	(6923.000,4007.000)

\path(6923,4007)	(6943.757,4056.650)
	(6982.408,4099.828)
	(7028.854,4134.091)
	(7073.000,4157.000)

\path(7073,4157)	(7123.459,4171.951)
	(7184.270,4179.455)
	(7245.696,4175.731)
	(7298.000,4157.000)

\path(7298,4157)	(7348.580,4086.597)
	(7363.559,4043.296)
	(7373.000,4007.000)

\path(7373,4007)	(7369.985,3970.276)
	(7373.000,3932.000)

\path(7373,3932)	(7421.180,3890.716)
	(7463.627,3873.791)
	(7523.000,3857.000)

\put(5648,632){\ellipse{150}{150}}
\put(6323,182){\ellipse{150}{150}}
\path(6098,707)	(6104.531,651.868)
	(6106.707,610.835)
	(6098.000,557.000)

\path(6098,557)	(6069.417,511.704)
	(6023.000,482.000)

\path(6023,482)	(5980.671,485.690)
	(5938.798,504.774)
	(5873.000,557.000)

\path(5873,557)	(5827.618,628.542)
	(5809.390,671.134)
	(5798.000,707.000)

\path(5798,707)	(5801.885,743.574)
	(5798.000,782.000)

\path(5798,782)	(5754.191,824.326)
	(5696.919,853.261)
	(5633.937,865.316)
	(5573.000,857.000)

\path(5573,857)	(5514.228,819.924)
	(5469.136,761.806)
	(5438.477,695.035)
	(5423.000,632.000)

\path(5423,632)	(5421.982,575.168)
	(5433.226,511.730)
	(5458.107,452.177)
	(5498.000,407.000)

\path(5498,407)	(5533.655,394.732)
	(5575.318,396.095)
	(5648.000,407.000)

\path(5648,407)	(5722.539,410.131)
	(5764.433,410.523)
	(5798.000,407.000)

\path(5798,407)	(5853.080,385.336)
	(5893.802,363.506)
	(5948.000,332.000)

\path(6398,632)	(6327.102,595.909)
	(6280.812,561.223)
	(6255.617,524.425)
	(6248.000,482.000)

\path(6248,482)	(6277.918,434.491)
	(6323.000,407.000)

\path(6323,407)	(6355.651,389.831)
	(6398.600,372.245)
	(6473.000,332.000)

\path(6473,332)	(6525.226,266.202)
	(6544.310,224.329)
	(6548.000,182.000)

\path(6548,182)	(6527.243,132.347)
	(6488.592,89.173)
	(6442.146,54.912)
	(6398.000,32.000)

\path(6398,32)	(6347.541,17.049)
	(6286.730,9.545)
	(6225.304,13.269)
	(6173.000,32.000)

\path(6173,32)	(6122.420,102.406)
	(6107.441,145.708)
	(6098.000,182.000)

\path(6098,182)	(6101.015,218.724)
	(6098.000,257.000)

\path(6098,257)	(6049.820,298.284)
	(6007.373,315.209)
	(5948.000,332.000)

\put(7748,632){\ellipse{150}{150}}
\put(7073,182){\ellipse{150}{150}}
\path(7298,707)	(7291.469,651.868)
	(7289.292,610.835)
	(7298.000,557.000)

\path(7298,557)	(7326.586,511.704)
	(7373.000,482.000)

\path(7373,482)	(7415.329,485.690)
	(7457.202,504.774)
	(7523.000,557.000)

\path(7523,557)	(7568.382,628.542)
	(7586.610,671.134)
	(7598.000,707.000)

\path(7598,707)	(7594.111,743.574)
	(7598.000,782.000)

\path(7598,782)	(7641.809,824.326)
	(7699.081,853.261)
	(7762.063,865.316)
	(7823.000,857.000)

\path(7823,857)	(7881.771,819.924)
	(7926.860,761.806)
	(7957.519,695.035)
	(7973.000,632.000)

\path(7973,632)	(7974.018,575.168)
	(7962.774,511.730)
	(7937.893,452.177)
	(7898.000,407.000)

\path(7898,407)	(7862.345,394.732)
	(7820.682,396.095)
	(7748.000,407.000)

\path(7748,407)	(7673.461,410.131)
	(7631.567,410.523)
	(7598.000,407.000)

\path(7598,407)	(7542.920,385.336)
	(7502.198,363.506)
	(7448.000,332.000)

\path(6998,632)	(7068.898,595.909)
	(7115.188,561.223)
	(7140.383,524.425)
	(7148.000,482.000)

\path(7148,482)	(7118.082,434.491)
	(7073.000,407.000)

\path(7073,407)	(7040.349,389.831)
	(6997.400,372.245)
	(6923.000,332.000)

\path(6923,332)	(6870.774,266.202)
	(6851.690,224.329)
	(6848.000,182.000)

\path(6848,182)	(6868.757,132.347)
	(6907.408,89.173)
	(6953.854,54.912)
	(6998.000,32.000)

\path(6998,32)	(7048.459,17.049)
	(7109.270,9.545)
	(7170.696,13.269)
	(7223.000,32.000)

\path(7223,32)	(7273.580,102.406)
	(7288.559,145.708)
	(7298.000,182.000)

\path(7298,182)	(7294.985,218.724)
	(7298.000,257.000)

\path(7298,257)	(7346.180,298.284)
	(7388.627,315.209)
	(7448.000,332.000)

\put(6701.120,2132.000){\arc{3156.246}{2.0150}{4.2682}}
\put(6923.000,2132.000){\arc{3000.000}{2.2143}{4.0689}}
\put(6323.000,2132.000){\arc{3000.000}{5.3559}{7.2105}}
\put(6619.880,2132.000){\arc{3156.246}{5.1565}{7.4098}}
\put(6640.505,3019.505){\arc{1379.787}{4.2341}{5.3900}}
\put(6698.000,1194.500){\arc{1275.000}{1.0808}{2.0608}}
\put(6660.500,2694.500){\arc{1576.785}{4.2700}{5.1548}}
\put(6585.500,1194.500){\arc{855.132}{0.6610}{2.2318}}
\put(1523,2132){\ellipse{3030}{3030}}
\path(923,3557)(2123,782)
\path(2123,3557)(923,782)
\path(6023,3332)(6923,932)
\path(6323,3407)(7223,932)
\path(7223,3332)(6848,2282)
\path(6998,3407)(6698,2657)
\path(6623,1532)(6323,857)
\path(6473,1907)(6098,932)
\path(3623,2132)(4523,2132)
\path(4388.000,2094.500)(4523.000,2132.000)(4388.000,2169.500)
\end{picture} }
        	\hspace{-1.9mm}
	\end{array}
 $$
\caption{For the chord diagram $\Ga$ shown on the left, the 
construction of
a handelbody $\S_{\Ga}$ on the right.}\lbl{E}
\end{figure}

 Note that the 
resulting surface $\S_{\Ga}$ will have genus $g=5n+1$. In each 
$\S 
(v)$ we have three copies $a(v) ,b(v) ,c(v)$ of the curves 
$a,b,c$ in 
$\S_2$ and we can assume they avoid the holes where the tubes $\{ 
T(e)\}$ meet $\S (v)$. Now label the edges of $\Ga$ with labels 
$a,b$ or $c$, so that all the internal edges have label $b$ and 
the external edges are labeled alternately $a$ and $c$ as we go 
round the external circle. Thus the three edges incident to any 
vertex 
have all different labels. See figure \ref{F}. 
Now for each edge $e$ that connects the vertices $v_e$ and $v_e 
'$ we take a 
band sum of 
the curves $x(v_e )$ and $x(v_e ')$, where $x \in \{a,b,c \}$ is 
the label of 
$e$, using a band which travels along the boundary of the tube 
$T(e)$. We will denote this band-sum curve by $\hat e$. 
Note that there are $6n$ labeled curves $a(v) ,b(v) ,c(v)$ in 
$\S_{\Ga}$
which, after the above mentioned band sum, yield $3n$
curves $\hat e$ in $\S_{\Ga}$.

\begin{figure}[htpb]
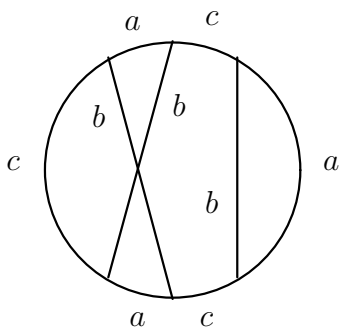

$$ \printname{F}
	\setlength{\unitlength}{0.03\standardunitlength}
	\begin{array}{c}  \hspace{-1.7mm}
        	\raisebox{-8pt}{\input draws/F.tex }
        	\hspace{-1.9mm}
	\end{array}
 $$
\caption{A labeling of the edges of a degree $3$ chord diagram 
with labels
$a,b,c$.}\lbl{F}
\end{figure}

Let $\gg_e$ be the diffeomorphism of $\S_{\Ga}$ defined by a Dehn 
twist along $\hat e$. Notice that $\gg_e\in\kg$ if $e$ has label 
$b$, since it can be arranged, by taking our connected sums 
correctly, 
that $\hat e$ bounds in $\S_{\Ga}$. See figure \ref{G}.

\begin{figure}[htpb]
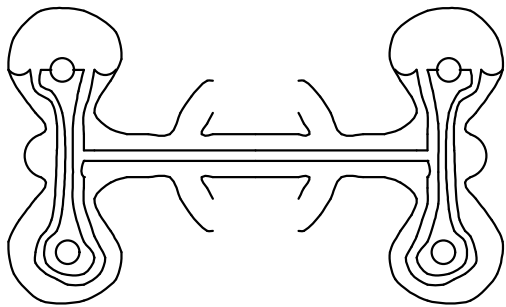

$$ \printname{G}
	\setlength{\unitlength}{0.03\standardunitlength}
	\begin{array}{c}  \hspace{-1.7mm}
        	\raisebox{-8pt}{\input draws/G.tex }
        	\hspace{-1.9mm}
	\end{array}
 $$
\caption{The connected sum of two $b$ labeled curves is a 
bounding
curve. Thus Dehn surgery along it represents an element of 
$\Kgg$.}\lbl{G}
\end{figure}
 
If we define our {\em Lagrangian} $L$ 
to be generated by the homology classes in each $\S (v)$ 
represented by 
$a(v)$ and $c(v)$ and, in addition, the meridians of the tubes 
$T(e)$, then $L$ is $f$-compatible  and each $\gg_e \in\Lg$.
It will
be a very important observation that whenever $\hat e_i$ and 
$\hat e_j$ are disjoint (for example if the two edges $e_i$ and
$e_j$ have the same label, or if $e_i$ has label 
$a$ 
and $e_j$ has label $c$) then  $\gg_{e_i}$ and $\gg_{e_j}$
commute. 
 .

\begin{itemize}
\item{{\bf Step 5}} 
\hspace{0.5cm}
Verification of claim \ref{claim.eL} for $n \neq 1$.
\end{itemize}

We begin with some preliminaries.
Choose any ordering $\I =\{ e_1 ,\cdots ,e_{3n}\}$ of the edges 
of 
$\Ga$. 
We now consider the element $\xi_{\I}^{\Ga}=(1-\gg_{e_1})\cdots 
(1-\gg_{e_{3n}})\in I^{3n}$, 
where $I=I\Lg$, and its image $\hfl (\xi_{\I}^{\Ga}) \in \GL 
{3n}$.
Recall from section \ref{sub.fti} that there is a map 
$\GL {3n} \to \Gas {3n}$ (induced by an inclusion map $\FL {3n} 
\subseteq
\Fas {3n}$); 
we will thus identify $\GL {3n}$
with its image in  $\Gas {3n}$. With this identification in mind,  
we will now 
describe the element 
$\hfl (\xi_{\I}^{\Ga})\in  \Gas {3n}$. Choose 
$3n$ concentric copies of $\S_{\Ga}$ in $\R^3$ and consider the 
curve $\hat e_i$ placed in the $i$-th copy of $\S_{\Ga}$ 
(counting from 
the outer to the inner). These curves define a $3n$ component 
link $K$ in
$S^3$ and we have that $[S^3 ,K] = \hfl (\xi_{\I}^{\Ga})$. If 
$e_i ,e_j ,e_k (i<j<k)$ are the edges incident to a vertex $v$, 
then 
the three curves $\hat e_i ,\hat e_j ,\hat e_k$ form a Borromean 
link 
if  $e_i ,e_j ,e_k$ have labels $a,b,c$, in that order, but form 
a 
{\em trivial link} if the labels are in any other order. We will 
call $v$ 
{\em $\I$-proper} in the former case and {\em $\I$-improper } in 
the latter 
case.
We call an ordering  $\I$, or the associated  element 
$\xi_{\I}^{\Ga}$  
proper if 
all vertices are $\I$-proper, otherwise improper.

Now recall the isomorphism $F ':\G_n\Ab\to\Gas {3n}$ of equation 
\eqref{eq.FF}
of section \ref{sub.fti}. 
Comparing the 
definition of this map with our description of $K$, we see that 
$\hfl (\xi_{\I}^{\Ga})=F '(\Ga_{\I})$
 where $\Ga_{\I}$ is the 
graph obtained from $\Ga$ by splitting every {\em improper} 
vertex $v$ into three univalent vertices. 
Observing  that $F '(\Ga_{\I}) \equiv F '(\Ga) \bmod F '(\Yb)$
(which follows from the fact that each trivalent graph with at 
least one
univalent vertex lies in $\Yb$, see the proof of lemma 
\ref{lemma.wb}),
it follows that 
\begin{equation}
\lbl{eq.coo}
\hfl (\xi_{\I}^{\Ga})=F '(\Ga_{\I}) \bmod F '(\Yb)
\end{equation}

For each label $x \in \{a,b,c \}$, let $x_i$ denote the edges 
with
label $x$, in any order. Fix an initial ordering  
$\I_0 = \{ a_1 ,\cdots ,a_n,b_1 ,\cdots 
,b_n,c_1 ,\cdots ,c_n \}$ of the edges of $\Ga$. 
To simplify the notation we will also write
$x_i$ when we really mean $\gg_{x_i}$. This should cause no
confusion. 

\begin{claim}
\lbl{claim.comm}
There is some formal commutator
$C$ in the $e_i$ such that 
$$ \xi_{\I_0}^{\Ga} \equiv 1-C+\text{ improper terms }\mod
 I^{3n+1}$$ 
\end{claim}

Note that equation \eqref{eq.coo} and the above claim 
imply that $ \hfl (\xi_{\I}^{\Ga})=F '(\Ga_b ) \bmod F '(\Yb ) $.
Using the definition of $F '$ (see equation \ref{eq.FF}) and the
definition of $\fl$ this 
imples that $F (\fl (\xi_{\I}^{\Ga}))=\Ga_b \bmod \Yb $ 
which finishes the proof of step 5.

\begin{pf}[of claim \ref{claim.comm}]
Choose any $a_r$ and $b_s$ which are incident. Then, using 
equation \eqref{eq.commute}, we have 
\begin{multline*}
1-\xi_{\I_0}^{\Ga}\equiv \pm\prod_{i\not= r}(1-a_i)(1-[a_r ,b_s 
])\prod_
{j\not= s}(1-b_j)\prod (1-c_i )\\
\pm\prod_{i\not= r}(1-a_i)(1-b_s)(1-a_r)\prod_
{j\not= s}(1-b_j )\prod (1-c_i )\mod I^{3n+1}
\end{multline*}
The second term on the right side is improper since $b_s$ 
precedes 
$a_r$. Next choose some $a_t$ which is incident to $b_s$. Then we 
have
\begin{multline*}
\prod_{i\not= r}(1-a_i)(1-[a_r ,b_s ])\prod_
{j\not= s}(1-b_j)\prod (1-c_i)\equiv \\ 
\shoveleft\prod_{i\not= r,t}(1-a_i)
(1-[a_t,[a_r ,b_s ]])\prod_{j\not= s}(1-b_j)\prod (1-c_i)\\-
\prod_{i\not= r,t}(1-a_i)(1-[a_r ,b_s ])(1-a_t)\prod_
{j\not= s}(1-b_j )\prod (1-c_i )\mod I^{3n+1} 
\end{multline*}
The second term on the right side can be further expanded by 
applying
equation \eqref{eq.commute} to $1-[a_r ,b_s ]$ and we obtain, 
$\mod I^{3n+1}$, a sum of two improper terms (since $b_s$ 
precedes
$a_t$ in both terms). We next look for some $b_u$ which is 
incident
to $a_t$ and, if it exists we can, in the same manner, write
\begin{multline*}
\prod_{i\not= r,t}(1-a_i)
(1-[a_t,[a_r ,b_s ]])\prod_{j\not= s}(1-b_j )\prod (1-c_i 
)\equiv\\
 \shoveleft
\prod_{i\not= r,t}(1-a_i )
(1-[[a_t,[a_r ,b_s ]],b_u ])\prod_{j\not= s,u}(1-b_j )\prod (1-
c_i )
\\ +\text{ improper terms }\mod I^{3n+1}\end{multline*}
Continuing in this way we eventually reach a point where we can 
write,
 after renumbering:
\begin{equation*}
\xi_{\I_0}^{\Ga}\equiv\pm \prod_{i\le p}(1-a_i)(1-C')\prod_{j\le 
q}(1-b_j )
\prod (1-c_i )+\text{ improper terms }\mod 
I^{3n+1}\end{equation*}
for some $p,q<n$ (actually $p=q$) and $C'$ is a formal commutator
in which each $a_i ,i>p$ and $b_j ,j>q$ appears once and no 
$a_i ,b_j$ with $i\le p,j\le q$ is incident to any $a_i ,b_j$ 
with
$i>p,j>q$. Thus $C'$ commutes with every $b_j$ and we have:
\begin{equation*}
\xi_{\I_0}^{\Ga}\equiv\pm\prod_{i\le p}(1-a_i)\prod_{j\le q}(1-
b_j )(1-C')
\prod (1-c_i )+\text{ improper terms }\mod 
I^{3n+1}\end{equation*}
We now play the same game with $c_i ,b_j$ and $C'$ that we just 
played with
$a_i ,b_j$. We will then eventually arrive at
\begin{equation}
\lbl{eq.prod}
\xi_{\I_0}^{\Ga}\equiv\pm\prod_{i\le p}(1-a_i )(1-C'')\prod_{j\le 
q}(1-b_j )
\prod_{i\le r}(1-c_i )+\text{ improper terms }\mod 
I^{3n+1}\end{equation}
for some new $q$ and $r$ and new commutator $C''$ involving all 
the $a_i ,b_i ,c_i$ 
not involved in the other terms on the right side, and no $b_i$ 
or
$c_i$ in $C''$ is incident to any not in $C''$. We now play the 
game
again with the $a_i ,b_i$ and $C''$. Going back and forth like 
this,
we eventually arrive at the point where we have an equation of 
the form
\eqref{eq.prod} where none of the edges in $C''$ are incident to 
any
of the edges outside $C''$. But, since $\Ga$ is connected, this 
is 
impossible unless there are no edges not in $C''$. 
\end{pf}

\begin{itemize}
\item{{\bf Step 6}} 
\hspace{0.5cm}
The Lagrangian $L$ and the \Heg\ surface $\S_{\Ga}$ do not depend 
on the 
choice of the chord diagram $\Ga$.
\end{itemize}

The equivalence (under isotopy) of the embeddings  
is not hard to see, assuming that we always choose the 
''natural'' 
embedding associated with each chord diagram, namely embed the 
external 
circle to coincide with the standard circle in $\R^3$ and embed 
the 
internal edges as straight lines between the endpoints on the 
external 
circle, introducing a small undercrossing or overcrossing when 
necessary 
to avoid intersecting another internal edge. See figure \ref{E} 
again. 
By sliding the ends of the 
tubes associated to the internal edges around on the torus 
associated to 
the external circle we can construct an isotopy between the 
embeddings 
associated to any two chord diagrams with the same number of 
edges. See 
figure \ref{H}. 

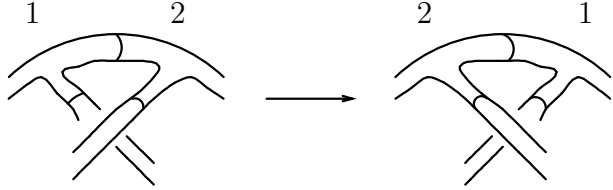
\begin{figure}[htpb]
$$ \printname{H}
	\setlength{\unitlength}{0.03\standardunitlength}
	\begin{array}{c}  \hspace{-1.7mm}
        	\raisebox{-8pt}{\begingroup\makeatletter\ifx\SetFigFont\undefined
% extract first six characters in \fmtname
\def\x#1#2#3#4#5#6#7\relax{\def\x{#1#2#3#4#5#6}}%
\expandafter\x\fmtname xxxxxx\relax \def\y{splain}%
\ifx\x\y   % LaTeX or SliTeX?
\gdef\SetFigFont#1#2#3{%
  \ifnum #1<17\tiny\else \ifnum #1<20\small\else
  \ifnum #1<24\normalsize\else \ifnum #1<29\large\else
  \ifnum #1<34\Large\else \ifnum #1<41\LARGE\else
     \huge\fi\fi\fi\fi\fi\fi
  \csname #3\endcsname}%
\else
\gdef\SetFigFont#1#2#3{\begingroup
  \count@#1\relax \ifnum 25<\count@\count@25\fi
  \def\x{\endgroup\@setsize\SetFigFont{#2pt}}%
  \expandafter\x
    \csname \romannumeral\the\count@ pt\expandafter\endcsname
    \csname @\romannumeral\the\count@ pt\endcsname
  \csname #3\endcsname}%
\fi
\fi\endgroup
\begin{picture}(8416,2405)(0,-10)
\thicklines
\put(1508.000,-67.000){\arc{4350.000}{3.9514}{5.4734}}
\put(6908.000,-67.000){\arc{4350.000}{3.9514}{5.4734}}
\put(1320.500,1920.500){\arc{530.330}{5.4978}{7.0686}}
\put(1770.500,1095.500){\arc{237.171}{4.3906}{6.6049}}
\put(1170.500,870.500){\arc{855.132}{3.8026}{4.4461}}
\put(6838.340,1899.065){\arc{440.473}{5.0342}{7.1374}}
\put(7358.000,1076.750){\arc{302.335}{4.1932}{6.4075}}
\put(6608.000,1133.000){\arc{212.132}{2.3562}{5.4978}}
\path(3608,1208)(4808,1208)
\path(4673.000,1170.500)(4808.000,1208.000)(4673.000,1245.500)
\path(8,1208)	(67.430,1254.955)
	(118.837,1295.211)
	(163.099,1329.429)
	(201.095,1358.266)
	(261.806,1402.438)
	(308.000,1433.000)

\path(308,1433)	(374.720,1481.446)
	(416.017,1502.892)
	(458.000,1508.000)

\path(458,1508)	(529.890,1472.918)
	(590.338,1408.381)
	(640.866,1337.404)
	(683.000,1283.000)

\path(683,1283)	(733.635,1234.747)
	(798.331,1174.021)
	(861.612,1111.535)
	(908.000,1058.000)

\path(908,1058)	(937.366,1007.210)
	(957.271,965.415)
	(983.000,908.000)

\path(3008,1208)	(2962.495,1238.332)
	(2923.275,1264.476)
	(2861.049,1305.957)
	(2816.048,1335.960)
	(2783.000,1358.000)

\path(2783,1358)	(2736.917,1394.420)
	(2679.800,1441.993)
	(2618.033,1485.069)
	(2558.000,1508.000)

\path(2558,1508)	(2501.640,1503.685)
	(2438.990,1483.546)
	(2379.595,1456.884)
	(2333.000,1433.000)

\path(2333,1433)	(2262.337,1388.512)
	(2220.908,1358.687)
	(2178.088,1326.283)
	(2135.881,1293.155)
	(2096.294,1261.162)
	(2033.000,1208.000)

\path(2033,1208)	(1964.744,1143.862)
	(1923.815,1103.104)
	(1881.080,1059.740)
	(1838.585,1016.158)
	(1798.376,974.748)
	(1733.000,908.000)

\path(1733,908)	(1688.437,863.437)
	(1633.857,808.857)
	(1572.598,747.598)
	(1508.000,683.000)
	(1443.402,618.402)
	(1382.143,557.143)
	(1327.563,502.563)
	(1283.000,458.000)

\path(1283,458)	(1216.904,391.904)
	(1176.065,351.065)
	(1133.000,308.000)
	(1089.935,264.935)
	(1049.096,224.096)
	(983.000,158.000)

\path(983,158)	(956.986,131.986)
	(908.000,83.000)

\path(908,458)	(953.501,503.501)
	(992.720,542.720)
	(1054.944,604.944)
	(1099.946,649.946)
	(1133.000,683.000)

\path(1133,683)	(1169.875,720.375)
	(1214.721,766.435)
	(1264.982,818.196)
	(1318.100,872.671)
	(1371.518,926.875)
	(1422.679,977.822)
	(1469.025,1022.525)
	(1508.000,1058.000)

\path(1508,1058)	(1573.137,1108.993)
	(1614.335,1139.090)
	(1658.000,1170.500)
	(1701.665,1201.910)
	(1742.863,1232.007)
	(1808.000,1283.000)

\path(1808,1283)	(1860.251,1329.789)
	(1924.767,1390.678)
	(1987.151,1453.977)
	(2033.000,1508.000)

\path(2033,1508)	(2087.814,1574.049)
	(2109.279,1615.504)
	(2108.000,1658.000)

\path(2108,1658)	(2066.781,1701.774)
	(2003.859,1721.761)
	(1936.757,1728.618)
	(1883.000,1733.000)

\path(1883,1733)	(1816.385,1737.767)
	(1775.525,1737.966)
	(1732.539,1737.238)
	(1689.610,1735.979)
	(1648.923,1734.589)
	(1583.000,1733.000)

\path(1583,1733)	(1545.926,1733.579)
	(1500.594,1734.984)
	(1449.734,1736.720)
	(1396.074,1738.291)
	(1342.342,1739.201)
	(1291.266,1738.953)
	(1245.576,1737.051)
	(1208.000,1733.000)

\path(1208,1733)	(1156.853,1719.917)
	(1094.101,1698.560)
	(1032.049,1675.673)
	(983.000,1658.000)

\path(983,1658)	(928.433,1648.458)
	(857.184,1639.122)
	(792.592,1620.476)
	(758.000,1583.000)

\path(758,1583)	(777.029,1517.273)
	(843.350,1456.141)
	(883.708,1428.106)
	(923.246,1402.189)
	(983.000,1358.000)

\path(983,1358)	(1016.048,1324.952)
	(1058.000,1283.000)
	(1099.952,1241.048)
	(1133.000,1208.000)

\path(1133,1208)	(1170.500,1170.500)
	(1208.000,1133.000)

\path(1208,1133)	(1234.014,1106.986)
	(1283.000,1058.000)

\path(1658,683)	(1703.501,637.499)
	(1742.720,598.280)
	(1804.944,536.056)
	(1849.946,491.054)
	(1883.000,458.000)

\path(1883,458)	(1935.043,405.958)
	(1976.524,364.476)
	(2033.000,308.000)

\path(1508,533)	(1564.476,476.524)
	(1605.957,435.042)
	(1658.000,383.000)

\path(1658,383)	(1691.048,349.952)
	(1733.000,308.000)
	(1774.952,266.048)
	(1808.000,233.000)

\path(1808,233)	(1841.048,199.952)
	(1883.000,158.000)
	(1924.952,116.048)
	(1958.000,83.000)

\path(1958,83)	(1984.014,56.986)
	(2033.000,8.000)

\path(8408,1208)	(8348.570,1254.955)
	(8297.163,1295.211)
	(8252.901,1329.429)
	(8214.905,1358.266)
	(8154.194,1402.438)
	(8108.000,1433.000)

\path(8108,1433)	(8041.280,1481.446)
	(7999.983,1502.892)
	(7958.000,1508.000)

\path(7958,1508)	(7886.108,1472.918)
	(7825.659,1408.381)
	(7775.130,1337.404)
	(7733.000,1283.000)

\path(7733,1283)	(7682.369,1234.747)
	(7617.672,1174.021)
	(7554.389,1111.535)
	(7508.000,1058.000)

\path(7508,1058)	(7478.634,1007.210)
	(7458.729,965.415)
	(7433.000,908.000)

\path(5408,1208)	(5453.505,1238.332)
	(5492.725,1264.476)
	(5554.951,1305.957)
	(5599.952,1335.960)
	(5633.000,1358.000)

\path(5633,1358)	(5679.083,1394.420)
	(5736.200,1441.993)
	(5797.967,1485.069)
	(5858.000,1508.000)

\path(5858,1508)	(5914.360,1503.685)
	(5977.010,1483.546)
	(6036.405,1456.884)
	(6083.000,1433.000)

\path(6083,1433)	(6153.662,1388.512)
	(6195.090,1358.687)
	(6237.909,1326.283)
	(6280.114,1293.155)
	(6319.701,1261.162)
	(6383.000,1208.000)

\path(6383,1208)	(6451.256,1143.862)
	(6492.185,1103.104)
	(6534.920,1059.740)
	(6577.415,1016.158)
	(6617.624,974.748)
	(6683.000,908.000)

\path(6683,908)	(6727.563,863.437)
	(6782.143,808.857)
	(6843.402,747.598)
	(6908.000,683.000)
	(6972.598,618.402)
	(7033.857,557.143)
	(7088.437,502.563)
	(7133.000,458.000)

\path(7133,458)	(7199.096,391.904)
	(7239.935,351.065)
	(7283.000,308.000)
	(7326.065,264.935)
	(7366.904,224.096)
	(7433.000,158.000)

\path(7433,158)	(7459.014,131.986)
	(7508.000,83.000)

\path(7508,458)	(7462.499,503.501)
	(7423.280,542.720)
	(7361.056,604.944)
	(7316.054,649.946)
	(7283.000,683.000)

\path(7283,683)	(7246.128,720.375)
	(7201.283,766.435)
	(7151.023,818.196)
	(7097.904,872.671)
	(7044.484,926.875)
	(6993.322,977.822)
	(6946.975,1022.525)
	(6908.000,1058.000)

\path(6908,1058)	(6842.863,1108.993)
	(6801.665,1139.090)
	(6758.000,1170.500)
	(6714.335,1201.910)
	(6673.137,1232.007)
	(6608.000,1283.000)

\path(6608,1283)	(6555.749,1329.789)
	(6491.233,1390.678)
	(6428.849,1453.977)
	(6383.000,1508.000)

\path(6383,1508)	(6328.186,1574.049)
	(6306.721,1615.504)
	(6308.000,1658.000)

\path(6308,1658)	(6349.219,1701.774)
	(6412.141,1721.761)
	(6479.243,1728.618)
	(6533.000,1733.000)

\path(6533,1733)	(6599.615,1737.767)
	(6640.475,1737.966)
	(6683.461,1737.238)
	(6726.390,1735.979)
	(6767.077,1734.589)
	(6833.000,1733.000)

\path(6833,1733)	(6870.074,1733.579)
	(6915.406,1734.984)
	(6966.266,1736.720)
	(7019.926,1738.291)
	(7073.658,1739.201)
	(7124.734,1738.953)
	(7170.424,1737.051)
	(7208.000,1733.000)

\path(7208,1733)	(7259.148,1719.917)
	(7321.903,1698.560)
	(7383.956,1675.673)
	(7433.000,1658.000)

\path(7433,1658)	(7487.562,1648.458)
	(7558.812,1639.122)
	(7623.406,1620.476)
	(7658.000,1583.000)

\path(7658,1583)	(7638.965,1517.273)
	(7572.642,1456.141)
	(7532.285,1428.106)
	(7492.749,1402.189)
	(7433.000,1358.000)

\path(7433,1358)	(7399.952,1324.952)
	(7358.000,1283.000)
	(7316.048,1241.048)
	(7283.000,1208.000)

\path(7283,1208)	(7245.500,1170.500)
	(7208.000,1133.000)

\path(7208,1133)	(7181.986,1106.986)
	(7133.000,1058.000)

\path(6758,683)	(6712.495,637.499)
	(6673.275,598.280)
	(6611.049,536.056)
	(6566.048,491.054)
	(6533.000,458.000)

\path(6533,458)	(6480.961,405.958)
	(6439.477,364.476)
	(6383.000,308.000)

\path(6908,533)	(6851.518,476.524)
	(6810.035,435.042)
	(6758.000,383.000)

\path(6758,383)	(6724.952,349.952)
	(6683.000,308.000)
	(6641.048,266.048)
	(6608.000,233.000)

\path(6608,233)	(6574.952,199.952)
	(6533.000,158.000)
	(6491.048,116.048)
	(6458.000,83.000)

\path(6458,83)	(6431.986,56.986)
	(6383.000,8.000)

\put(233,2258){\makebox(0,0)[lb]{$1$}}
\put(2258,2258){\makebox(0,0)[lb]{$2$}}
\put(7958,2258){\makebox(0,0)[lb]{$1$}}
\put(5708,2258){\makebox(0,0)[lb]{$2$}}
\end{picture} }
        	\hspace{-1.9mm}
	\end{array}
 $$
\caption{An isotopy of the handlebody on the left to the 
handlebody on the 
right via a handle slide.}\lbl{H}
\end{figure}

What happens to $L$ under such an isotopy? The generators of 
$L$ are either curves on the $\S (v)$ or meridians on the tubes 
associated to internal edges. In either case the isotopy 
preserves the 
curve and so there is no problem. This finishes the proof of 
theorem
\ref{thm.L}.
\end{pf}

We now discuss the analogous result for the Torelli group $\tg$ 
and 
Johnson's group $\kg$ instead of the Lagrangian group $\Lgg$. 
Recall  from section \ref{sub.fti} the maps $\phi^J_f$ for
$\J=\Tgg,\Kgg,\Lgg$, for our fixed \Heg\ splitting $f$, and the 
$f$-compatible Lagrangian $L$ of theorem \ref{thm.L}. For 
simplicity,
we drop the dependence on $f$ from the notation of the above 
maps.  
We now have the following corollary:

\begin{corollary}
\lbl{cor.L}
Assuming $f$ to be the standard genus $g$ \Heg\ splitting of 
$S^3$, 
the maps $\fk$ and $\ft$ are onto for $g\ge 5n+1$.
\end{corollary}

\begin{pf} 
Suppose $\Ga$ is a connected trivalent vertex oriented graph of 
degree $n$. 
We have 
constructed above an element $C\in (\Lg )_{3n}$ such that $\hfl 
(1-C)=\Ga \in \Gas {3n}$. 
Recall that $C$ is a formal commutator in the elements 
$a_i ,b_i ,c_i$ for $1\le i\le n$. Also recall that $b_i \in\kg$ 
and so, since $\kg$ is normal in $\Lg$, it follows that $C\in 
(\kg 
)_n$. Since  $\hfk$ (respectively, $\phi^L$) is induced by 
$\Phi^K$ 
(respectively, $\Phi^L$) it follows by lemma \ref{lemma.incl}
that $\phi^L(C)=\hfk(C)=\hfl(1-C) = \Ga$, which proves the 
corollary for $\kg$.

For $\tg$ we need to recall from section \ref{sub.johnson}
the fact that $\Kgg= \Tgg(2)$, where $\Tgg(2)$ is the second 
quotient
of the rational central series of $\Tgg$. Since 
$C\in (\kg )_n \sub\tg (2n)$, we have $C^r \in (\tg )_{2n}$ for 
some 
$r$. This and equation \eqref{eq.power} show that 
$1-C^r\equiv r(1-C)\mod (I\kg )^{n+1}$ and so 
$\ft (C^r )=\fk (C^r )=\Phi^K (1-C^r )=\Phi^K (r(1-C))=r\Ga$ and 
we are done.
\end{pf}

\section{Results for the subgroups $\Kgg$, $\Lgg$ of the 
mapping class group}
\lbl{sec.KL}

\subsection{Cocycles for $\Kgg$}
\lbl{subsec.K}

In this section we discuss some similar constructions of
 cohomology classes for the Johnson subgroup $\Kgg$, discussed in
section \ref{sub.johnson}. 
Following \cite{GL3}, recall  from  section \ref{sub.fti} that 
given an
 admissible surface $f:\S_g \hookrightarrow M$, 
the induced map $\Phi^K_f :\Q\Kgg \to\M$ maps the powers 
$(I\Kgg )^m$ into $\FK m$. Moreover, it was shown  that 
$\FK m \sub\FT {2m}=\Fas {3m} $ thus inducing a map of 
associated 
graded spaces $\GK m \to \GT {2m} \simeq\Gas {3m}
\simeq \G_m\A(\phi)$,  where the last isomorphism 
is given
by the fundamental theorem of \cite{LMO}, as was explained in 
the introduction. In remark \ref{rem.LKT} we pointed out that the 
composite map $\GK m \to \Gas {3m}$ is onto, a fact that we 
will not use here.
If we combine corollary \ref{cor.tt}, for $G=\Kgg , q=2$,
 with $\Phi^K_f ((I\Kgg )^m )$, we obtain cocycles 
$C^K_{f,m}\in 
C^m (\Kgg /\Kgg (2);\G_m\A(\phi))$. Unlike the case of the 
Torelli
 group, we {\em cannot} conclude that the pull-back of 
$[C^K_{f,m}]$ 
to 
$H^m (\Kgg ;\G_m \A (\phi ))$ is trivial.

The little that we can say is summarized in the following 
proposition.
\begin{proposition}
\lbl{prop.Kg}
\begin{itemize}
\item $C^K_{f,m}$ is multilinear and $\Gamma_g$-equivariant, 
i.e.,  satisfies the following
property (for $\alpha_i \in\Kgg /\Kgg (2), h \in \Ga_g$):
\begin{equation}
C_{h f,m}^K (h_\ast \alpha_1 ,  \ldots h_\ast \alpha_m )= 
C_{f,m}^K (\alpha_1, \ldots \alpha_m)
\end{equation}
\item If $f$ is an admissible Heegaard surface, then 
$C^K_{f,m}$ 
depends only 
on the "Lagrangians" $L^+ ,L^- 
\sub\pi /\pi (3)$, where $\pi=\pi_1 (\S_g )$ and 
$L^{\pm}=\text{Ker}\{ i_{\pm}\circ f_{\ast}:\pi /\pi (3)\to\pi_1 
(M_{\pm})/\pi_1 
(M_{\pm})(3)\}$.
\end{itemize}
\end{proposition}

\begin{remark}  
Let $\GKg$ denote the quotient group 
$\Gamma_g /\Kgg$. Then, from the work of Johnson, we have a 
short exact sequence:
$$ 1\to U\to\GKg\to Sp(H)\to 1  $$
$\GKg$ acts on $\Kgg /\Kgg (2)$ by conjugation, and 
on $\pi /\pi (3)$ by definition. It therefore acts on the 
Lagrangian pair in $\pi /\pi (3)$ and the equivariance 
property of $C^K_{f,m}$ is really an $\GKg$-equivariance. 
Moreover the subgroup of $\GKg$ which 
preserves the Lagrangian pair acts trivially.
\end{remark}

\begin{pf} 
The multilinearity and $\Gamma_g$-equivariance 
follows exactly as for the Torelli group.

The proof of the second assertion follows the same lines as 
the 
proof of
the analogous result in theorem \ref{thm.00} for the Torelli 
group. 
We need 
the following analogue of lemma \ref{lem.hndl}.

\begin{lemma} 
\lbl{lem.hndlK} 
Suppose that $Q$ is a handlebody. Set $\theta=\pi_1 (Q), 
\pi =\pi_1 (\partial Q)$ and $L=\text{Ker}\{ \pi /\pi (3)\to 
\theta 
/\theta(3)\}$.
 If $\a$ is an automorphism of $\pi /\pi (3)$ such that $\a 
(L)=L$, 
then there exists a diffeomorphism $h$ of $Q$ such that 
$(h|\partial Q)_{\ast}=\a$ (modulo inner automorphisms).
\end{lemma}
If we assume this lemma, then the rest of the proof proceeds 
as in 
the
 proof of theorem \ref{thm.00} with the following changes. 

The analogue of
lemma \ref{lem.c1} is established with $H_1 (\S )$ replaced 
by 
$\pi /\pi (3)$.
The proof only needs to be modified by observing that 
$g_{\e}$ 
belongs to 
$\Kgg$, since Johnson's result says that an element 
$g\in\Gamma_g$ 
belongs to $\Kgg$ if and only if it induces the identity on 
$\pi 
/\pi (3)$.

The analogue of lemma \ref{lem.c2}, with $H_1 (\S )$ replaced 
again 
by 
$\pi /\pi (3)$, is established by the same proof, with the 
extra 
observation 
that we may choose $h\in\Kgg$ using the result of Morita 
\cite{Mo2} 
that any 
\ihs\ is of the form $S^3_h$ for some $h\in\K_{g'}$, for some 
$g'$.

\begin{pf*}{Proof of lemma \ref{lem.hndlK}}Before we begin, 
lift 
$\a$ to an
automorphism $F/F(3)$, where $F=\pi_1 (\S_0 )$ is a free 
group. 
Recall 
the classical fact
 that the induced automorphism $\a_{\ast}$ on $H_1 (\S )$ is 
symplectic and 
so we can apply lemma \ref{lem.hndl} to find a diffeomorphism 
$g$ 
of 
$Q$ so
 that $(g|\partial Q)_{\ast}=\a_{\ast}$ on $H_1 (\S )$. Thus 
we 
can 
assume 
that $\a_{\ast}=\text{ identity}$. The effect of $\a$ on 
$F/F(3)$ 
is
measured by $\tau (\a )$, where $\tau$ is the Johnson 
homomorphism. 
We will
now use Morita's lemma \ref{lemma.mor}. We can choose a set 
of 
generators 
$\{ x_1 ,\ldots ,x_g ,y_1 ,\ldots ,y_g \}$ for $F$ 
so that $L$ is normally generated by $\{ y_1 ,\ldots ,y_g \} $. 
Then 
Morita's lemma
says that $\a$ is induced by a diffeomorphism of $\S_0$ if 
and 
only 
if 
$\tau (\a )$ belongs to the subgroup $W=\text{Ker}\{ \L3 H\to\L3 
H/\L3 
L\}$. Recall 
the definition of  $\tau (\a )$. For any $h\in F/F(3)$ we
can write $\a (h)=h\l (h)$. The assignment $h\to\l (h)$ 
defines a 
homomorphism
 $\l :H\to\Lambda^2 H\simeq F(2)/F(3)$. Now consider the 
element
\begin{equation}
\lbl{eq.tau}
 \tau(\a ) =\sum_i (x_i \otimes \l (y_i )-y_i \otimes\l (x_i 
))\in\L3 H\sub
 H\otimes\Lambda^2 H 
\end{equation}
The inclusion $\L3 H\sub H\otimes\Lambda^2 H $ is defined by 
$$a\wedge b\wedge c\to a\otimes (b\wedge c)+
b\otimes (c\wedge a)+c\otimes (a\wedge b) $$
Now it is easy to see that $W=W'\cap\L3 H$, where $W'\sub 
H\otimes\Lambda^2 H$
is the kernel of the projection $H\otimes\Lambda^2 H\to (H/L) 
\otimes\Lambda^2 (H/L)$. The condition that $\a (L)=L$ 
implies 
that 
$\l (h)\in\text{Ker}\{ \Lambda^2 H\to\Lambda^2 (H/L) \}$. 
Remembering 
that 
$L$
is generated by $\{ y_i \}$, we see that the first terms in 
equation 
\eqref{eq.tau} lie in $H\otimes\text{Ker}\{ \Lambda^2 
H\to\Lambda^2 (H/L) 
\}$ while 
the second terms lie in $L\otimes\Lambda^2 H$. Thus $\tau (\a 
)\in 
W$. 
\end{pf*} 
This completes the proof of proposition \ref{prop.Kg}.
\end{pf}

\subsection{Cocycles for $\Lg$}
\lbl{subsec.L}
Given an admissible surface $f:\S \hookrightarrow M$,
and an $f$-compatible Lagrangian $L$,  
recall from section \ref{sub.fti} (see also \cite{GL3}) the
Lagrangian  subgroup 
$\Lg\sub\Ga_g$ of the mapping class group and the 
associated  map $\Phi^L_f :\Lg\to\M$.
 It is proved in \cite{GL3} that $\Fas m$ is the 
union of the images $\Phi_f (I\Lg )^m )$ over all $f$ and $f$-
admissible
Lagrangians $L$. Recall 
also that $\Fas m$ is a $3$-step filtration and that $\Fas 
{3m}/\Fas {3m+1}\simeq\G_m\A(\phi)$. The general 
results of section \ref{sub.gen} yield the following. 

\begin{proposition}
\lbl{prop.Lg}
Let $f:\S_g\hookrightarrow M$ be an admissible surface,
$L$ a $f$-compatible Lagrangian  and $m$ be a nonegative integer. 
\begin{itemize}
\item There exist cocycles $C_{f,3m}^L \in C^{3m}(\Lg /\Lg 
(2) ;\G_m\A(\phi))$ for all such\ \ $m$, $C_{f,3m-1}^L \in 
C^{3m-1}(\Lg /\Lg (3) ;\G_m\A(\phi))$ for odd $m$ and   
$C_{f,3m-2}^L \in C^{3m-2}(\Lg /\Lg (4) ;\G_m\A(\phi))$ for 
even $m$.
\item The pullback of $C_{f,3m}^L$ to $C^{3m}(\Lg /\Lg (3) 
;\G_m\A(\phi))$ under the projection \newline $\Lg /\Lg (2)\to\Lg 
/\Lg 
(3)$ is a coboundary if $m$ is even. The pullback of 
$C_{f,3m-1}^L$ to  $C^{3m-1}(\Lg /\Lg (4) ;\G_m\A(\phi))$ 
under the projection $\Lg /\Lg (3)\to\Lg /\Lg (4)$ is a 
coboundary if $m$ is odd.
\end{itemize}
\end{proposition}
\begin{pf} The first statement follows from corollary 
\ref{cor.tt} and the paragraph preceding proposition 
\ref{prop.coh}. The second statement follows from proposition 
\ref{prop.coh}.
\end{pf}

\section{Discussion}
\lbl{sec.dis}

In this section we discuss  the results of the present paper 
in 
comparison with the work of R. Hain \cite{Hain} and S. Morita 
\cite{Mo5}, which  has been a source of motivation and 
inspiration for our results.

\subsection{Finite type invariants of knots and integral 
homology 3-spheres}
\lbl{sub.knots}

One of the results of the present paper is the construction 
of a cocycle $C_{f,m}: \otimes^{2m}U  \to \G_m\A(\phi)$ given an
admissible  surface $f: \S \hookrightarrow M$, see theorem 
\ref{thm.0}.
There is a well known (dictionary) correspondence
 between invariants of \ihs s and invariants of knots. 
For several statements  using the above
dictionary, see \cite{Hain}. We caution the reader however, 
that 
the 
above
mentioned dictionary is helpful in stating results, but {\em 
not 
necessarily}
in proving them.

In this section we discuss a related map after we replace 
\ihs s by knots and admissible surfaces by admissible braids.
%%Let $\sigma \in P_n$ be a pure braid on $n$ strands, and 
%%$\hat{\sigma} \subseteq S^3$
%%be a closure of it that represents a {\em knot};
Let $\sigma$ be a braid whose associated permutation
 is transitive, i.e. whose closure is a knot: 
such a braid will be called {\em admissible}. Let
$\A(S^1)$ be the vector space over $\BQ$ on the set
of admissible trivalent graphs with additional univalent 
vertices that
lie on a circle, divided out by the $AS$ and $IHX$ relations, 
see \cite{B-N}.
Using the  definition of the map of equation \eqref{eq.map}, 
and replacing
$\M$ (the vector space over $\BQ$ of \ihs s) by $\K$ (the 
vector space over $\BQ$ of oriented knots  in $S^3$), 
$\Fas \ast$ by the Vassiliev 
filtration $\F_\ast \K$, $\A(\phi)$ by $A(S^1)$, admissible 
surfaces 
by
admissible braids, and the fundamental theorem of \fti s of 
\ihs s 
by
the  fundamental theorem of \fti s of knots,  we can define a 
map: 
\begin{equation}
\lbl{eq.cknot}
C_{\sigma ,m}: \otimes^m (P_n/P_n(2)) \to \G_m\A(S^1)
\end{equation}
where $P_n$ is the pure braid group in $n$ strands. 
One can show that the above map coincides with the following 
one:
recall first that the abelianization $ P_n/P_n(2)$
of the pure braid group is a free abelian group in
generators $x_{ij}$ with $i,j=1,2, \ldots  , n$ and relations 
$x_{ii}=0$, $x_{ij}=x_{ji}$. 
Thus the tensor algebra  
$T(P_n/P_n(2))$ is a free (noncommutative) 
algebra
in generators $x_{ij}$ for $1 \leq i < j \leq n$. Monomials 
in 
this algebra
are represented in the left hand side of figure 
\ref{cdvertical} 
and
will be called chord diagrams on $n$ vertical strands.
The map $C_{\sigma,m} : \otimes^m (P_n/P_n(2)) \to 
\G_m\A(S^1)$
defined above coincides with the map that closes a degree $m$ 
monomial 
(thought of as a degree $m$ vertical chord diagram on $n$ 
strands)
to a chord diagram on $S^1$. The above mentioned closure of 
course 
depends
on the admissible braid $\sigma$ but only in a mild way: 
one can show that
it depends only on the image of the associated permutation.

As we discussed above, in the case of knots,  the map 
\eqref{eq.cknot} 
is well understood.
This is due to the fact that there is a nice presentation of 
the
pure braid group, and the fact that
the $I$-adic completion of the rational group
ring of the pure braid group $P_n$ is equal to a quotient of 
the 
tensor
algebra $T(P_n/P_n(2))$ modulo the ideal generated by the 
$4$-term
relation.

In the case of the Torelli group though, this is {\em not} 
the 
case.
To begin with, it is still unknown whether the Torelli group 
is
finitely presented. On the other hand,  the $I$-adic completion
of the  rational group ring of the Torelli group has been
recently  calculated
by R. Hain \cite{Hain}
using the {\em transcendental} theory of Mixed Hodge 
Structures. 
No
combinatorial proof of the result is known. The map 
$C_{\lpm,m}$ 
of 
equation \eqref{eq.map} may help us understand the structure 
of 
the 
Torelli
group in a combinatorial way, 
and, in the other direction, help us understand the space of
finite type invariants of
\ihs s. It may also be a first step in understanding Hain's
calculation.

We can now give the following dictionary
 between the case of knots and \ihs s, summarized
in the following table:
$$
\vbox{
\offinterlineskip
\halign{\strut#&&\vrule#&\quad\hfil#\hfil\quad\cr
\noalign{\hrule}
&&   && Knots && 3-Manifolds        &\cr
\noalign{\hrule}
&& $G$  && $P_n$ && $\Tgg$   &\cr
\noalign{\hrule}
&&  $G/G(2)$ && $Free(x_{ij})_{i < j}$  && $U$    &\cr
\noalign{\hrule}
&& admissible objects && braids && surfaces &\cr
\noalign{\hrule}
&& Graphical interpretation && Figure \ref{cdvertical}  && 
Figure
\ref{3wedge}   &\cr
\noalign{\hrule}
&& Cocycles && $C_{\sigma,m}$ && $C_{f,m}$ &\cr
\noalign{\hrule}
&& Chord diagrams  && $\A(S^1)$      && $\A(\phi)$    &\cr
\noalign{\hrule}
}}
$$

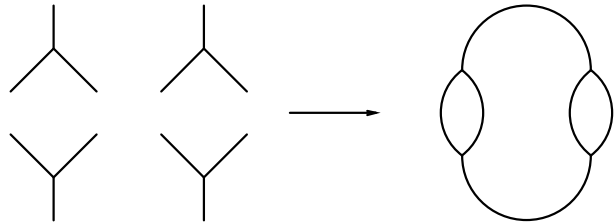
\begin{figure}[htpb]
$$ \printname{3wedge}
	\setlength{\unitlength}{0.03\standardunitlength}
	\begin{array}{c}  \hspace{-1.7mm}
        	\raisebox{-8pt}{\begingroup\makeatletter\ifx\SetFigFont\undefined
% extract first six characters in \fmtname
\def\x#1#2#3#4#5#6#7\relax{\def\x{#1#2#3#4#5#6}}%
\expandafter\x\fmtname xxxxxx\relax \def\y{splain}%
\ifx\x\y   % LaTeX or SliTeX?
\gdef\SetFigFont#1#2#3{%
  \ifnum #1<17\tiny\else \ifnum #1<20\small\else
  \ifnum #1<24\normalsize\else \ifnum #1<29\large\else
  \ifnum #1<34\Large\else \ifnum #1<41\LARGE\else
     \huge\fi\fi\fi\fi\fi\fi
  \csname #3\endcsname}%
\else
\gdef\SetFigFont#1#2#3{\begingroup
  \count@#1\relax \ifnum 25<\count@\count@25\fi
  \def\x{\endgroup\@setsize\SetFigFont{#2pt}}%
  \expandafter\x
    \csname \romannumeral\the\count@ pt\expandafter\endcsname
    \csname @\romannumeral\the\count@ pt\endcsname
  \csname #3\endcsname}%
\fi
\fi\endgroup
\begin{picture}(8421,3039)(0,-10)
\thicklines
\put(7212.000,2112.000){\arc{1800.000}{3.1416}{6.2832}}
\put(7212.000,912.000){\arc{1800.000}{6.2832}{9.4248}}
\put(6762.000,1512.000){\arc{1500.000}{2.2143}{4.0689}}
\put(5862.000,1512.000){\arc{1500.000}{5.3559}{7.2105}}
\put(7662.000,1512.000){\arc{1500.000}{5.3559}{7.2105}}
\put(8562.000,1512.000){\arc{1500.000}{2.2143}{4.0689}}
\path(612,3012)(612,2412)(12,1812)
\path(612,2412)(1212,1812)
\path(2712,3012)(2712,2412)(2112,1812)
\path(2712,2412)(3312,1812)
\path(12,1212)(612,612)(1212,1212)
\path(612,612)(612,12)
\path(2112,1212)(2712,612)(3312,1212)
\path(2712,612)(2712,12)
\path(3912,1512)(5112,1512)
\path(4992.000,1482.000)(5112.000,1512.000)(4992.000,1542.000)
\end{picture} }
        	\hspace{-1.9mm}
	\end{array}
 $$
\caption{On the left, $4$ trivalent vertices, and on the 
right
a particular closing to a trivalent graph.}\lbl{3wedge}
\end{figure}

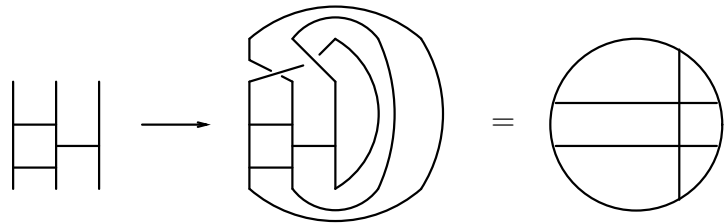
\begin{figure}[htpb]
$$ \printname{cdvertical}
	\setlength{\unitlength}{0.03\standardunitlength}
	\begin{array}{c}  \hspace{-1.7mm}
        	\raisebox{-8pt}{\begingroup\makeatletter\ifx\SetFigFont\undefined
% extract first six characters in \fmtname
\def\x#1#2#3#4#5#6#7\relax{\def\x{#1#2#3#4#5#6}}%
\expandafter\x\fmtname xxxxxx\relax \def\y{splain}%
\ifx\x\y   % LaTeX or SliTeX?
\gdef\SetFigFont#1#2#3{%
  \ifnum #1<17\tiny\else \ifnum #1<20\small\else
  \ifnum #1<24\normalsize\else \ifnum #1<29\large\else
  \ifnum #1<34\Large\else \ifnum #1<41\LARGE\else
     \huge\fi\fi\fi\fi\fi\fi
  \csname #3\endcsname}%
\else
\gdef\SetFigFont#1#2#3{\begingroup
  \count@#1\relax \ifnum 25<\count@\count@25\fi
  \def\x{\endgroup\@setsize\SetFigFont{#2pt}}%
  \expandafter\x
    \csname \romannumeral\the\count@ pt\expandafter\endcsname
    \csname @\romannumeral\the\count@ pt\endcsname
  \csname #3\endcsname}%
\fi
\fi\endgroup
\begin{picture}(9920,3033)(0,-10)
\thicklines
\put(3912.000,1509.000){\arc{2418.677}{5.2315}{7.3348}}
\put(4512.000,2109.000){\arc{1500.000}{3.7851}{5.6397}}
\put(4512.000,909.000){\arc{1500.000}{0.6435}{2.4981}}
\put(2824.500,1509.000){\arc{5033.947}{5.8529}{6.7135}}
\put(4512.000,1184.000){\arc{3650.000}{3.9948}{5.4299}}
\put(4512.000,1834.000){\arc{3650.000}{0.8533}{2.2883}}
\put(4062.000,1509.000){\arc{3911.521}{5.7165}{6.8499}}
\put(8712,1359){\ellipse{2400}{2400}}
\path(1812,1359)(2712,1359)
\path(2592.000,1329.000)(2712.000,1359.000)(2592.000,1389.000)
\path(12,1959)(12,459)
\path(612,1959)(612,459)
\path(1212,1959)(1212,459)
\path(12,1359)(612,1359)
\path(612,1059)(1212,1059)
\path(12,759)(612,759)
\path(3312,1959)(3312,459)
\path(3912,1959)(3912,459)
\path(4512,1959)(4512,459)
\path(3312,1359)(3912,1359)
\path(3912,1059)(4512,1059)
\path(3312,759)(3912,759)
\path(4512,1959)(3912,2559)
\path(3312,1959)(4062,2184)
\path(4287,2334)(4512,2559)
\path(3912,1959)(3762,2034)
\path(3612,2109)(3312,2259)(3312,2559)
\path(7587,1659)(9837,1659)
\path(7587,1059)(9837,1059)
\path(9312,2409)(9312,309)
\put(6687,1284){\makebox(0,0)[lb]{$=$}}
\end{picture} }
        	\hspace{-1.9mm}
	\end{array}
 $$
\caption{On the left chord diagrams on $3$ vertical strands, 
and on the right
the resulting a chord diagram on $S^1$ ontained by closing the 
chord diagram on the left.}\lbl{cdvertical}
\end{figure}

\subsection{Comparison with the results of Morita}
\lbl{sub.comp}

In this section we review briefly a very recent and important
paper \cite{Mo5} of Morita.
A common problem addressed in both Morita's recent paper and 
ours 
is
the one of constructing cocycles in various subgroups of the 
mapping 
class
group.
Morita \cite[theorem p.3]{Mo5}
uses a map\footnote{recall that we denote $U$ by $\L3 H/H$}
 $\rho_1: \Ga_{g} \to \frac{1}{2}U\rtimes 
Sp(H)$, to construct for each $i$ an $Sp(H)$-invariant 
element
 $\beta_i \in H^{2i}(\frac{1}{2}U, \BQ)$
with the property that, if $\bar{\beta_i}$ is the class in 
$H^{2i}(U\rtimes 
Sp(H);\Q )$ naturally associated to $\beta_i$, then  
$[\rho_1^{\ast}(\bar{\beta_i})]= e_i \in 
H^{2i}(\Ga_g, 
\BQ)$,
where $e_i$ is a certain cohomology class studied by Mumford 
and 
Morita.
By definition, the classes $\beta_i$ have the following 
properties:
\begin{itemize}
\item
They are defined on the {\em cocycle} level.
\item
The pullback $[\rho_1^{\ast}(\bar{\beta_i})]$ represents 
cohomology 
classes
in the full mapping class group.
\item
They are   $Sp(H)$ invariant cocycles.
\item
The pullback cocycles  in $\Kgg$ vanish.
\item
The coefficients of the cocycles  are rational numbers.
\end{itemize}
 
On the other hand, the cocycles $C_{\lpm,m}$ that we defined 
in 
theorems
\ref{thm.0} and \ref{thm.00}
have the following properties:
\begin{itemize}
\item
They are cocycles.
\item
They are defined in the abelianization of the Torelli group.
\item
The pullback of these cocycles to the Torelli group $\Tgg$
 represent trivial cohomology classes.
\item
They depend on a choice of Lagrangian pair $(L^+, L^-)$ and 
thus
are only $GL(L^+)$ invariant, and not $Sp(H)$ invariant.
\item
The pullback of these cocycles to $\Kgg$ vanish.
\item
The coefficients of these cocycles are the finite dimensional 
vector
spaces of manifold weight systems $\G_m\A(\phi)$.
\end{itemize}
 
\section{An epilogue or a beginning?}
\lbl{sec.epil}

We end this paper with the following  question.
 Recall from theorem \ref{thm.homo} the construction of a linear 
map
$D_m: (\G_{2m}\Tgg \otimes \BQ)^{Sp_g} \to \G_m\A^{conn}(\phi)$. 
This map
is stable with respect to the genus and for $m=\text{ even }$
it follows from  theorem \ref{thm.homo} that it is trivial. However,
for $m=1$ it 
was shown in theorem \ref{thm.homo}
to be a vector space isomorphism (of one dimensional vector 
spaces).
The authors now ask the following question:

\begin{question}
Is the map $D_m$ stably an isomorphism for $m=\text{ odd }$?
\end{question}

Note that a positive answer would connect several different areas 
together.
We will come back to the above question  in the near future.

%%%%%%%%%%%%%%%%%%%%\inplude{refs}

\ifx\undefined\bysame
	\newcommand{\bysame}{\leavevmode\hbox 
to3em{\hrulefill}\,}
\fi


\begin{thebibliography}{[EMSS]}


\bibitem[AM]{AM} S. Akbulut, J. C.  McCarthy, 
        {\em Casson's invariant for oriented homology 3-
spheres: 
        an exposition},  Princeton Math Notes, Princeton, 
1990.

\bibitem[B-N]{B-N} D. Bar-Natan, 
        {\em On the Vassiliev knot invariants}, Topology {\bf 
34} 
(1995)
        423-472.

\bibitem[CE]{CE} E. Cartan, Eilenberg,
        {\em Homological Algebra}, Princeton University 
Press, 
1956.

\bibitem[FH]{FH} W. Fulton, J. Harris,
        {\em Representation theory, a first course}, GTM {\bf 
129},
        Springer-Verlag, 1991.

\bibitem[Ga]{Ga1} S. Garoufalidis,
        {\em On finite type 3-manifold invariants I},  J. Knot
Theory and its Ramifications {\bf 5}, no. 4 (1996), p.441-462.
 
\bibitem[GL1]{GL1} S. Garoufalidis, J. Levine,
       {\em On finite type 3-manifold invariants II}, 
Brandeis 
Univ. and
       M.I.T. preprint June 1995, to appear in Math. Annalen.

\bibitem[GL2]{GL2} \bysame,
       {\em On finite type 3-manifold invariants IV: 
       comparison of definitions}, 
       Brandeis Univ. and M.I.T. preprint September 1995, to 
       appear in Proc. Camb. Phil. Soc.

\bibitem[GL3]{GL3} \bysame,
       {\em Finite type 3-manifold invariants, the mapping 
       class group and blinks}, Brandeis Univ. and
       M.I.T. preprint March 1996.

\bibitem[GO1]{GO1} S. Garoufalidis, T. Ohtsuki,
       {\em On finite type 3-manifold invariants III: 
       manifold weight systems}, 
       Tokyo Institute of Technology  and
       M.I.T. preprint August 1995.

\bibitem[Ha1]{Hain1} R. Hain,
       {\em Completions of the mapping class group and the cycle
       $C - C^{-}$},
       Mapping class groups and Moduli Spaces of Riemann surfaces,
       Contemp. Math {\bf 150} (1993) 75-105.
 
\bibitem[Ha2]{Hain} R. Hain,
       {\em Infinitesimal presentations of the Torelli groups},
       preprint December 1995.

%%%\bibitem[Jo1]{Jo1} D. Johnson, 
%%%       {\em Homeomorphisms of a surface which act 
%%%trivially on 
%%%homology},
%%%       Proc. Amer. Math. Soc., {\bf 75} (1979) 118-125.

\bibitem[Jo1]{Jo1} D. Johnson,
       {\em An abelian quotient of the mapping class group},
       Math. Ann. {\bf 249} (1980) 225-242.

\bibitem[Jo2]{Jo2} \bysame,
       {\em On the structure of the Torelli group III: the 
abelianization
       of} $\cal T$, Topology {\bf 24} (1085) 127-144.
 
%%%\bibitem[Jo2]{Jo2}\bysame ,
%%%       {\em A survey of the Torelli group},
%%%       Contemporary Math. {\bf 20} (1983) 163-179.

\bibitem[KM]{KM} N. Kawazumi, S. Morita,
       {\em The primary approximation to the cohomology of the 
moduli
       space of curves and stable characteristic classes}, 
preprint
       July 1996.
 
\bibitem[KK]{KK} K. Koike, I. Terada,
       {\em Young-diagramatic methods for the representation 
theory
       of the classical groups of type} $B_n, C_n, D_n$,
       Journal of Algebra, {\bf 107} (1987) 466-511.

\bibitem[Ko1]{Ko1} M. Kontsevich,
       {\em Formal (non)-commutative symplectic geometry},
       Gelfand Math. Seminars, 1990-92, Birkhauser, Boston, 
(1993) 
173-188.

\bibitem[Ko2]{Ko2} \bysame,
       {\em Feynmann diagrams and low-dimensional topology},
       Proceedings of the first European Congress of 
Mathematicians,
       vol. 2, Progress in Math. {\bf 120} Birkhauser, 
Boston, 
(1994)
       97-121.

\bibitem[LMO]{LMO} T.T.Q.  Le, J. Murakami, T. Ohtsuki,
   {\em A universal quantum invariant of 3-manifolds},
   preprint, November 1995.

\bibitem[L]{L} T.T.Q. Le, 
        {\em An invariant of \ihs s which is universal for 
all 
\fti s},
        preprint January 1996.

%%\bibitem[Lj]{Lj} ....

\bibitem[Mac]{Mac} S. MacLane, {\em Homology}, 
        Springer-Verlag, 1963. 

\bibitem[Mo1]{Mo0} S. Morita, 
        {\em Characteristic classes of surface bundles},
        Inventiones Math. {\bf 90} (1987) 551-577.

\bibitem[Mo2]{Mo1} \bysame, 
        {\em Casson's invariant for homology 3-spheres
         and characteristic classes of vector bundles I},
         Topology, {\bf 28} (1989) 305-323.

\bibitem[Mo3]{Mo2} \bysame, 
        {\em On the structure of the Torelli  group and the 
        Casson invariant},
         Topology, {\bf 30} (1991) 603-621.

\bibitem[Mo4]{Mo4} \bysame,
        {\em The structure of the mapping class group and 
        characteristic         classes of vector bundles},
        Contemporary Math. {\bf 150} (1993) 303-315.

\bibitem[Mo5]{Mojo} \bysame,
        {\em Abelian quotients of subgroups of the mapping class 
group
        of surfaces}, Duke Math. Journal, {\bf 70} (1993) 699-
726.

\bibitem[Mo6]{Mo5} \bysame,
        {\em A linear representation of the mapping class 
        group of orientable surfaces and characteristic 
        classes of vector bundles},
        preprint April, 1996.


%\bibitem[Mu]{Mu} D. Mumford,
%        {\em Towards an enumerative geometry of the moduli 
%        space of curves},
%        Arithmetic and Geometry, Progr. Math. {\bf 36} (1983) 
%        271-328.
      
\bibitem[Oh]{Oh} T. Ohtsuki,
        {\em Finite type invariants of integral homology 3-
        spheres}, 
        J. Knot Theory and its  Rami. {\bf 5} (1996) 101-115. 

%%\bibitem[Pe]{Pe} R.C. Penner,
%%        {\em The Poincare dual of the Weil-Petersson Kahler 
%%two 
%%form},
%%        preprint?......

\bibitem[Qu]{Qu} D. Quillen,
        {\em On the associated graded ring of a group ring},
        Journal of Algebra {\bf 10} (1968) 411-418.

\bibitem[W]{W} H. Weyl, {\em The classical groups }, Second 
Edition, 
Princeton U. Press 1946.

\bibitem[Wi1]{Wi1} E. Witten,
        {\em Two dimensional gravity and intersection theory 
on 
moduli space},
        Surveys in Diff. Geom. {\bf 1} (1991) 243-310.

\bibitem[Wi2]{Wi2} \bysame,
        {\em On the Kontsevich model and other models of two 
dimensional gravity},
        Proc. Conf. Diff. Geom. Methods in Physics, (S. Cato 
and 
A. Rocha Editors),
        Baruch College (1991).

\end{thebibliography}
\end{document}